\documentclass[aos,preprint]{imsart}

\RequirePackage[numbers]{natbib}
\RequirePackage[colorlinks,citecolor=blue,urlcolor=blue]{hyperref}

\usepackage{graphicx}
\usepackage{amsthm,amsmath,natbib}
\usepackage{amsfonts}
\usepackage{color}

\startlocaldefs

\newcommand{\be}{\begin{equation}}
\newcommand{\ee}{\end{equation}}
\newcommand{\bes}{\begin{equation*}}
\newcommand{\ees}{\end{equation*}}
\newcommand{\beqn}{\begin{eqnarray}}
\newcommand{\eeqn}{\end{eqnarray}}
\newcommand{\beqns}{\begin{eqnarray*}}
\newcommand{\eeqns}{\end{eqnarray*}}
\newcommand{\bea}{\begin{align}}
\newcommand{\eea}{\end{align}}
\newcommand{\beas}{\begin{align*}}
\newcommand{\eeas}{\end{align*}}

\newcommand{\lkr}{\left(}
\newcommand{\lkv}{\left[}
\newcommand{\rkv}{\right]}
\newcommand{\rkr}{\right)}
\newcommand{\lfi}{\left\{}
\newcommand{\rfi}{\right\}}

\newcommand{\fr}[1]{(\ref{#1})}

\newcommand{\del}{\delta}
\newcommand{\Del}{\Delta}

\newcommand{\al}{\alpha}
\newcommand{\af}{\alpha}

\newcommand{\ga}{\gamma}
\newcommand{\te}{\theta}
\newcommand{\om}{\omega}
\newcommand{\lam}{\lambda}
\newcommand{\Up}{\Upsilon}

\newcommand{\sig}{\sigma}

\newcommand{\Lam}{\Lambda}
\newcommand{\Om}{\Omega}
\newcommand{\Sig}{\Sigma}

\newcommand{\EE}{\ensuremath{{\mathbb E}}}

\newcommand{\II}{\ensuremath{{\mathbb I}}}

\newcommand{\PP}{\ensuremath{{\mathbb P}}}

\newcommand{\RR}{\ensuremath{{\mathbb R}}}

\newcommand{\vect}{\mbox{vec}}
\newcommand{\Pen}{\mbox{Pen}}

\newcommand{\Var}{\mbox{Var}}

\newcommand{\diag}{\mbox{diag}}
\newcommand{\supp}{\mbox{supp}}
\newcommand{\etal}{{\it et  al. }}

\newcommand{\Tr}{\mbox{Tr}}

\newcommand{\bal}{\mbox{\small bal}}

\newtheorem{theorem}{Theorem}
\newtheorem{lemma}{Lemma}
\newtheorem{corollary}{Corollary}

\newtheorem{remark}{Remark}

\newcommand{\ba}{\mathbf{a}}
\newcommand{\bb}{\mathbf{b}}

\newcommand{\bd}{\mathbf{d}}
\newcommand{\boe}{\mathbf{e}}
\newcommand{\bof}{\mathbf{f}}
\newcommand{\bg}{\mathbf{g}}

\newcommand{\bq}{\mathbf{q}}
\newcommand{\bt}{\mathbf{t}}

\newcommand{\bv}{\mathbf{v}}

\newcommand{\bx}{\mathbf{x}}
\newcommand{\by}{\mathbf{y}}
\newcommand{\bz}{\mathbf{z}}

\newcommand{\bA}{\mathbf{A}}
\newcommand{\bB}{\mathbf{B}}
\newcommand{\bC}{\mathbf{C}}
\newcommand{\bD}{\mathbf{D}}
\newcommand{\bF}{\mathbf{F}}
\newcommand{\bG}{\mathbf{G}}
\newcommand{\bI}{\mathbf{I}}
\newcommand{\bH}{\mathbf{H}}

\newcommand{\bQ}{\mathbf{Q}}
\newcommand{\bS}{\mathbf{S}}
\newcommand{\bU}{\mathbf{U}}
\newcommand{\bV}{\mathbf{V}}
\newcommand{\bW}{\mathbf{W}}
\newcommand{\bX}{\mathbf{X}}
\newcommand{\bY}{\mathbf{Y}}
\newcommand{\bZ}{\mathbf{Z}}

\newcommand{\bzero}{\mathbf{0}}
\newcommand{\bone}{\mathbf{1}}

 \newcommand{\blam}{\mbox{\mathversion{bold}$\lam$}}
 \newcommand{\bte}{\mbox{\mathversion{bold}$\te$}}
\newcommand{\bxi}{\mbox{\mathversion{bold}$\xi$}}
\newcommand{\boeta}{\mbox{\mathversion{bold}$\eta$}}

\newcommand{\bzeta}{\mbox{\mathversion{bold}$\zeta$}}

\newcommand{\bom}{\mbox{\mathversion{bold}$\om$}}

\newcommand{\bPhi}{\mbox{\mathversion{bold}$\Phi$}}
\newcommand{\bUp}{\mbox{\mathversion{bold}$\Up$}}

\newcommand{\bLam}{\mbox{\mathversion{bold}$\Lambda$}}
\newcommand{\bSig}{\mbox{\mathversion{bold}$\Sigma$}}
\newcommand{\bTe}{\mbox{\mathversion{bold}$\Theta$}}
\newcommand{\bXi}{\mbox{\mathversion{bold}$\Xi$}}
\newcommand{\bPi}{\mbox{\mathversion{bold}$\Pi$}}

\newcommand{\calB}{{\mathcal{B}}}
\newcommand{\calC}{{\mathcal{C}}}

\newcommand{\calF}{{\mathcal{F}}}
\newcommand{\calG}{{\mathcal G}}

\newcommand{\calK}{{\mathcal{K}}}

\newcommand{\calM}{{\mathcal M}}

\newcommand{\calP}{{\mathcal{P}}}
\newcommand{\calS}{{\mathcal{S}}}
\newcommand{\calT}{{\cal T}}

\newcommand{\calX}{{\cal{X}}}

\newcommand{\calZ}{{\cal{Z}}}

\newcommand{\hbLam}{\widehat{\bLam}} 
\newcommand{\hbd}{\widehat{\bd}}
\newcommand{\hbC}{\widehat{\bC}}
\newcommand{\hbS}{\widehat{\bS}}
\newcommand{\hbq}{\widehat{\bq}}
\newcommand{\hbte}{\widehat{\bte}}
\newcommand{\hbTe}{\widehat{\bTe}}
\newcommand{\hbW}{\widehat{\bW}}
\newcommand{\hbV}{\widehat{\bV}}
\newcommand{\hbZ}{\widehat{\bZ}}
\newcommand{\hbUp}{\widehat{\bUp}}
\newcommand{\hbUphJ}{\hbUp_{\hbC,\hJ}}

\newcommand{\tilbC}{\tilde{\bC}}

\newcommand{\tilbZ}{\tilde{\bZ}}

\newcommand{\hm}{\widehat{m}}
\newcommand{\hJ}{\widehat{J}}
\newcommand{\hM}{\widehat{M}}
\newcommand{\hrho}{\widehat{\rho}}

\newcommand{\bds}{{\bd^*}}
\newcommand{\bDs}{\bD^*}
\newcommand{\bQs}{{\bQ^*}}
\newcommand{\bCs}{{\bC^*}}
\newcommand{\bqs}{{\bq^*}}
\newcommand{\btes}{{\bte^*}}
\newcommand{\bLams}{{\bLam^*}}
\newcommand{\bWs}{{\bW^*}}
\newcommand{\bPhis}{{\bPhi^*}}
\newcommand{\bTes}{{\bTe^*}}
\newcommand{\bSs}{{\bS^*}}

\newcommand{\ms}{{m^*}}

\newcommand{\Ms}{{M^*}}%{{\ms(\ms+1)/2}}
\newcommand{\tilS}{\tilde{S}}

\newcommand{\PJ}{\bPi_{\bC,J}}

\newcommand{\hPhJ}{\widehat{\bPi}_{\hbC,\hJ}}

\newcommand{\PCJ}{\bPi_{\bC,J}}
\newcommand{\PCJs}{\bPi_{\bCs,J}}
\newcommand{\PCJso}{\bPi^{\bot}_{\bC^*,J}}
\newcommand{\PCJo}{\bPi_{\bC,J}^{\bot}}
\newcommand{\hPChJ}{\widehat{\bPi}_{\hbC,\hJ}}
\newcommand{\hPChJo}{\widehat{\bPi}_{\hbC,\hJ}^{\bot}}

\newcommand{\PSJ}{\bPi_{\bS,J}}

\newcommand{\hPShJ}{\widehat{\bPi}_{\hbS,\hJ}}

\newcommand{\rhons}{\rho_n^{*}}

\newcommand{\bl}{b_l}
\newcommand{\blo}{b_{l+1}}
\newcommand{\gl}{g_l}
\newcommand{\glo}{g_{l+1}}

\newcommand{\tCo}{\tilde{C}_0}
 
\long\def\ignore#1{}

\endlocaldefs

\begin{document}

\begin{frontmatter}

\title
{\bf      DYNAMIC NETWORK MODELS AND GRAPHON ESTIMATION}

\runtitle{DYNAMIC NETWORK MODELS AND GRAPHON ESTIMATION} 
\author{\fnms{Marianna Pensky}\thanksref{t1}\ead[label=e1]{Marianna.Pensky@ucf.edu}}
\runauthor{M. Pensky}
\thankstext{t1}{Supported in part by National Science Foundation (NSF),
grants  DMS-1407475 and DMS-1712977}

\affiliation{University of Central Florida}

\address
{Marianna Pensky\\
Department of Mathematics \\
University of Central Florida \\
Orlando FL 32816-1354, USA \\
\printead{e1}}

\begin{abstract}
In the present paper we consider a dynamic stochastic network model.
The objective is estimation of the tensor of connection probabilities $\bLam$
when it is generated by a Dynamic Stochastic Block Model (DSBM) or a dynamic graphon.
In particular, in the context of the DSBM, we derive a penalized least squares   estimator  $\hbLam$  of $\bLam$ 
and show that $\hbLam$ satisfies an oracle inequality and also attains minimax lower bounds for the risk.  
We extend  those results to estimation of $\bLam$ when it is generated by a dynamic graphon function.
The estimators constructed in the paper are adaptive to the unknown number of blocks in the context of the DSBM or
to the smoothness of the graphon function.  
The technique relies on the vectorization of the model and leads to 
much simpler mathematical arguments  than the ones used previously in the stationary set up.
In addition, all  results in the paper are non-asymptotic and allow a variety of  extensions.  
\end{abstract}

%%%%%%%%%%%%%%%%%%%%%%%%%%%%%%%%%%%%%%%%%%%%%%%%%%%%%%%%%%%%%%%%%%%%%%%%%%%%%%%%%%%%%%%%%%%%%%%%%%%%%%%%% 



\end{frontmatter}

%%%%%%%%%%%%%%%%%%%%%%%%%%%%%%%%%%%%%%%%%%%%%%%%%%%%%%%%%%%%%%%%%%%%%%%%%%%%%%%%%%%%%%%%%%%%%%%%%%%%%%%%% 
%%%%%%%%%%%%%%%%%%%%%%%%%%%%%%%%%%%%%%%%%%%%%%%%%%%%%%%%%%%%%%%%%%%%%%%%%%%%%%%%%%%%%%%%%%%%%%%%%%%%%%%%% 

\section{Introduction}
\label{sec:introduction}

Networks arise in many areas of research such as sociology, biology, genetics, ecology, information technology
to list a few. An overview of statistical modeling of random graphs can be found in, e.g.,   Kolaczyk (2009) and 
Goldenberg \etal (2011). While  static network models are relatively well understood, the literature on the dynamic 
network models is fairly recent.

In this paper, we consider a dynamic network defined as an undirected graph with $n$ nodes
with connection probabilities changing in time.
Assume that we observe the values of a tensor $\bB_{i,j,l} \in \{ 0,1\}$
at times  $t_l$ where $0 < t_1 < \cdots < t_L =T$. 
For simplicity, we assume that time instants are equispaced and the time interval is scaled to one, i.e. $t_l = l/L$.
Here $\bB_{i,j,l}=1$ if a connection between nodes $i$ and $j$ is observed at time $t_l$ and 
$\bB_{i,j,l}=0$  otherwise. 
We set $\bB_{i,i,l}=0$ and  $\bB_{i,j,l} = \bB_{j,i,l}$ for any $i,j = 1, \cdots n$ and $l=1, \cdots, L$,
and assume that $\bB_{i,j,l}$ are independent Bernoulli random variables with $\bLam_{i,j,l} = \PP(\bB_{i,j,l}=1)$
and $\bLam_{i,i,l}=0$.
Below, we study two types of objects:   a Dynamic Stochastic Block Model (DSBM) and a dynamic graphon.

The DSBM can be viewed as a natural extension of the Stochastic Block Model  (SBM) which, according to
Olhede and Wolfe  (2014), provides an  universal tool for description of time-independent stochastic network data.
In a DSBM,    all $n$ nodes are grouped into $m$ classes $\Om_1, \cdots, \Om_m$,
and probability of a connection $\bLam_{i,j,l}$ is entirely  determined by the  groups to which the nodes $i$ and $j$ 
belong at the moment $t_l$. In particular,  if $i \in \Om_{k}$ and $j \in \Om_{k'}$, then 
$\bLam_{i,j,l} = \bG_{k,k',l}$. Here, $\bG$ is the {\it connectivity tensor} at time $t_l$  with $\bG_{k,k',l} = \bG_{k',k,l}$.  
Denote by $n_k^{(l)}$ the number of nodes in class $k$ at the moment $t_l$, $k=1, \ldots, m$, $l=1, \ldots, L$.

 A dynamic graphon can be defined as follows.
Let $\bzeta = (\zeta_1, \cdots, \zeta_n)$ be a random vector sampled from a distribution $\PP_\zeta$
supported on $[0,1]^n$. Although the most common choice for $\PP_\zeta$ is the i.i.d. uniform distribution 
for each $\zeta_i$, we do not make this assumption in the present paper. 
We further assume that there exists a   function $f: [0,1]^3 \to [0,1]$ such that 
for any $t$ one has $f(x,y,t) = f(y,x,t)$ and  
% If we set 
\be \label{graphon}
\bLam_{i,j,l} = f(\zeta_i, \zeta_j, t_l), \quad i,j = 1, \cdots,n, \ l=1, \cdots, L.
\ee
Then, function $f$ summarizes behavior of the network and can be  called   {\it dynamic graphon},
similarly to the graphon  in the situation of a stationary network. This formulation allows   
to study a different set of  stochastic network models than the DSBM.

It is known that graphons play an important role in the theory of graph limits
described in Lov\'{a}sz and Szegedy (2006) and   Lov\'{a}sz (2012).
The definition of the dynamic graphon above fully agrees with their theory.
Indeed, for every   $l=1, \cdots, L$, the limit of  $\bLam_{*,*,l}$ as $n \to \infty$ 
is $f(\cdot,\cdot, t_l)$. We shall further elaborate on the notion of the dynamic graphon 
in Section~\ref{sec:dyn_graphon}.

\ignore{
Given an observed adjacency tensor $\bB$ sampled according to model \fr{graphon}, 
the graphon function $f$ is not identifiable since the topology of  a network 
is invariant with respect to any change of labeling of its nodes. Therefore, 
for any $f$ and any measure-preserving bijection $\mu: [0,1]\to [0,1]$
(with respect to Lebesgue measure), the functions $f(x,y,t)$ and $f(\mu(x), \mu(y),t)$
define the same probability distribution on random graphs. For this reason, we are 
considering equivalence classes of graphons. Note that in order for it to be possible to compare 
clustering of nodes across time instants, we introduce an assumption 
% (which, to the best of our knowledge, first appeared in Matias and Miele (2015)) 
that there are no label switching in time, that is, every  
node carries the same label at any time $t_l$, so that function $\mu$ is independent of $t$.

In addition, in our model we assume that, for every $x$ and $y$,
functions $f(x,y,\ldot)$ are smooth.   
}

In the last few years, dynamic network models attracted a great deal of attention
(see, e.g., Durante \etal (2015), Durante \etal   (2016), Han \etal (2015), 
Kolar \etal (2010), Leonardi \etal (2016), Matias and  Miele (2015),  Minhas \etal (2015),
Xing \etal (2010), Xu (2015), Xu and Hero III (2014) and Yang \etal  (2011) among others).
Majority of those paper describe   changes in the connection probabilities and 
group memberships via various  kinds of   Bayesian or Markov random field models 
and carry out the inference using the EM or iterative optimization algorithms.
While procedures described in those papers show good computational properties, 
they come without guarantees for the estimation precision. The only paper 
known to us that is concerned with estimation precision in the  dynamic  setting is by
Han \etal (2015) where the authors study consistency of their procedures 
when $n \to \infty$ or $L \to \infty$.

On the other hand,  recently,  several authors carried out  minimax studies
 in the context of stationary network  models. In particular, Gao \etal (2015) 
developed   upper and   minimax lower bounds for the risk of estimation of the   
matrix of connection probabilities. In a subsequent paper, Gao \etal (2016)
generalized the results to a somewhat more general problem of estimation of matrices with bi-clustering structures.
In addition, Klopp \etal (2017) extended these results to the case when the network is sparse in a sense 
that probability of connection is uniformly small and tends to zero as $n \to \infty$.
Also, Zhang and  Zhou (2016) investigated minimax rates  of community detection in the two-class 
stochastic block model.

The present paper has several objectives. First, we describe the non-parametric DSBM
 model  that allows  for smooth evolution of the tensor $\bG$ of connection probabilities  as well
as changes in group memberships in time.  Second,  we introduce   vectorization of the model
that enables us to take advantage of   well studied methodologies in nonparametric regression estimation. 
Using these techniques, we derive penalized least squares   estimators $\hbLam$ of $\bLam$ 
and show that they satisfy oracle inequalities. These inequalities do not require any assumptions on the mechanism
that drives evolution of the group memberships of the nodes in time and can be applied 
under very mild conditions. Furthermore, we consider a particular situation where only 
at most $n_0$ nodes can change their memberships between two consecutive time points.
Under the latter assumption, we derive minimax lower bounds 
for the risk of an  estimator of $\bLam$ and confirm that the estimators constructed in the paper  
attain those lower bounds.  Moreover, we extend those results to estimation 
of the   tensor $\bLam$ when it is generated by a graphon function.
We show that, for  the graphon, the estimators are minimax optimal within a logarithmic factor of $L$.
Estimators, constructed in the paper,  do not require   knowledge of the number of classes $m$ 
 in the context of the DSBM, or a degree of smoothness of the 
graphon function $f$ if $\bLam$ is generated by a dynamic graphon.

Note that unlike in Klopp \etal (2016) we do not consider a  network that is sparse 
in a sense that probabilities of connections between classes are  uniformly small. 
However, since our technique is based on model selection, it allows to study a network 
where some groups do not communicate with each other and obtain more accurate results.
Moreover, as we show in Section \ref{sec:sparse}, by adjusting the penalty, one can 
provide adaptation to uniform sparsity assumption if the number of nodes in each class is large enough.

The present paper makes several key contributions. 
First, to the best of our knowledge, the time-dependent networks are  usually handled via generative models that 
assume some probabilistic mechanism which governs the evolution  of the network  in time. The present paper offers 
the first fully non-parametric  model for the time-dependent networks which does not make any of such assumptions. 
It treats connection probabilities for each group as  functional data, allows group membership switching   
and enables one to exploit  stability in the group memberships over time.
Second, the paper provides the first minimax study of estimation of the tensor of connection probabilities in a dynamic setting.
The estimators constructed in the paper  are adaptive to the number of blocks in the context of the DSBM and 
to the smoothness of the graphon function in the case of a dynamic graphon. Moreover, the approach of the paper is  non-asymptotic, 
so it can be used irrespective of how large the number of nodes $n$, the number of groups $m$ and the number of time instants $L$ are and 
what the relationship between these parameters is.  
Third, in order to handle the tensor-variate functional data, we use vectorization of the model. 
This technique allows to reduce the problem of estimation of an unknown tensor of connection probabilities 
to a solution of a  functional linear regression problem with sub-gaussian errors. The technique is very potent and is used in a novel way. 
In particular, it leads to much more simple mathematics than in Gao \etal (2015)  and Klopp \etal (2017).
In the case of a time-independent SBM, it immediately reduces the SBM to a linear regression setting. 
In addition, by using the properties of the Kronecker product, we are able to reduce   
the  smoothness assumption on the  connection probabilities to sparsity assumption on 
their coefficients in one of the common  orthogonal transforms (e.g., Fourier or wavelet). 
Fourth, we use the novel structure of the penalty a part of which is proportional to the logarithm of
the  cardinality of the set of all possible clustering matrices over $L$ time instants.
The latter  allows to accommodate various group membership switching scenarios and is based on 
the Packing lemma (Lemma~4) which can be viewed as a version 
% the Packing lemma (Lemma~\ref{lem:packing}) which can be viewed as a version 
of the Varshamov-Gilbert  lemma for clustering matrices. In particular, while all papers that studied the SBM 
dealt with the situation where no restrictions are placed on the set of clustering matrices,
our approach allows to impose those restrictions.  
Finally,  the methodologies of the paper admit  various generalizations. 
For example, they can be adapted to a situation where the number of nodes 
in the network depends on time, or the connection probabilities have jump discontinuities,  
or when some of the groups have no connection with each other.
Section \ref{sec:sparse}  shows that the technique can be adapted to an additional uniform sparsity 
considered in Klopp \etal (2017)  if the number of nodes in each class is large enough.

The rest of the paper is organized as follows.
In Section~\ref{sec:notations_data}, we   introduce   the notations and  
describe   the vectorization of the model.
In Section~\ref{sec:assump_inference}, we  construct  the penalized least squares   estimators $\hbLam$
of the tensor $\bLam$. In Section~\ref{sec:oracle},  we derive the oracle inequalities for their risks.
In Section \ref{sec:DSBM_lower bounds}, we obtain the minimax lower bounds for the risk
that confirm that the estimators $\hbLam$ are minimax optimal. 
Section~\ref{sec:sparse} shows how our technique provides 
 adaptation to uniform sparsity assumption studied in Klopp \etal (2017).
Section~\ref{sec:dyn_graphon}    develops the nearly minimax optimal (within a logarithmic factor of $L$) 
estimators of $\bLam$ when the network is generated by a graphon.
Finally,   Section~\ref{sec:discussion}, provides a  discussion of  various generalizations of the techniques 
proposed in the paper. The proofs of all statements   are placed into   the Supplemental Material.

\section{Notation, discussion of the model   and data structures} 
\label{sec:notations_data}

\subsection{Notation} 
\label{sec:notations}

% For any positive integer $d$, denote $[d]= \{1,2, \cdots, d\}$.
% For any $a,b \in \RR$, denote $a \vee b = \max(a,b)$, $a \wedge b = \min(a,b)$.
For any two positive sequences $\{ a_n\}$ and $\{ b_n\}$, $a_n \asymp b_n$ means that 
there exists a constant $C>0$ independent of $n$ such that $C^{-1} a_n \leq b_n \leq C a_n$
for any $n$. For any set $\Om$, denote cardinality of $\Om$ by $|\Om|$.
For any $x$, $[x]$ is the largest integer no larger than $x$.

For any vector $\bt \in \RR^p$, denote  its $\ell_2$, $\ell_1$, $\ell_0$ and $\ell_\infty$ norms by, 
respectively,  $\| \bt\|$, $\| \bt\|_1$,  $\| \bt\|_0$ and $\| \bt\|_\infty$.
Denote by $\|\bt_1 - \bt_2\|_H$ the Hamming distance between vectors $\bt_1$ and $\bt_2$.
Denote by $\bone$ and $\bzero$ the vectors that have, respectively, only unit or zero elements.
Denote by $\boe_j$ the vector with 1 in the $j$-th position and all other elements equal to zero.

For a matrix $\bA$, its $i$-th row and $j$-th columns are denoted, respectively, by
$\bA_{i, *}$ and $\bA_{*, j}$. Similarly, for a tensor $\bA \in \RR^{n_1 \times n_2\times n_3}$,
we denote its $l$-th $(n_1 \times n_2)$-dimensional sub-matrix by $\bA_{*,*,l}$.
Let $\vect(\bA)$ be the vector obtained from matrix $\bA$ by sequentially stacking its columns. Denote  
by $\bA \otimes \bB$ the Kronecker product of matrices $\bA$ and $\bB$. Also, $\bI_k$ is the identity matrix of size $k$.
For any subset $J$ of indices, any vector $\bt$ and any matrix $\bA$, denote the restriction of $\bt$ to indices in $J$ by $\bt_J$
and the restriction of $\bA$ to columns $\bA_{*,j}$ with $j \in J$ by $\bA_J$. Also, denote by $\bt_{(J)}$  the modification of vector 
$\bt$ where all elements $\bt_j$ with $j \notin J$ are set to zero.

% Notation $\bA >0$ or $\bA \geq 0$ means, respectively,  that $\bA$ is positive or non-negative definite.
% Denote determinant of $\bA$ by $|\bA|$ and the largest, in absolute value, element of $\bA$ by $\| \bA\|_{\infty}$.
% Denote the Moore-Penrose inverse of matrix $\bA$ by $\bA^{+}$.
%

For any matrix $\bA$,  denote its spectral and Frobenius norms by, respectively,  $\| \bA \|_{op}$ and $\| \bA \|$.
Denote 
% Define the Hamming and the $\ell_0$ norms of a matrix $\bA$ by, respectively, 
  $\| \bA \|_H \equiv \| \vect(\bA) \|_H$, $\| \bA \|_\infty = \| \vect(\bA) \|_\infty$ and   $\| \bA \|_0 \equiv \| \vect(\bA) \|_0$.  
For  any tensor $\bA \in \RR^{n_1 \times n_2\times n_3}$, denote 
$\| \bA \|^2 = \sum_{l=1}^{n_3} \|\bA_{*,*,l}\|^2$.

Denote by $\calM (m,n)$ a collection of {\it membership} (or {\it clustering}) matrices $\bZ \in \{0,1\}^{n\times m} $, i.e. 
matrices  such that $\bZ$ has exactly one 1 per row
and $\bZ_{ik} =1$ iff  a node $i$ belongs to the class $\Om_k$ and is zero otherwise.
Denote by $\calC(m,n,L)$  a set of clustering matrices   such that 
\be  \label{clustmatr}
\calC(m,n,L) \subseteq \prod_{l=1}^L \calM (m,n).
\ee

\subsection{Discussion of the model} 
\label{sec:mod_discuss}

Note that   the values of  $\bB_{i,j,l}$ are independent given the values of $\bLam_{i,j,l}$,
that is,  $\bB_{i,j,l}$ are independent in the sense that their deviations from $\bLam_{i,j,l}$ are independent from each other.
Therefore, the values of  $\bB_{i,j,l}$ are linked to each other in the same way as observations 
of a continuous function with independent Gaussian  errors are related to each other.  
Moreover, in majority of papers treating dynamic block models (see, e.g., Durante \etal (2015), 
Han \etal (2015), Matias and Miele (2017), Yang \etal  (2011) among others), 
similarly  to the present paper, the authors assume that observations $\bB_{i,j,l}$ 
are independent given   $\bLam_{i,j,l}$. Note that this is not an artificial construct: 
Durante \etal (2015), for example, use the model for studying international relationships between countries over time.

The only difference between the present paper and the papers cited above is that we assume that the underlying 
connection probabilities $\bG_{*,*,l}$ are functionally linked (e.g., smooth) rather than being   
 probabilistically  related. Indeed,   many papers  that treat dynamic block models assume  
some Bayesian generative mechanism on the values of connection probabilities as well as on evolution 
of clustering matrices. In particular, they impose some prior distributions that relate  
 $\bG_{*,*,l+1}$ to $\bG_{*,*,l}$   and $\tilbZ^{(l+1)}$ to $\tilbZ^{(l)}$, the matrices of 
underlying probabilities and the clustering matrices for consecutive time points. 
Since the proposed generative mechanism may be invalid, we avoid making assumptions about the probabilistic
structures that generate connection probabilities and group memberships, and treat the network as a given object.
% we do not impose any  specific generative mechanism for the network and treat the network as a given object, 
% we avoid making such assumptions since they may not  be true. 
However, our model enforces, in a sense, a more close but yet flexible relation between the values of $\bB_{i,j,l}$ 
since $\bG_{*,*,l}$ are functionally (and not stochastically) related. 
Moreover, our theory allows to place any restrictions on the set of clustering matrices.

To illustrate this point, consider just one pair of nodes  $(i,j)$ and assume that these nodes do not switch their 
memberships between times   $t_l$ and $t_{l+1}$ and also that $\bG_{i,j,l}$ is continuous at $t_l$ . It is easy to see that if  
$\bG_{i,j,l}$ is close to zero (or one), then $\bG_{i,j,l+1}$ is also close to zero (or one) 
and, hence,  $\bB_{i,j,l}$ and $\bB_{i,j,l+1}$ are likely to be equal to zero (or one) simultaneously. 
This relationship takes place in general.

To simplify the narrative, just for this paragraph, 
denote $\bl = \bB_{i,j,l}$, $\blo = \bB_{i,j,l+1}$,  $\gl = \bG_{i,j,l}$  and $\glo = \bG_{i,j,l+1}$.
In order we are able to assert conditional probabilities $\PP (\bB_{i,j,l+1}=1|\bB_{i,j,l}=1) \equiv \PP(\blo =1|\bl=1)$
and $\PP (\bB_{i,j,l+1}=0|\bB_{i,j,l}=0) \equiv \PP(\blo =0|\bl=0)$, consider the situation 
where $\gl$ and $\glo$ are random variables with the joint pdf $p(\gl, \glo)$ 
such that, given $\gl$, on the average  $\glo$ is equal to $\gl$: $\EE(\glo|\gl) = \gl$.
Assume, as it is done in the present paper,  that, given $g_l$, values of $b_l$ are independent Bernoulli variables, so that 
$$
p(\bl, \blo|\gl, \glo)= \gl^{\bl}(1-\gl)^{1-\bl}\, \glo^{\blo}(1-\glo)^{1-\blo}.
$$
It is straightforward to calculate marginal probabilities $\PP(\bl=1) = \EE(\gl)$,
$\PP(\blo=1) = \EE(\glo) = \EE[\EE(\glo|\gl)] = \EE(\gl)$ and the joint probability
$\PP(\bl =1, \blo =1) = \EE(\glo \gl)= \EE[\EE(\glo \gl|\gl)] = \EE(\gl^2)$ which yields
\bes
\PP(\blo =1|  \bl  =1) - \PP(\blo=1)  = \frac{\EE(\gl^2)}{\EE(\gl)} - \EE(\gl) =  \frac{\Var(\gl)}{\EE(\gl)} >0
\ees
unless $\Var(\gl) =0$. The latter means that, even in the presence of the assumption of the conditional independence, 
the probability of interaction at the moment $t_{l+1}$ is larger if there were an interaction at the moment $t_l$
than it would be in the absence of this assumption. Similarly, repeating the calculation with $\gl$ and $\glo$
replaced by  $1-\gl$ and $1-\glo$, obtain 
\bes
\PP(\blo =0|  \bl  =0) - \PP(\blo=0) = \frac{\Var(\gl)}{\EE(1-\gl)} >0. 
\ees
In the absence of  the probabilistic assumptions on $\gl$ and $\glo$, we cannot evaluate 
those conditional probabilities but the relationship persists in this situation as well.

\subsection{Vectorization of the model} 
\label{sec:vectorization}

Note that tensor $\bLam$ of connection probabilities has a lot of structure. On one hand,  
it is easy to check that 
\be \label{eq:bLam}
\bLam_{*,*,l} = \tilbZ^{(l)} \bG_{*,*,l} (\tilbZ^{(l)})^T, \quad \bB_{i,j,l} \sim \mbox{Bernoulli} (\bLam_{i,j,l}),
\ee 
where $\tilbZ^{(l)} \in \calM (m,n)$ is the clustering matrix at the moment $t_l$. On the other hand,
for every $k_1$ and $k_2$, vectors $\bG_{k_1,k_2,*} \in \RR^L$ are comprised of values of some smooth functions
and, therefore, have low  complexity. Usually, efficient representations of such vectors are achieved by applying 
some orthogonal transform  $\bH$   (e.g., Fourier or wavelet transform), however, we cannot apply this transform to the original data 
tensor for two reasons. First, the errors in the model are not Gaussian, so application of $\bH$ will convert 
the data tensor with independent Bernoulli components into a data tensor with dependent entries that are not Bernoulli variables any more.
In addition, application of this transform to the original data will not achieve our goals since, although vectors
$\bG_{k_1,k_2,*}$ represent smooth functions, vectors $\bLam_{i,j,*}$ do not, due to possible switches  in the group memberships.
In addition, for every $l$, matrix  $\bLam_{*,*,l}$ in \fr{eq:bLam} forms the  so called bi-clustering structure 
(see, e.g., Gao \etal (2016)) which makes recovery of $\bG_{*,*,l}$ much harder than in the case of a  usual regression model.

In order  to handle all these intrinsic difficulties, we apply  operation of vectorization to $\bLam_{*,*,l}$. Denote
\be \label{vec_l}
\blam^{(l)} = \vect(\bLam_{*,*,l}),\quad \bb^{(l)} = \vect(\bB_{*,*,l}),\quad 
\bg^{(l)} = \vect(\bG_{*,*,l}).
\ee
Then, Theorem 1.2.22(i) of Gupta and Nagar (2000) yields 
{\small{
\be  \label{laml_full}
\blam^{(l)} = (\tilbZ^{(l)} \otimes \tilbZ^{(l)}) \bg^{(l)}, \ 
 \bb_i^{(l)} \sim {\rm Bernoulli} (\blam_i^{(l)}), \  i=1, \cdots, n^2,\ l=1, \cdots, L.
\ee 
}}
Note that  $\bb_i^{(l)}$ in \fr{laml_full} are independent for different values of $l$ but not $i$
due to the symmetry. In addition, the values of $\bb_i^{(l)}$ and $\blam_i^{(l)}$
that are corresponding to diagonal elements of matrices $\bB_{*,*,l}$ and   $\bLam_{*,*,l}$, 
are equal to zero by construction. Since all those values are not useful for estimation, 
we remove redundant entries from  vectors $\blam^{(l)}$ and $\bb^{(l)}$ for every $l=1, \cdots, L$. 
Specifically, in \fr{laml_full}, we remove the elements  in $\blam^{(l)}$ and 
the rows in $(\tilbZ^{(l)} \otimes \tilbZ^{(l)})$  corresponding, respectively,  
to   $\bLam_{{i_1},{i_2},l}$ and $(\tilbZ_{i_1,*}^{(l)} \otimes \tilbZ_{i_2,*}^{(l)})$ with $i_1 \geq i_2$.
We denote the reductions of vectors $\blam^{(l)}$, $\bb^{(l)}$ and matrices $(\tilbZ^{(l)} \otimes \tilbZ^{(l)})$
by, respectively, $\bte^{(l)}$, $\ba^{(l)}$ and $\tilbC^{(l)}$ obtaining 
{\small
\be \label{bern2}
\bte^{(l)} = \tilbC^{(l)} \bg^{(l)}, \quad \ba_i^{(l)} \sim {\rm Bernoulli} (\bte_i^{(l)}), \   
i=1, \cdots,   n(n-1)/2,\ l=1, \cdots, L.
\ee 
} 
Note that unlike in the case of $\bb^{(l)}$, elements $\ba_{i}^{(l)}$ and $\ba_{i'}^{(l')}$
are independent whenever $i \neq i'$ or $l \neq l'$. The interesting thing here is that matrices 
$\tilbC^{(l)}$ are still clustering matrices, i.e., $\tilbC^{(l)} \in \calM(n(n-1)/2,m^2)$. Indeed, 
$\tilbC^{(l)}$ are binary matrices such that, 
for  $i$ corresponding to $(i_1, i_2)$ with  $i_1 < i_2$ and $k$ corresponding to $(k_1, k_2)$ 
in $(\tilbZ_{i_1,k_1}^{(l)} \otimes \tilbZ_{i_2,k_2}^{(l)})$ one has 
$\tilbC^{(l)}_{i,k}=1$ if and only if the nodes $i_1 \in \Om_{k_1}$ and  $i_2 \in \Om_{k_2}$.

\begin{figure}
 % \centering
%\[ \hspace*{6mm} \includegraphics[height=6cm]{pic123.eps}  \hspace{3mm}\includegraphics[height=6cm]{pic5.eps}  \]
\[   \includegraphics[height=6.2cm]{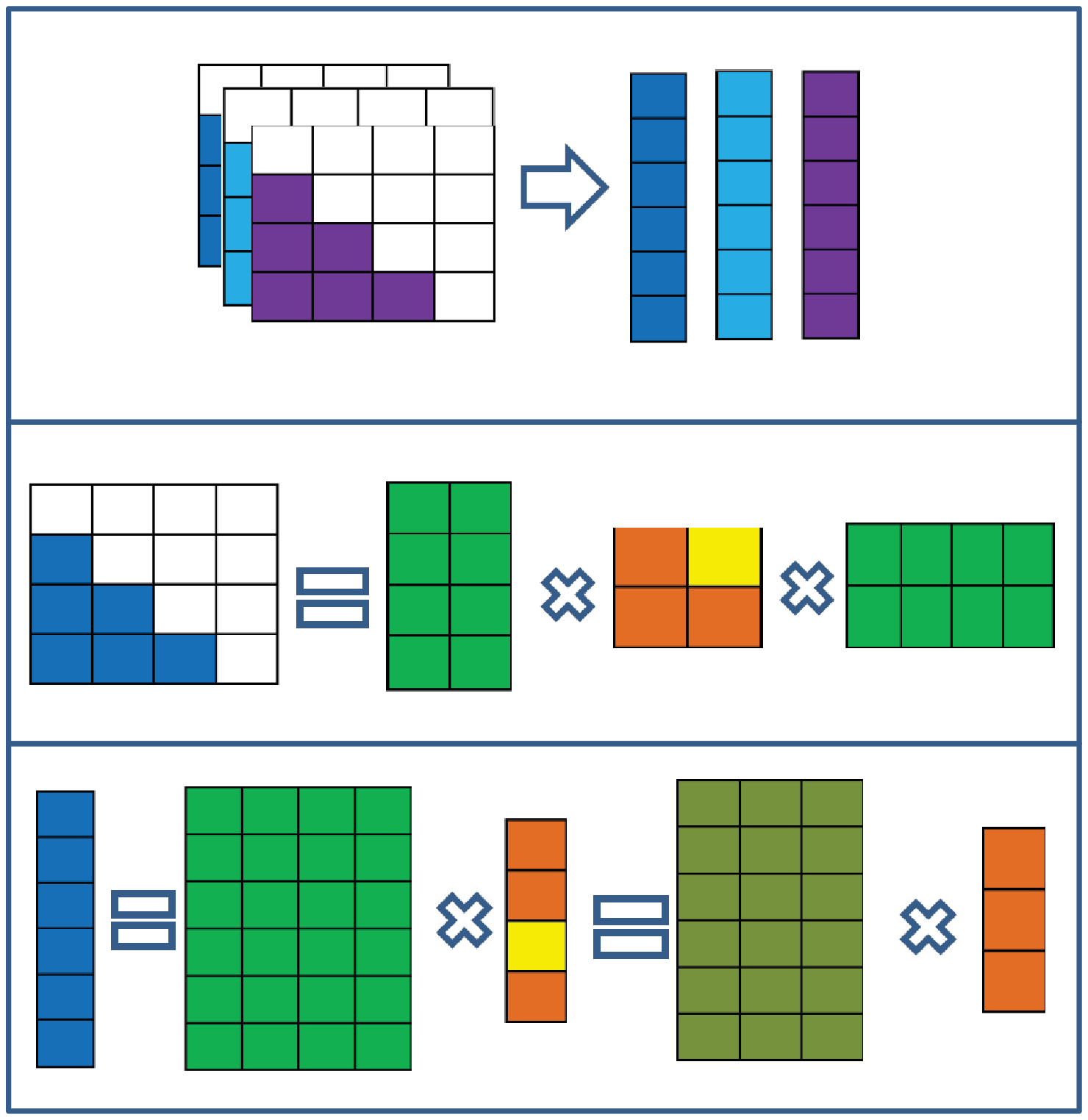}  \hspace{3mm}\includegraphics[height=6cm]{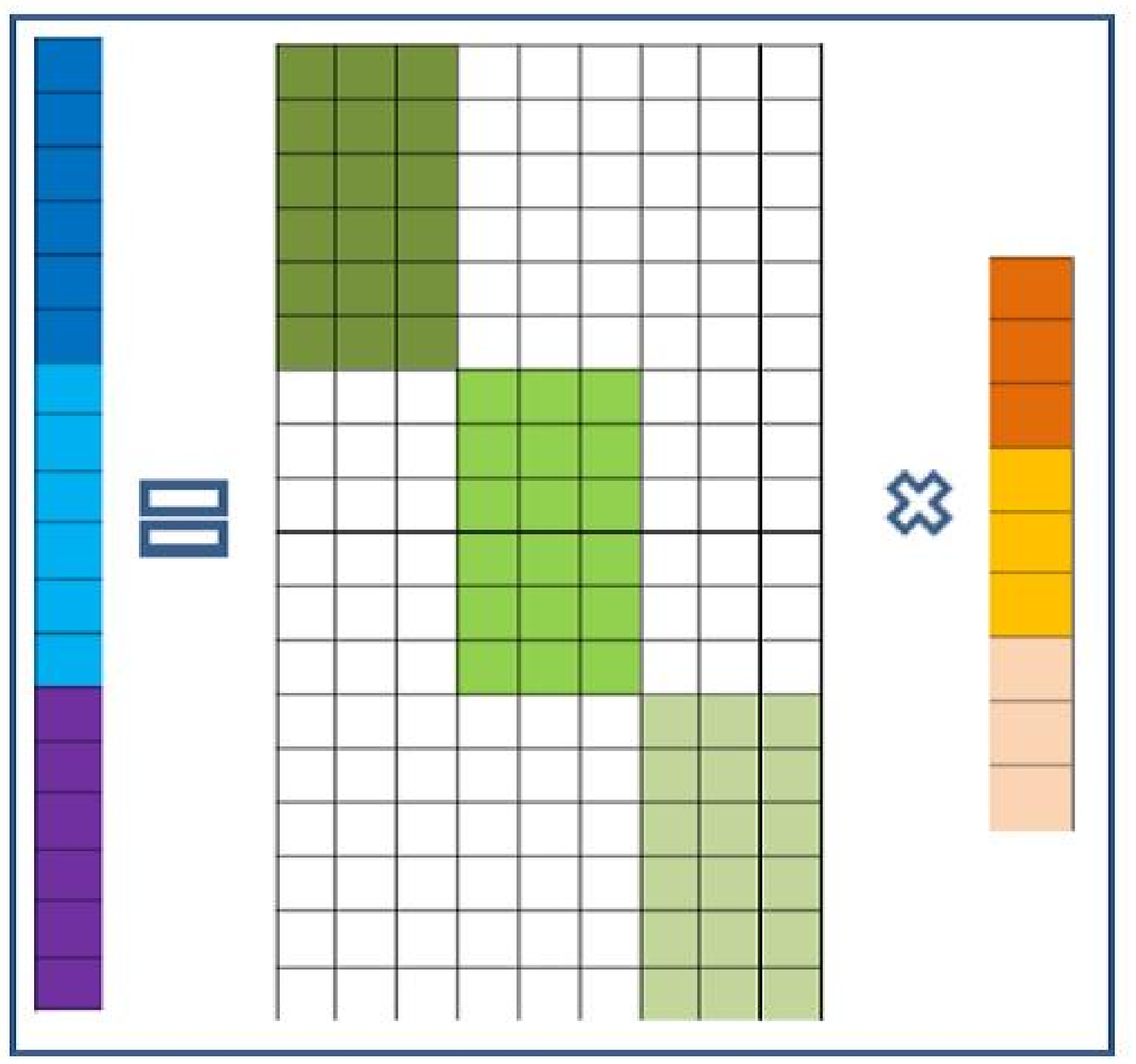}  \]
%\[   \includegraphics[height=6cm]{mfig1.eps}  \hspace{3mm}\includegraphics[height=6cm]{mfig2.eps}  \]
%
\caption{\small \label{figure2}   
 Vectorization of the probability tensor $\bLam$ with $n = 4$, $m=2$, $N=n(n-1)/2=6$, $M=m(m+1)/2=3$ and $L=3$.  
 Left panel, top: transforming $\bLam_{*,*,l}$ into $\bte^{(l)}$, $l=1,2,3$. 
 Left panel, middle: $\bLam_{*,*,1} = \tilbZ^{(1)} \bG_{*,*,1} (\tilbZ^{(1)})^T$.  
 Left panel, bottom: $\bte^{(1)} = \tilbC^{(1)} \bg^{(1)} =  \bC^{(1)} \bq^{(1)}$. 
In the left panel, redundant elements of $\bLam$ are white, redundant elements of $\bG$ are yellow. 
 Right panel: $\bte = \bC \bq$. }
\end{figure}

Observe that although we removed the redundant elements from vectors $\blam^{(l)}$ and $\bb^{(l)}$,
we have not done so for the vectors $\bg^{(l)}$. Indeed, since matrices $\bG_{*,*,l}$ are symmetric, 
the elements of vectors  $\bg^{(l)}$ corresponding to $\bG_{k_1,k_2,l}$ and $\bG_{k_2,k_1,l}$ with $k_1 \neq k_2$ are equal 
to each other. For the sake of eliminating such redundancy (and, hence, the need of tracing the equal elements 
in the process of estimation), for indices $k$ corresponding to pairs of classes $(k_1, k_2)$   with $k_1> k_2$, 
we remove entries $\bg^{(l)}_k$  from vectors $\bg^{(l)}$  and denote the resulting vectors by $\bq^{(l)}$. 
In order an equivalent of the relation \fr{bern2} still holds with vectors $\bq^{(l)}$ instead of $\bg^{(l)}$, 
we add together columns of matrices $\tilbC^{(l)}$ corresponding to $(k_1, k_2)$ and $(k_2, k_1)$ with $k_1 < k_2$, obtaining new matrices  $\bC^{(l)}$.
It is easy to see that, for every $l$,  since  $\bC^{(l)}$ is obtained from $\tilbC^{(l)}$ by adding columns together and 
since each row of $\tilbC^{(l)}$ has exactly one unit element with the rest of them being zeros, 
$\bC^{(l)}$ is  again a clustering matrix of size $[n(n-1)/2] \times [m(m+1)/2]$. In particular, for indices $i$ and $k$ corresponding to 
nodes $(i_1,i_2)$ and classes $(\Om_{k_1},\Om_{k_2})$ with $i_1 < i_2$ and $k_1 \leq k_2$, one has $\bC^{(l)}_{i,k}=1$
if  $i_1 \in \Om_{k_1}$ and $i_2 \in \Om_{k_2}$ or $i_1 \in \Om_{k_2}$ and $i_2 \in \Om_{k_1}$;
 $\bC^{(l)}_{i,k}=0$ otherwise. The process of vectorization of the model and removing redundancy is presented in Figure 1.

Using $\bC^{(l)}$ and $\bq^{(l)}$, one can rewrite equations \fr{bern2} as  
\begin{align} \label{model_l}
& \ba^{(l)} = \bte^{(l)} + \bxi^{(l)}  \quad \mbox{with} \quad \bte^{(l)} = \bC^{(l)} \bq^{(l)},  \quad l=1, \cdots, L,
% \label{param}
% &  \bte^{(l)} \in \RR^N, \ \bq^{(l)}\in \RR^M,\ \bC^{(l)}\in \calM(M,N).
%\  N=  n(n-1)/2,\ M =  m(m+1)/2   
\end{align}
where $\bC^{(l)}\in \calM(M,N)$,  $\bte^{(l)} \in \RR^N$, $\bq^{(l)}\in \RR^M$, $N=  n(n-1)/2$ and $M =  m(m+1)/2$. 
Here, for every $i$ and $l$, components $\ba_i^{(l)}$ of vector $\ba^{(l)}$ are independent Bernoulli variables with $\PP (\ba_i^{(l)}=1) = \bte_i^{(l)}$, 
so that components of vectors $\bxi^{(l)}$ are also independent for different values of $i$ or $l$.

If we had the time-independent SBM ($L=1$) and the clustering matrix were known, equation \fr{model_l} would reduce 
estimation of $\bq^{(1)}$ to the linear regression problem with independent sub-gaussian (Bernoulli) errors. 
Since in the case of the DSBM, for each $i$, the elements $\bg_i^{(l)}$, $l=1, \cdots, L$, of vector $\bg_i$ 
represent the values of a smooth function, we combine vectors in \fr{model_l} into matrices.
Specifically, we consider matrices $\bA, \bTe, \bXi \in  \RR^{N\times L}$ and $\bQ\in \RR^{M\times L}$
with columns $\ba^{(l)}$, $\bte^{(l)}$, $\bxi^{(l)}$ and $\bq^{(l)}$, respectively. Note that if the group memberships of the nodes were 
constant in time, so that $\bC^{(l)} \in \{0,1\}^{N \times M}$ were independent of $l$, formula \fr{model_l} would imply  
\be \label{con-memb}
\bA = \bTe + \bXi, \quad \bTe = \bZ \bQ 
\quad \mbox{if} \quad \bC^{(l)} = \bZ, \ l=1, \cdots, L.
\ee
However, we consider the situations where nodes can switch group memberships in time and \fr{con-memb} is not true.

For this reason, we proceed with further vectorization. We denote $\ba  = \vect(\bA)$, $\bte = \vect(\bTe)$ and  $\bq = \vect(\bQ)$
and observe that vectors $\ba, \bte \in \RR^{NL}$ and $\bq \in \RR^{ML}$ are obtained 
by stacking  vectors $\ba^{(l)}$, $\bte^{(l)}$ and $\bq^{(l)}$ in  \fr{model_l}   vertically for $l=1, \cdots, L$.  
Define a block diagonal matrix $\bC \in \{0,1\}^{NL \times ML}$ with blocks  $\bC^{(l)}$, $l=1, \cdots, L,$
on the diagonal. Then,   \fr{model_l} implies that
\be \label{full_model}
\ba  = \bte  + \bxi   \quad \mbox{with} \quad \bte  = \bC  \bq = \bC\, \vect(\bQ),
\ee
where $\ba_i$ are independent Bernoulli$(\bte_i)$ variables, $i = 1, \cdots,  NL$.

Observe  that if the  matrix $\bC$ were known, then 
equations in \fr{full_model} would represent a regression model with independent
 Bernoulli errors. Moreover, matrix $\bC^T \bC$ is diagonal 
since matrices $(\bC^{(l)})^T \bC^{(l)} = (\bS^{(l)})^2,\ l=1, \cdots, L,$ are diagonal with 
$\bS_{{k_1},{k_2}}^{(l)} = \sqrt{N^{(l)}_{{k_1},{k_2}}}$, where 
$N^{(l)}_{{k_1},{k_2}}$ is the number of pairs $(i_1, i_2)$ of nodes 
such that $i_1 < i_2$ and one node is in class $\Om_{k_1}$ while another is in class 
$\Om_{k_2}$ at time instant $t_l$: 
\be \label{eq:Nkl} 
N^{(l)}_{{k_1},{k_2}} = \lfi 
\begin{array}{ll}
n^{(l)}_{k_1} n^{(l)}_{k_2}, & \mbox{if}\ k_1 \neq k_2;\\
& \\
n^{(l)}_{k_1} (n^{(l)}_{k_1} -1), & \mbox{if}\ k_1 = k_2.
 \end{array} \right.
\ee

\begin{remark} \label{rem:direct}
{\bf (Directed graph). }
{\rm Similar vectorization algorithm can be used when the dynamic  network is constructed 
from directed graphs or  graphs with self loops. In the former case, the only redundant entries 
of matrices $\bLam_{*,*,l}$ would be the   diagonal ones while, in the latter case, $\bLam$ 
has no redundant elements and no row removal is necessary.
} 
\end{remark}

\begin{remark} \label{rem:biclust}
{\bf (Biclustering structures). }
{\rm 
Vectorization presented above can significantly simplify the inference in the so called biclustering
models  considered, for example, by Lee \etal (2010) and Gao \etal (2016). In those models, one needs to recover 
matrix $\bX$ from observations of matrix $\bY$ given by
$\bY = \bU_1 \bX \bU_2 + \bXi$ where matrices $\bU_1$ and $\bU_2$ are known and matrix $\bXi$ has independent 
zero-mean Gaussian or sub-gaussian entries. 
As long as there are no structural assumptions on matrix $\bX$ (such as, e.g., low rank),
one can apply vectorization and reduce the problem to the familiar non-parametric regression problem 
of the form $\by = \bU \bx + \bxi$ where matrix $\bU = \bU_1 \otimes \bU_2$ is known, 
$\bxi = \vect(\bXi)$ is the vector with independent components  
and one needs to recover $\bx = \vect(\bX)$ from observations $\by = \vect(\bY)$. 
}
\end{remark}

%%%%%%%%%%%%%%%%%%%%%%%%%%%%%%%%%%%%%%%%%%%%%%%%%%%%%%%%%%%%%%%%%%%%%%%%%%%%%%%%%%%%%%%%%%%%%%%%

\section{Assumptions and estimation for the DSBM}
\label{sec:assump_inference}

It is reasonable to assume that the values of the probabilities   $\bq^{(l)}$ 
of connections do not change dramatically from one  time instant   to another.
Specifically, we assume that  for various   $k = 1, \cdots, M$, 
vectors $\bq_{k} = (\bq_{k}^{(1)}, \cdots, \bq_{k}^{(L)})$  
represent   values of some smooth functions, so that $\bq_{k}^{(l)} = f_k (t_l)$, $l=1, \cdots, L$.
In order to quantify this phenomenon, we assume that vectors  $\bq_{k}$ have sparse
representation in some orthogonal basis $\bH \in \RR^{L \times L}$ with $\bH^T \bH = \bH \bH^T = \bI_L$,
so that vector  $\bH \bq_{k}^T$ is sparse: it has only few large coefficients, the rest of the coefficients
are small or equal to zero. This is a very common assumption in 
functional data analysis. For example,  if $\bH$ is the matrix of the Fourier transform and    $f_k$
belongs to a Sobolev space  or $\bH$ is a matrix of a wavelet transform and $f_k$ belongs to a Besov space, 
the coefficients $\bH \bq_{k}^T$ of $\bq_{k}^T$   decrease rapidly and, hence, vector $\bH \bq_{k}^T$
is sparse. In particular, one needs only few elements in vector $\bH \bq_{k}^T$ to represent $\bq_{k}$
with high degree of accuracy. The extreme case occurs when the connection probabilities do not change in time,
so that vector $\bq_{k}$ has constant components: then, for the Fourier or a periodic wavelet transform, 
the vector  $\bH \bq_{k}^T$  has only one non-zero element.

Denote $\bD= \bQ \bH^T$ where matrix $\bQ$ is defined in the previous section and $\bd = \vect(\bD)$.
Observe that vector $\bd$ is obtained by stacking together the columns of matrix $\bD = \bQ \bH^T$
while its transpose $\bD^T = \bH \bQ^T$ has vectors $\bH \bq_{k}^T$ as its columns.
% while vector  $\bd$ is obtained by stacking together the columns of its transpose $\bD = \bQ \bH^T$.  
%
Then, sparsity of the matrix $\bD$ can be controlled by imposing a complexity penalty on
$\|\bd \|_0 = \|\bD \|_0 = \|\bD^T\|_0$ on matrix $\bD$. Note that complexity penalty does  not require the actual 
matrix $\bD$ to have only few non-zero elements, it merely forces the procedure to keep only few 
large elements in $\bD$ while setting the rest of the   elements to zero and, hence,  acts as a kind of hard thresholding.
Note that by Theorem 1.2.22 of Gupta and Nagar (2000), one has
\be \label{vecQHT}
\bd = \vect(\bQ \bH^T ) = (\bH \otimes \bI_M) \vect(\bQ) = (\bH \otimes \bI_M) \bq = \bW \bq,
\ee 
where   $\bW = (\bH \otimes \bI_M)$   is an orthogonal matrix such
that $\bW^T \bW = \bW \bW^T = \bI_{ML}$. 
Denote  % coefficients of expansion of $\bq$ in $\bW$ by $\bd$ and the set of indices corersponding to nonzero components of $\bd$ by $J$:
\be \label{transformed}
% \bd = \bW \bq, \quad \bd \in \RR^{ML}, \ 
J \equiv J_M = \lfi j:\   \bd_j \neq 0 \rfi,\ \bd_{J^C} = \bzero,
%\quad \bZ = \bC \bW, %\quad \bd \in \RR^{ML}, 
\ee
so that $J$  is the set of indices corresponding to nonzero elements of the vector $\bd$.

Consider a set of clustering matrices $\calC(m,n,L)$ satisfying \fr{clustmatr}. 
At this point we   impose very mild assumption on $\calC(m,n,L)$: 
\be \label{clust_assump} 
\log(|\calC(m,n,L)|) \geq 2 \log m.
\ee
Assumption \fr{clust_assump} is used just for simplifying expression for the penalty.
Indeed, until now, we allowed any collection of clustering matrices, so potentially, we can 
work with the case where all cluster memberships are fixed in advance (although this would be a totally trivial case).
Condition \fr{clust_assump} merely means that at least two nodes at some point in time can be assigned arbitrarily to any of 
$m$ classes. 
Later,   we shall consider some special cases such as fixed membership (no membership switches over time)  or limited change 
(only at most $n_0$ nodes can change their memberships between two consecutive time points).
\ignore{
Observe that condition \fr{clust_assump} is extremely mild: if there are $n$ nodes that are grouped 
initially into $m$ classes, $n/m$ nodes per class and the group memberships do not change in time,
the inequality \fr{clust_assump} still holds since
\bes
\log(|\calC(m,n,L)|) \geq \log {n \choose n/m} \geq  \frac{n}{m} \log m \geq 2 \log m
\ees
provided $n \geq 2 m$.
}

We find $m, J, \bd$ and $\bC$  as one of the solutions of the following penalized least squares 
optimization problem
\be \label{opt_problem}
(\hm, \hJ,\hbd, \hbC) = \underset{m,J,\bd,\bC}{\operatorname{argmin}} \lkv \|\ba - \bC \bW^T \bd\|^2 + \Pen(|J|,m)\rkv
\quad \mbox{s.t.}\    \bd_{J^c}=\bzero  
%\arg\min \lfi m,J =J_M,\bd,\bC, \bd_{J^c}=0, \bC \in \calC(m,n,L): 
%\|\ba - \bC \bW^T \bd\|^2 + \Pen(|J|,m)\rfi,
\ee
where $\bC \in \calC(m,n,L)$, $\ba$ is defined in \fr{full_model},  $\bd \in \RR^{ML}$, $\bW \in \RR^{ML \times ML}$, $M=m(m+1)/2$ and 
\be \label{penalty}
\Pen(|J|,m) = 11 \,\log(|\calC(m,n,L)|) +   \frac{11}{2}\,|J| \log\lkr \frac{25\, m^2L}{|J|} \rkr.
\ee
Observe that the penalty in \fr{penalty} consists of two parts. The first part accounts for the complexity
of clustering and, therefore, allows one to obtain an estimator adaptive to the number of unknown groups $m$
as long as the we can express the complexity of clustering in terms of $m,n$ and $L$. The second term represents the price 
of estimating $|J|$ elements of vector $\bd$ and finding those $|J|$ elements in this vector of length $m (m+1) L/2$.

Note that since minimization is carried out also  with respect to  $m$,
optimization problem \fr{opt_problem}  should be solved separately for every $m =1, \cdots, n$,    
yielding $\hbd_M, \hbC_M$ and  $\hJ_M$. After that, one needs to select   
the value $\hM = \hm(\hm+1)/2$ that delivers the minimum in 
\fr{opt_problem}, so that
\be \label{est_val}
\hbd = \hbd_{\hM},\quad \hbC = \hbC_{\hM}, \quad \hJ = \hJ_{\hM}.
\ee 
Finally, due to \fr{transformed}, we  set  $\hbW = (\bH \otimes \bI_{\hM})$ and calculate 
\be \label{est}
 \hbq = \hbW^T \hbd,\quad  \hbte = \hbC \hbq.
\ee
We obtain $\hbLam$ by packing vector $\hbte$ into the  tensor and taking the symmetries into account.

%%%%%%%%%%%%%%%%%%%%%%%%%%%%%%%%%%%%%%%%%%%%%%%%%%%%%%%%%%%%%%%%%%%%%%%%%%%%%%%%%%%%%%%%%%%%%%%%

\section{Oracle inequalities for the DSBM}
\label{sec:oracle}

Denote the true value of tensor $\bLam$ by $\bLams$. Also, denote by $\ms$ the true number of groups, 
by $\bqs$ and  $\btes$ the true values of $\bq$   and $\bte$ in \fr{full_model}
and by  $\bCs$ the true value of $\bC$. Denote by $\bDs$ and $\bds$ the true 
values of matrix $\bD$ and vector $\bd$, respectively.
Let $\Ms = \ms(\ms+1)/2$ and $\bWs = (\bH \otimes \bI_{\Ms})$ 
be true values of $M$ and $\bW$.
Note that vector $\btes$ is obtained by vectorizing  $\bLams$ 
and then removing the redundant entries. Then, it follows from \fr{full_model} that 
\be \label{true_model}
\ba  = \btes  + \bxi   \quad \mbox{with} \quad \btes  = \bCs  \bqs = \bCs (\bWs)^T \bds.
\ee
Due  to the relation between the $\ell_2$ and the Frobenius norms, one has  
% for any vector $\bte$ and any tensor $\bLam$, one has 
\be \label{symrel}
\|\bte - \btes\|^2 \leq \|\bLam - \bLams\|^2 \leq 2 \|\bte - \btes\|^2,   
\ee 
and the following statement holds.

\begin{theorem} \label{th:oracle}
Consider a DSBM with a true matrix of probabilities $\bLams$ and the estimator $\hbLam$ obtained according to 
\fr{opt_problem}--\fr{est}. Let $\calC(m,n,L)$ be a set of clustering matrices satisfying conditions 
\fr{clustmatr} and \fr{clust_assump}.
Then, for any $t>0$, with probability at least $1 - 9 e^{-t}$, one has
{\small
\be \label{oracle_prob} 
\frac{\|\hbLam - \bLams\|^2}{n^2\,L} \leq \min_{\stackrel{m,J,\bd}{\bC \in \calC(m,n,L)}}
\lkv     \frac{6\, \| \bC \bW^T \bd_{(J)} - \btes \|^2}{n^2\,L} +  \frac{4\, \Pen(|J|,m)}{n^2\,L} \rkv + 
\frac{38\, t}{n^2\,L}  
\ee
}
 and 
{\small
\be \label{oracle_expec}
\EE  \lkr   \frac{\|\hbLam - \bLams\|^2}{n^2\,L}  \rkr \leq \min_{\stackrel{m,J,\bd}{\bC \in \calC(m,n,L)}}
\lkv  \frac{6\, \| \bC \bW^T \bd_{(J)} - \btes \|^2}{n^2\,L} +  \frac{4\, \Pen(|J|,m)}{n^2\,L}   + \frac{342}{n^2\,L}\rkv,
\ee}
where $\bd_{(J)}$ is the modification of vector 
$\bd$ where all elements $\bd_j$ with $j \notin J$ are set to zero.
\end{theorem}

The proof of   Theorem~\ref{th:oracle} is given in the Supplementary Material. % Section \ref{sec:suppl}.
Here, we just explain its idea. 
Note that if the values of $m$ and $C$ are fixed, the problem \fr{opt_problem}
reduces to a regression problem with a   complexity penalty $\Pen(|J|,m)$.
Moreover, if $J$ is known, the optimal estimator $\hbd$ of $\bd^*$ is just a projection estimator.
Indeed, denote $\bUp_{\bC} = \bC \bW^T$ and let $\bUp_{\bC,J} = (\bC \bW^T)_J$ be the reduction of matrix $\bC \bW^T$ to columns $j \in J$.
Given $\hm$, $\hJ$ and $\hbC$, one obtains $\hM = \hm(\hm+1)/2$,
$\hbW = (\bH \otimes \bI_{\hM})$,  $\bUp_{\bC,J} = (\bC \bW^T)_J$
%  $\hbUp = \hbC \hbW^T$ 
and  $\hbUp_{\hbC,\hJ} = (\hbC \hbW^T)_{\hJ}$.
Let  
\be \label{projection}
\PJ=   \bUp_{\bC,J} (\bUp_{\bC,J}^T \bUp_{\bC,J})^{-1} \bUp_{\bC,J}^T,\  % \quad \mbox{and} \quad
\hPhJ =  \hbUphJ (\hbUphJ^T \hbUphJ)^{-1} \hbUphJ^T
\ee 
be the projection matrices on the column spaces of $\bUp_{\bC,J}$ and $\hbUphJ$, respectively. 
Then, it is easy to see that $\hbUphJ\, \hbd = \hPhJ\, \ba$ and vector $\hbd$ is of the form 
\be \label{hbd_solution}
\hbd = (\hbUphJ^T \hbUphJ)^{-1}\, \hbUphJ^T\, \ba.
\ee
Hence, the values of $\hm$, $\hJ$ and $\hbC$ can be obtained as a solution of the following optimization problem
\bes 
(\hbC, \hm, \hJ) = \underset{m,J, \bC}{\operatorname{argmin}} \lkv \|\ba - \PJ \ba \|^2 + \Pen(|J|,m)\rkv
\quad \mbox{s.t.}\  % J \equiv J_M,\ %\bd_{J^c}=0,\ 
\bC \in \calC(m,n,L),
% \arg\min \lfi m,J =J_M, \bC,   \bC \in \calC(m,n,L): 
%\|\ba - \PhJ \ba \|^2 + \Pen(|J|,m)\rfi,
\ees 
where $\PJ$ and   $\Pen(|J|,m)$  are defined in \fr{projection} and \fr{penalty}, respectively. 
After that, we use the    arguments that are relatively standard in the proofs of oracle inequalities 
for the penalized least squares estimators.

Note that $\| \bC \bW^T \bd_{(J)} - \btes \|^2$ in  the right-hand sides of expressions \fr{oracle_prob} and \fr{oracle_expec},
 is the  bias term that quantifies how well one can estimate the true values of probabilities $\btes$ 
by   blocking them together, averaging the values in each block and simultaneously
setting all but $|J|$ elements of vector $\bd$ to zero. 
If $|J|$ is too small, then $\bd$ will not be well represented by its truncated version 
$\bd_{(J)}$ and the bias will be large. 
% coefficients of the expansions of the vector $\bq$ in \fr{full_model} in the basis $\bW$ to zero.
The penalty represents the stochastic error and constitutes the ''price'' for choosing too many blocks and coefficients.
In particular, the second term  $(11/2)\  |J| \log\lkr  25\, m^2L/|J|  \rkr$
% $(n^{2} L)^{-1}\, |J| \, \log(4 m^2L/|J|)$ 
in \fr{penalty} is due to the need of finding and estimating  $|J|$ 
elements of the $L m(m+1)/2$-dimensional vector. 
The  first term, $\log(|\calC(m,n,L)|)$, accounts for the difficulty of clustering
and is due to application of the union bound in probability.

Theorem \ref{th:oracle} holds for any collection  $\calC(m,n,L)$ of clustering matrices satisfying assumption \fr{clust_assump}. 
In order to obtain some specific results, denote by $\calZ (m,n,n_0,L)$ the collection of clustering matrices 
corresponding to the situation where   at most $n_0$ nodes can change their memberships between any two consecutive time points,
so that 
\be \label{card_clust}
|\calZ (m,n,n_0,L)| = m^n \lkv {n \choose n_0} m^{n_0} \rkv^{L-1}, 
% |\calZ (m,n,0,L)|= m^n; \quad |\calZ (m,n,n,L)| = m^{nL}.
\ee
yielding $|\calZ (m,n,0,L)|= m^n$ and $|\calZ (m,n,n,L)| = m^{nL}$.
Note that the case of $n_0=0$ corresponds to the scenario where  the group memberships of the nodes are 
constant and do  not depend on time while the case of $n_0=n$ means that memberships of all nodes can 
change arbitrarily from one time instant to another. 
% Then, one obtains the following corollary of Theorem \ref{th:oracle}.
Since 
\bes
\log \lkv {n \choose n_0} m^{n_0} \rkv \leq n_0  \,\log\lkr \frac{mne}{n_0}\rkr,  
\ees
formulae \fr{penalty} and  \fr{card_clust} immediately yield the following corollary.

\begin{corollary} \label{cor:upper_DSBM}
Consider a DSBM with a true matrix of probabilities $\bLams$ and estimator $\hbLam$ obtained according to 
\fr{opt_problem}--\fr{est} where $\calC(m,n,L) = \calZ (m,n,n_0,L)$. Then, inequalities \fr{oracle_prob} and \fr{oracle_expec} hold with 
\be \label{specific_pen}
\small{
\Pen(|J|,m) = 11 \lkv n \log m +   n_0 (L-1)\log\lkr \frac{mne}{n_0}\rkr + \frac{|J|}{2}  \log\lkr \frac{25\, m^2L}{|J|}\rkr \rkv
}
% \frac{\Pen(|J|,m)}{n^2\,L} = \frac{11\, \log m}{nL} + 
% \frac{11\, n_0 (L-1)}{n^2 L} \,\log\lkr \frac{mne}{n_0}\rkr + \frac{11\,  |J|}{2 n^2 L}\,  \log\lkr \frac{25\, m^2L}{|J|}\rkr  
\ee
\end{corollary}
 
It is easy to see that the first term in \fr{specific_pen} accounts for  the uncertainty of the initial clustering,
the second term is due to the changes in the group memberships of the nodes over time (indeed, if $n_0 =0$, this term just vanishes)
while the last term is identical to the second term in the expression for the generic penalty \fr{penalty}.
While we elaborate only on the special case where the collection of clustering matrices is given by \fr{card_clust},
one can easily produce results similar to Corollary \ref{cor:upper_DSBM} for virtually any nodes' memberships scenario.

\begin{remark} \label{rem:SBM}
{\bf (The SBM). }
{\rm Theorem \ref{th:oracle} provides an oracle inequality in the case of 
a   time-independent  SBM ($L=1$). Indeed, in this case, by taking $\bH = 1$ and $\bW = \bI_M$, obtain
for any $t>0$
\beqn  \label{oracle_probL1}
\frac{\EE  \|\hbLam - \bLams\|^2}{n^2} & \leq  & \min_{\stackrel{m,J,\bq}{\bC \in \calM(m,n)}}
\lkv     \frac{6  \| \bC \bq_{(J)} - \btes \|^2}{n^2} +  \frac{44 \, \log m}{n} \right. \\
& + &
\left. \frac{22 |J|}{n^2}  \log \lkr\frac{25\, m^2}{|J|}\rkr \rkv + \frac{342}{n^2}\nonumber  
\eeqn  
and a similar result holds for the probability.
Note that if $|J|=m(m+1)/2$, our result coincides with the one of Gao \etal (2015). However, if many groups 
have zero probability of connection, then   $|J|$ is small and the right-hand of \fr{oracle_probL1} 
can be asymptotically smaller than $n^{-1}\, \log m  + n^{-2} m^2$
obtained in  Gao \etal (2015). In addition, our oracle inequality is non-asymptotic 
and the  estimator is  naturally adaptive to the unknown number of classes.
(Gao \etal (2016) obtained adaptive estimators but not via an oracle inequality).
}
\end{remark}

Corollary \ref{cor:upper_DSBM} quantifies the stochastic error term in Theorem \ref{th:oracle}. 
The size of the bias depends on the level of sparsity of coefficients  of functions $\bq_k$ in the basis $\bH$
and on the constitution of classes. While one can study a variety of scenarios, in order to be specific,
we consider the case   of a {\it balanced network model} where the sizes of all the classes are proportional to each other, 
in particular, for some absolute constants $0 < \aleph_1 \leq 1 \leq  \aleph_2 < \infty$, one has
\be \label{balanced}
\aleph_1\, \frac{n}{m} \leq n_k^{(l)} \leq \aleph_2 \, \frac{n}{m}, \quad  k=1, \ldots, m,\ l=1, \ldots, L,
\ee
where $n_k^{(l)}$ the number of nodes in class $k$ at the moment $t_l$.

Note that the condition \fr{balanced} is very common in studying random network models
(see, e.g.,  Gao \etal (2017) or Amini and Levina  (2018)  among others). 
In addition, if class memberships are generated from the multinomial distribution
with the vector of probabilities $(\pi_1, \cdots, \pi_m)$,  
and $C_1/m \leq \pi_i \leq C_2/m$ for some constants $0<C_1<C_2<\infty$,
as it is done in, e.g., Bickel and Chen  (2009), condition   \fr{balanced} holds with high probability.

In particular, we consider networks that satisfy condition \fr{balanced} but yet allow only $n_0$ nodes 
switch their memberships between time instances. We denote the corresponding set of clustering matrices 
by $\calZ_{\bal} (m,n,n_0,L,\aleph_1, \aleph_2)$. It would seem that condition \fr{balanced} should make clustering much simpler. 
However, as  Lemma \ref{lem:balanced} below shows, this reduction does not makes estimation significantly easier since
the complexity  of the set of balanced clustering matrices $\log |\calZ_{\bal} (m,n,n_0,L,\aleph_1,\aleph_2)|$ 
is smaller than the complexity  of the set of unrestricted clustering matrices  $\log |\calZ (m,n,n_0,L)|$
only by, at most, a constant factor.

\begin{lemma} \label{lem:balanced} {\bf(Balanced network model complexity)}
% {\bf(Complexity of the balanced model).}
If $n \geq \sqrt{e\,n_0^3}$, then 
\be \label{card_balanced}
 % \frac{1}{4}  \lkr n \log m +  n_0  \,\log\lkr \frac{mne}{n_0}\rkr \rkr \leq 
%\log |\calZ_{\bal} (m,n,n_0,L,\aleph_1,\aleph_2)| \leq n \log m +  (L-1)\, n_0  \,\log\lkr \frac{mne}{n_0}\rkr.
\log |\calZ_{\bal} (m,n,n_0,L,\aleph_1,\aleph_2)| \geq \frac{1}{4} \lkv n \log m +  (L-1)\, n_0  \,\log\lkr \frac{mne}{n_0}\rkr \rkv.
\ee 
\end{lemma}

% We denote by $\calK_M$ the set of pairs of nodes that have non-zero connection probability.
% such that the rows $\bQs_{k,*}$, $k=1, \cdots, M$, of matrix  $\bQs$ 
% (and, therefore, the rows $\bDs_{k,*}$, $k=1, \cdots, M$, of matrix  $\bDs$) are not equal to identical zero.
% Let $M_0 = |\calK_M|$ be the cardinality of this set  (with $M_0=M$ if all groups in the network communicate with one another).
Then, one can use the same penalty that was considered in Corollary  \ref{cor:upper_DSBM}, so that Theorem~\ref{th:oracle}  yields the following result.

\begin{theorem} \label{th:upper_smooth_DSBM}
Consider a balanced DSBM satisfying condition \fr{balanced}. 
% Let $\calK^{*}_{\Ms}$ be the set of pairs of nodes such that $\bDs_{k,*} \equiv \bzero$ if $k \not\in \calKsMs$. 
Let  $\bLams$ be the true matrix of probabilities,  $\ms$ be the true number of classes, $\Ms = \ms(\ms +1)/2$, $\bQs$ be
the true matrix of probabilities of connections for pairs of classes and $\bDs = \bQs \bH$.
If $n \geq \sqrt{e\,n_0^3}$ and the  estimator $\hbLam$ is obtained as a solution of optimization problem \fr{opt_problem} with 
the penalty \fr{specific_pen} where  
\be \label{eq:J_union}
J = \bigcup_{k=1}^M J_k,
\ee 
then, for any $t>0$, with probability at least $1 - 9 e^{-t}$, one has
\be \label{oracle_prob_specific} 
{\small \frac{\|\hbLam - \bLams\|^2}{n^2 L} \leq  
 \min_{J} \lfi   \frac{6 \aleph_2^2}{(\ms)^2 L}   \sum_{k=1}^{\Ms}   \sum_{l \notin J_k} (\bDs_{k,l})^2  +  \frac{4  \Pen(|J|,\ms)}{n^2 L} \rfi  + 
\frac{38  t}{n^2 L} }
% \min_{J} \lfi   \frac{6 \aleph_2^2}{(\ms)^2\,L}\,  \sum_{k=1}^{\Ms} \, \sum_{l \notin J_k} (\bDs_{k,l})^2  +  \frac{4\, \Pen(|J|,\ms)}{n^2\,L} \rfi  + 
% \frac{38\, t}{n^2\,L} } 
\ee
 and a similar result holds for the expectation. 
\end{theorem}

In order to obtain  specific upper bounds in \fr{oracle_prob_specific}, we need to impose some assumptions on 
the smoothness of functions $\bQ_{k, *}^*$, $k=1, \cdots, \Ms$.
For  the sake of brevity, we assume that all   vectors  $\bDs_{k,*},$ $k=1, \cdots, \Ms$, behave similarly with respect to the basis $\bH$ 
(generalization to the case where this is not true is rather pedestrian but very cumbersome
as we point out in Section \ref{sec:discussion}, Discussion).
\\

{\bf (A0). } There exist   absolute constants $\nu_0$ and $K_0$ such that  
\be \label{bias_coef_cond}
\sum_{l=1}^L (l-1)^{2 \nu_0}\, (\bDs_{k,l})^2     \leq   K_0, \quad k =1, \cdots, \Ms. 
\ee  

\begin{corollary} \label{cor:upper_smooth_DSBM}
Let conditions of Theorem \ref{th:upper_smooth_DSBM} hold and 
$\bDs_{k,*}$ satisfy assumption \fr{bias_coef_cond}.   
If  the  estimator $\hbLam$ is obtained as a solution of optimization problem \fr{opt_problem}, 
then for any $t>0$, with probability at least $1 - 9 e^{-t}$, one has
\begin{align}  \nonumber 
\frac{\|\hbLam - \bLams\|^2}{n^2\,L} & \leq \tilde{K}_0\, \lkr  
\min \lfi \frac{1}{L}\, \lkv \lkr \frac{\ms}{n}\rkr^2 \, \log \lkr \frac{n}{\ms} \rkr \rkv^{\frac{2\nu_0}{2 \nu_0 +1}},
\lkr \frac{\ms}{n}\rkr^2 \rfi \right. \\
& + \left.  \frac{\log \ms}{nL} + \frac{n_0}{n^2} \log \lkr \frac{\ms\, n e}{n_0} \rkr + \frac{t}{n^2 L} \rkr
\label{oracle_prob_smooth}
 \end{align}
 and a similar result holds for the expectation. Here, $\tilde{K}_0$ is an absolute constant that depends on 
$K_0$, $\nu_0$, $\aleph_1$ and $\aleph_2$ only.
\end{corollary}

\ignore{ 
\begin{remark} \label{rem:VarCoef}
{\bf (Relation to the sparse varying coefficient model). }
{\rm  If $m=n$, the model studied in this paper reduces to the sparse varying coefficient model studied in, e.g., 
Klopp and Pensky  (2015). Note, however, that the problem studied in this paper is much more difficult. To start with, 
there is a block structure imposed on all variables. In addition, the number of blocks is unknown and 
the block memberships of the nodes are unknown and  can change in time. To appreciate the degree of the added complexity, compare the paper of 
Gao \etal (2015) to a paper on a well known classical regression (and Gao \etal (2015) does not have an added difficulty
of block memberships  switching in time and the number of blocks being unknown). 
}
\end{remark}
}

%%%%%%%%%%%%%%%%%%%%%%%%%%%%%%%%%%%%%%%%%%%%%%%%%%%%%%%%%%%%%%%%%%%%%%%%%%%%%%%%%%%%%%%%%%%%%%%%%

\section{The lower bounds for the risk  for the DSBM}
\label{sec:DSBM_lower bounds}

In order to prove that the estimator obtained as a solution of optimization problem 
\fr{opt_problem} is minimax optimal, we need to show that  the upper bounds in Corollaries~\ref{cor:upper_DSBM}
and \ref{cor:upper_smooth_DSBM} coincide with the minimax lower bounds obtained under similar constraints. 
For the sake of derivation of   lower bounds for the error,   we impose   mild conditions on   
the orthogonal   matrix $\bH$  as follows: for any binary vector
$\bom \in \{0,1\}^L$ one has 
\be \label{H_assump}
\|\bH^T \bom \|_\infty \leq \|\bom\|_1 /\sqrt{L} \quad \mbox {and}\quad
\bH \, \bone = \sqrt{L} \boe_1,
%\quad  \bone = (1,1, \cdots, 1)^T,\ \boe_1 = (1,0, \cdots, 0)^T.
\ee
where $\bone = (1,1, \cdots, 1)^T$ and $\boe_1 = (1,0, \cdots, 0)^T$.
Assumptions \fr{H_assump} are not restrictive. In fact, they are satisfied for a variety of common 
orthogonal transforms such as the Fourier transform or a periodic wavelet transforms.

First, we derive the lower bounds for the risk under the assumption that vector $\bd$ is $l_0$-sparse
and has only $s$ nonzero components. 
Let $\calG_{m,L,s}$ be a collection of tensors such that $\bG \in  \calG_{m,L,s}$ 
implies that the vectorized versions $\bq$ of $\bG$ can be written as 
$\bq = \bW^T \bd$ with $\|\bd\|_0 \leq s$. 
In order to be more specific, we consider the collection of clustering matrices 
$\calZ (m,n,n_0,L)$ with cardinality given by \fr{card_clust} that corresponds to the situation 
where   at most $n_0$ nodes can change their memberships between consecutive time instants.
In this case, $\Pen(|J|,m)$ is defined in \fr{specific_pen}.

\begin{theorem} \label{th:DSBM_lower_bounds}
Let orthogonal matrix $\bH$ satisfy condition \fr{H_assump}.
Consider the DSBM   where $\bG \in  \calG_{m,L,s}$  with $s \geq \kappa m^2$ where $\kappa>0$ is 
independent of $m$, $n$ and $L$. Denote $\ga = \min(\kappa, 1/2)$ and assume that 
$L \geq 2$, $n \geq 2m$, $n_0 \leq \min(\ga n, 4/3\, \ga n\,m^{-1/9})$ and $s$ is such that
\be \label{Js_cond}
s^2 \log(2 L M/s) \leq 68 L M n^2.
\ee 
Then
\begin{align} \label{lower_DSBM}
\inf_{\hbLam} \sup_{\stackrel{\bG \in \calG_{m,L,s}}{\bC \in \calZ (m,n,n_0,L)}}
\PP_{\bLam} & \lfi \frac{\|\hbLam - \bLam \|^2}{n^2\,L} \right.   \geq   C(\ga)   \lkr 
\frac{\log m}{n L} + \frac{n_0}{n^2}  \log \lkr \frac{m n e}{n_0} \rkr    \right.\\
%\frac{\log|\calZ (m,n,n_0,L)|}{n^2 L}
&+   \left. \left. \frac{s\, \log (L  m^2/s)}{n^2 L} \rkr \rfi \geq \frac{1}{4}, \nonumber
\end{align}
 where $\hbLam$ is any estimator of $\bLam$, $\PP_{\bLam}$ is the probability under the true value 
of the tensor $\bLam$ and $C(\ga)$ is an absolute constant that depends on $\ga$ only.
% and $|\calZ (m,n,n_0,L)|$ is provided in  \fr{card_clust}.
 \end{theorem}

Theorem \ref{th:DSBM_lower_bounds} ensures that if vector $\bd$ has only $s$ nonzero components,  
then the upper bounds in Corollary~\ref{cor:upper_DSBM} are optimal up to a constant.
In order to provide a similar assertion in the case of Corollary \ref{cor:upper_smooth_DSBM}, we assume that rows of  matrix $\bD$ are 
$l_2$-sparse. For this purpose, we consider a collection of tensors $\calG_{m,L,\nu_0}$ such that $\bG \in  \calG_{m,L,\nu_0}$ 
implies that $\bQ = \bD \bH$ and rows $\bD_{k,*}$ of matrix $\bD$ satisfy condition \fr{bias_coef_cond}.
Let as before $\calZ_{\bal} (m,n,n_0,L,\aleph_1, \aleph_2)$ be a collection of clustering matrices 
satisfying condition \fr{balanced} and  such that at most $n_0$ nodes change their memberships between two 
consecutive time instances. The following statement ensures that the upper bounds in Corollary \ref{cor:upper_smooth_DSBM}
are minimax optimal up to a constant factor.

\begin{theorem} \label{th:smooth_DSBM_lower_bounds}
Let orthogonal matrix $\bH$ satisfy condition \fr{H_assump}.
Consider the DSBM   where $\bG \in  \calG_{m,L,\nu_0}$  with $\nu_0 > 1/2$,
$L \geq 2$ and  $n \geq 2m$. 
Then, for any absolute constants $0 < \aleph_1 \leq 1 \leq  \aleph_2 < \infty$, one has
\begin{align} \label{smooth_lower_DSBM}
\inf_{\hbLam} \sup_{\stackrel{\bG \in \calG_{m,L,s}}{\bC \in \calZ_{\bal}}}
\PP_{\bLam} & \lfi \frac{\|\hbLam - \bLam \|^2}{n^2\,L} \right.   \geq 
   C\, \lkv \min \lfi \frac{1}{L} \lkv \lkr \frac{m}{n}\rkr^2  \rkv^{\frac{2 \nu_0}{2 \nu_0+1}};   
\lkr \frac{m}{n}\rkr^2 \rfi \right. \\
& \left. \left.  + \frac{\log m}{nL} + \frac{n_0}{n^2}  \log \lkr \frac{m n e}{n_0} \rkr \rkv \rfi \geq \frac{1}{4},
\nonumber
\end{align}
 where $\calZ_{\bal}$ stands for $\calZ_{\bal} (m,n,n_0,L,\aleph_1, \aleph_2)$,
$\hbLam$ is any estimator of $\bLam$, $\PP_{\bLam}$ is the probability under the true value 
of the tensor $\bLam$ and $C$ is an absolute constant independent of $n$, $m$ and $L$.
 \end{theorem}

Theorems~\ref{th:DSBM_lower_bounds} and \ref{th:smooth_DSBM_lower_bounds} 
confirm that the estimator constructed above is {\bf minimax optimal up to a constant}
if $\bG \in \calG_{m,L,s}$ and $\bC \in \calZ (m,n,n_0,L)$, or $\bG \in  \calG_{m,L,\nu_0}$
and $\bC \in \calZ_{\bal} (m,n,n_0,L,\aleph_1, \aleph_2)$.

Note that the   terms $\log m/(nL)$  and $n_0 n^{-2}\, \log(mne/n_0)$   in \fr{lower_DSBM}   and \fr{smooth_lower_DSBM}
correspond  to, respectively, the error of initial clustering and the clustering error due to membership changes.
The remaining  terms are due to   nonparametric estimation and model selection.
Assumptions \fr{H_assump} and \fr{Js_cond} are purely technical and are necessary to 
ensure that  the ``worst case scenario'' tensor $\bG$ 
of connection probabilities has nonnegative components. As we mentioned earlier, 
conditions \fr{H_assump} are totally non-restrictive. 
Condition \fr{Js_cond} in Theorem~\ref{th:DSBM_lower_bounds} holds 
whenever representation of the tensor of probabilities in the basis $\bH$ is at least somewhat sparse.
Indeed, if there is absolutely no sparsity (which is a very implausible  scenario when smooth functions are represented in a basis) 
and $s \approx M L$, then condition \fr{Js_cond} reduces to  $m (m+1) L   \leq C n^2$ and will still be true if $L$ is relatively small. 
If $L$ is large, the situation where $s \approx M L$ is very unlikely. 
Assumption that $s \geq \kappa m^2$ for some  $\kappa>0$ independent of  
$m$, $n$ and $L$,  restricts the sparsity level and ensures that one does not 
have too many classes where nodes have no interactions with each other or members
of other classes.

Finally, it is also worth  keeping in mind that 
all assumptions in Theorems~\ref{th:DSBM_lower_bounds} and \ref{th:smooth_DSBM_lower_bounds}   are used for the derivation 
of the minimax lower bounds for the risk and are not necessary
for either the  construction of the estimator $\hbLam$ of $\bLam$ in \fr{opt_problem} 
 or for the assessment of its  precision in Theorems~\ref{th:oracle} and \ref{th:upper_smooth_DSBM}.

%%%%%%%%%%%%%%%%%%%%%%%%%%%%%%%%%%%%%%%%%%%%%%%%%%%%%%%%%%%%%%%%%%%%%%%%%%%%%%%%%%%%%%%%%%%%%%%%%

\section{The uniformly sparse DSBM}
% \section{Sparse DSBMs}
\label{sec:sparse}

In the current literature, the notion of the sparse SBM  refers to the case where the entries of the matrix 
of the connection probabilities are uniformly small: $\bLam = \rho_n \bLam^{(0)}$ with $\| \bLam^{(0)} \|_\infty =1$ and 
$\rho_n \to 0$ as $n  \to \infty$. The concept is based on the idea that when the number of nodes in a network grow,
the probabilities of connections between them decline. The minimax study of the sparse SBM has been carried out
by Klopp \etal (2017). The logical generalization of the sparse SBM of this type would be the sparse DSBM where
the elements of the tensor $\bLam$ are bounded above by $\rho_n$ where  $\rho_n \to 0$ as $n$ grows.
We refer to this kind of network as {\it uniformly sparse}. 

On the other hand,  not all networks become uniformly sparse as  $n  \to \infty$.
Indeed, in the real world, when a  network grows, the number of communities increase and, while 
the probabilities of connections for majority of pairs of groups become very small, some of the of pairs groups 
will still maintain  high connection probabilities. We refer to this type of network as {\it non-uniformly sparse}.  
 The idea of such a  network has been  elaborated in the recent paper of
Borgs \etal (2016). The authors  considered  heavy-tailed sparse graphs such that, in the context of the SBM, 
one still has $\bLam = \rho_n \bLam^{(0)}$ but the elements of $\bLam^{(0)}$ are no longer bounded by one
but by a quantity that grows with $n$.

While distinguishing between very small probabilities might be essential in a clustering problem,
it is not so necessary in the problem of estimation of the tensor of the connection probabilities 
studied in the present paper. Indeed, it is a common knowledge that, in the nonparametric regression model, 
in order to obtain the best error rates, one needs to replace   small elements of the vector of interest 
by zeros rather than estimating them. Similarly, if the network is non-uniformly sparse, i.e., some pairs of
groups have probabilities of connections  equal or very close to zero, one would obtain an estimator with   
better overall precision by setting those very small connection probabilities to zeros.   
Although nowhere in the present paper  we   make  an assumption that a network is sparse and, moreover, consideration of 
the non-uniformly sparse SBM or DSBM is not one of its objectives, the paper naturally provides the tools
for minimax optimal statistical estimation  in such models that deliver results with very little  
additional work.

In addition, the techniques developed in this paper allow, with some additional work, to extend results obtained 
in Klopp \etal (2017) to the dynamic setting. However, majority of their results depend upon solution of 
optimization problem \fr{opt_problem} under the restriction that $\|\bW^T \bd\|_{\infty} \leq \rho_n$ 
which requires representation of the estimator via a different projection operator and will result in more cumbersome calculations. 
Therefore, we avoid studying this new optimization problem and only extend Corollary 2.2    of 
Klopp \etal (2017)  that handles the case of the balanced model without placing the above-mentioned restriction.  
For this purpose, consider a small $\rho_n$  and denote
\be \label{eq:rnm}
r_n(m) = \max(\rho_n, m^2/n^2).
\ee
Similarly to \fr{opt_problem}, we find $m, J, \bd$ and $\bC$  as one of the solutions of the following penalized least squares 
optimization problem
{\small 
\be \label{opt_problem_balanced}
 (\hm, \hJ,\hbd, \hbC) = \underset{m,J,\bd,\bC}{\operatorname{argmin}} \lkv \|\ba - \bC \bW^T \bd\|^2 +  \lam_0  r_n(m)   \Pen(|J|,m)\rkv
\  \mbox{s.t.}\    \bd_{J^c}=\bzero    
\ee}

\noindent
where $\bC \in \calZ_{\bal} (m,n,n_0,L,\aleph_1,\aleph_2)$, $\ba$ is defined in \fr{full_model},  
$\bd \in \RR^{ML}$, $\bW \in \RR^{ML \times ML}$, $M=m(m+1)/2$, $\Pen(|J|,m)$  is defined in \fr{specific_pen}
and $\lam_0$ is a tuning parameter that is bounded above and below by a constant.
\\

In order the estimator has the uniform sparsity property, we need to make sure that transformation $\bH$ is such that, 
whenever it is used for sparse representation of smooth functions, the maximum absolute value of the estimator  
obtained by truncation of the vector of coefficients is bounded above by a constant factor of the  
maximum absolute value of the original function. In particular, we denote the projection matrix 
on the column space of matrix $(\bC \bW^T)_J$ by $\PCJ$ and impose the following condition on the transformation matrix $\bH$:
\\

{\bf (A1). }   There exists an absolute  constant $B_0$  such that  for any 
$\bC \in \calZ_{\bal} (m,n,n_0,L,\aleph_1,\aleph_2)$ and any vector $\bte$ 
\be \label{uni_sparse_cond}
\| \PCJo \bte \|_\infty = \| \te - \PCJ \te \|_\infty \leq B_0 \| \te\|_\infty.
\ee

% Note that assumption \fr{uni_sparse_cond}  always holds in the case of $L=1$. 

Let, as before, $\bLams$ be the true matrix of probabilities,    $\ms$  be the true number of classes, 
$\Ms = \ms(\ms +1)/2$, $\bCs$  be the true clustering matrix, $\bQs$  be the true matrix of 
probabilities of connections for pairs of classes, $\bDs = \bQs \bH$, $\bds = \vect(\bDs)$,
$\btes = \bCs (\bWs)^T \bds$ and $\bWs = \bH \otimes \bI_{\Ms}$. 

\begin{theorem} \label{th:sparse}
Consider a balanced DSBM satisfying condition \fr{balanced}.  Let matrix $\bH$ be such that  condition 
\fr{uni_sparse_cond} is satisfied  and $\| \bLams \|_\infty \leq \rhons$.
If $\rho_n \geq \rhons$, $n \geq \sqrt{e\,n_0^3}$ and the  estimator $\hbLam$ is obtained as a solution of optimization problem \fr{opt_problem_balanced},
% with the penalty \fr{specific_pen} where  $J$ is given by \fr{eq:J_union} and $\Pen(|J|,m)$  is defined in \fr{specific_pen}.
then, for an absolute constant $\tilde{C}_0$ and any $t>0$, with probability at least $1 - 9 e^{-t}$, one has
\be \label{oracle_prob_sparse1} 
\frac{\|\hbLam - \bLams\|^2}{n^2\,L} \leq  \tilde{C}_0\, 
 \min_{J} \lfi   \frac{\| {\bPi}_{\bC^*,J}^{\bot} \btes \|^2}{n^2\,L}   +  
\frac{r_n(\ms)   \, [\Pen(|J|,\ms)   + t]}{n^2\,L} \rfi
\ee
where ${\bPi}_{\bC^*,J}$ is the projection matrix on the column space of $(\bCs \bWs^T)_J$ and $\tilde{C}_0$
is an absolute constant that depends on  $B_0$, $\aleph_1$ and $\aleph_2$ only.

In particular, if condition \fr{bias_coef_cond} holds with $K_0$ replaced with $\rhons K_0$, then
 \begin{align}  \nonumber 
\frac{\|\hbLam - \bLams\|^2}{n^2\,L} & \leq \tilde{K}_0  r_n(\ms) \, \lkr  
\min \lfi \frac{1}{L}\, \lkv \lkr \frac{\ms}{n}\rkr^2 \, \log \lkr \frac{n}{\ms} \rkr \rkv^{\frac{2\nu_0}{2 \nu_0 +1}},
\lkr \frac{\ms}{n}\rkr^2 \rfi \right. \\
& + \left.  \frac{\log \ms}{nL} + \frac{n_0}{n^2} \log \lkr \frac{\ms\, n e}{n_0} \rkr + \frac{t}{n^2 L} \rkr
\label{oracle_prob_sparse2}
 \end{align}
Here, $\tilde{K}_0$ is an absolute constant that depends on $B_0$, $K_0$, $\nu_0$,  $\aleph_1$ and $\aleph_2$ only. 
Results similar to \fr{oracle_prob_sparse1}  and \fr{oracle_prob_sparse2} hold for the expectations.
\end{theorem}

%%%%%%%%%%%%%%%%%%%%%%%%%%%%%%%%%%%%%%%%%%%%%%%%%%%%%%%%%%%%%%%%%%%%%%%%%%%%%%%%%%%%%%%%%%%%%%%%%

\section{Dynamic graphon estimation}
\label{sec:dyn_graphon}

Consider the situation where   tensor  $\bLam$ is generated by a dynamic graphon $f$,
so that $\bLam$ is given by expression \fr{graphon}
% \be \label{graphon}
% \bLam_{i,j,l} = f(\zeta_i, \zeta_j, t_l), \quad i,j = 1, \cdots,n, \ l=1, \cdots, L,
% \ee
where function $f: [0,1]^3 \to [0,1]$ is such that  
$f(x,y,t) = f(y,x,t)$  for any $t$  and 
$\bzeta = (\zeta_1, \cdots, \zeta_n)$ is a random vector sampled from a distribution $\PP_\zeta$
supported on $[0,1]^n$.

Given an observed adjacency tensor $\bB$ sampled according to model \fr{graphon}, 
the graphon function $f$ is not identifiable since the topology of  a network 
is invariant with respect to any change of labeling of its nodes. Therefore, 
for any $f$ and any measure-preserving bijection $\mu: [0,1]\to [0,1]$
(with respect to Lebesgue measure), the functions $f(x,y,t)$ and $f(\mu(x), \mu(y),t)$
define the same probability distribution on random graphs. For this reason, we are 
considering equivalence classes of graphons. Note that in order for it to be possible to compare 
clustering of nodes across time instants, we introduce an assumption 
% (which, to the best of our knowledge, first appeared in Matias and Miele (2015)) 
that there are no label switching in time, that is, every  
node carries the same label at any time $t_l$, so that function $\mu$ is independent of $t$.

Under this condition, we further assume that probabilities $\bLam_{i,j,l}$ do not change drastically 
from one time point to another, i.e. that, for every $x$ and $y$,
functions $f(x,y,t)$ are smooth in $t$. We shall also assume that $f$ is piecewise smooth
in $x$ and $y$. 
In order to quantify those assumptions, for each $x, y \in [0,1]^2$,  we consider a vector 
$\bof(x,y) = (f(x,y,t_1), \cdots,f(x,y,t_L))^T$ and   an orthogonal transform $\bH$ used in the previous sections.
We assume that elements $\bv_l(x,y)$ of vector $\bv(x,y) = \bH \bof(x,y)$ satisfy the following assumption:

{\bf (A2). } There exist constants $0 = \beta_0 < \beta_1 < \cdots < \beta_r =1$ and  $\nu_1, \nu_2, K_{1}, K_2>0$  
such that for any $x,x' \in (\beta_{i-1}, \beta_i]$ and  $y,y' \in (\beta_{j-1}, \beta_j]$, $1 \leq i,j \leq r$,
one has
\beqn    
 [\bv_l(x,y)- \bv_l(x',y')]^2   & \leq &   K_1  [|x - x'| + |y - y'|]^{2 \nu_1}, \label{A1}\\ 
\sum_{l=1}^L (l-1)^{2 \nu_2}\, \bv_l^2(x,y)  & \leq & K_2. \label{A2}
\eeqn 
Note that,  for a graphon corresponding to the  DSBM model,
on each of the rectangles $(\beta_{i-1}, \beta_i] \times (\beta_{j-1}, \beta_j]$,
functions $\bv_l(x,y)$ are constant, so that $\bv_l(x,y) = 0$ for $l=2, \cdots, L$, and $\nu_1 = \infty$. 
\\

% \noindent 
% Observe that Condition \fr{A2} is somewhat similar to \fr{bias_coef_cond} in Assumption {\bf (A1)}
% since both assumptions 

We denote the class of graphons satisfying assumptions \fr{graphon}, \fr{A1} and \fr{A2}   by $\Sig(\nu_1, \nu_2, K_1, K_2)$.
In order to estimate  the dynamic graphon, we approximate it by an appropriate DSBM   
and then estimate the probability tensor of the DSBM.   
Note that, since $\nu_1, \nu_2, K_{1}$ and $K_2$
in Assumption {\bf A} are independent of $x$ and $y$, one can simplify the optimization procedure in \fr{opt_problem}.

Let $\bQ$ be the matrix defined in \fr{con-memb} and \fr{full_model}. 
Note that since random variables $\zeta_1, \ldots, \zeta_n$ are time-independent,
we can approximate the graphon by a DSBM where group memberships of the nodes 
do not change in time. Hence, , matrices $\bC^{(l)}$ are independent of $l$, 
so that $\bC^{(l)}= \bZ$,   \fr{con-memb} holds and $\bTe = \bZ \bQ$. 
Denote $\bX = \bA \bH^T$. 
Denote by $\bV$ and $\bPhi$ the matrices of the coefficients of $\bQ$ and $\bTe$ in the transform $\bH$:
$\bV = \bQ \bH^T$ and $\bPhi = \bTe \bH^T$. 
% Note that $\bV$ and $\bPhi$ represent the coefficients of, respectively, $\bQ$ and $\bTe$ in the transform $\bH$.  
Then, by \fr{con-memb}, $\bTe =   \bZ \bV \bH$ and $\bPhi =  \bZ \bV$.
Note that each row of the matrices $\bV$ and $\bPhi$ corresponds to one spatial location.
Since, due to \fr{A2}, the coefficients in the transform $\bH$ decrease uniformly irrespective of the location,
one can employ   $L_1 <L$ columns instead of $L$ columns  in the final representations of $\bQ$ and $\bTe$.
In order to simplify our presentation, we denote $L_1 = L^\rho$ where $0 < \rho \leq 1$
and use the optimization procedure \fr{opt_problem} to find $m,\rho, \bV^{(\rho)}$ and $\bZ$
where $\bV^{(\rho)}$ is the submatrix of $\bV$ with columns $\bV_{*,j}$, $1 \leq j \leq L^\rho$.
Due to $|J| = M L^\rho$ and $0.5\, m^2 L^\rho \leq |J|  \leq m^2 L^\rho$,  
 in this case   optimization problem \fr{opt_problem} can be reformulated as
\begin{align} \label{graphon_opt}
(\hm, \hrho,\hbV^{(\hrho)}, \hbZ) & = \underset{\stackrel{m,\rho,\bV^{(\rho)}}{\bZ\in \calZ (m,n,0,L)}}
{\operatorname{argmin}} 
\lkv \|\bX^{(\rho)} - \bZ \bV^{(\rho)} \|^2 \right. \\
& \left.+ 11 n \log m + \frac{11}{2}  m^2 L^\rho \log(25  L^{1-\rho})\rkv \nonumber
% \quad \mbox{s.t.}\   \bZ \in \calZ (m,n,0,L)
\end{align}
where $\calZ (m,n,0,L)$ is defined in \fr{card_clust}. 
Then the estimation algorithm appears as follows:

\begin{enumerate}
\item
Apply  transform $\bH$ to the data matrix $\bA$ obtaining matrix $\bX = \bA \bH^T$.
\item
Consider a set $\Re = \lfi \rho \in (0,1]:\ L^{\rho} \quad \mbox{is an integer} \rfi$.
For every $\rho \in \Re$, remove all columns $\bX_{*,l}$ with $l \geq L^\rho+1$  obtaining matrix $\bX^{(\rho)}$ with 
$\EE \bX^{(\rho)} = \bZ \bV^{(\rho)} \equiv \bPhi^{(\rho)}$
where matrix $\bV^{(\rho)}$ has $L^\rho$ columns. 
\item
Find $(\hm, \hrho,\hbV^{(\hrho)}, \hbZ)$ as a solution of the optimization problem \fr{graphon_opt}.
\item 
Choose $\hbTe = \hbZ \hbV^{(\hrho)} \bH$ and obtain $\hbLam$ by packing $\hbTe$ into a tensor. 
\end{enumerate}

\noindent
Note  that construction of the estimator $\hbLam$ does not require  knowledge 
of  $\nu_1, \nu_2, K_{1}$ and $K_2$, so the estimator is fully adaptive.
The following statement provides a  minimax  upper bound for the risk of   $\hbLam$.

\begin{theorem} \label{th:upper_graphon}
Let $\Sig \equiv \Sig(\nu_1, \nu_2, K_1, K_2)$ be the class of   graphons 
satisfying Assumptions \fr{graphon}, \fr{A1} and \fr{A2}. 
If $\hbLam$ is   obtained as a solution of optimization problem  \fr{graphon_opt} as described above, then
\begin{align} \label{graph_upper}
 \sup_{f \in \Sig}  \frac{\EE \|\hbLam - \bLams\|^2}{n^2\,L} & \leq 
C \min_{\stackrel{1 \leq h \leq n-r}{0 \leq \rho \leq 1}} \lfi
\frac{L^{\rho-1}}{h^{2 \nu_1}} + \frac{I(\rho <1)}{L^{2 \rho \nu_2+1}} \right. \\
& + \left.
\frac{(h+r)^2 (1 + (1 - \rho) \log L) }{n^2 L^{1-\rho}}  
+ \frac{\log(h+r)}{n  L}  \rfi, \nonumber
\end{align}
where  the constant $C$ in \fr{graph_upper} depends on $\nu_1, \nu_2, K_{1}$ and $K_2$ only. 
\end{theorem}

Note that $h$ in \fr{graph_upper} stands for $h = m-r$ where $m$ is the number  of
blocks in the DSBM which approximates the graphon, hence, $h \leq n-r$.  
On the other hand, $h\geq 0$ since one needs at least $r$ blocks to approximate the graphon 
that satisfies condition \fr{A1}. Since the  expression in the right hand side 
of \fr{graph_upper} is  rather complex and is hard to analyze,
we shall consider only two regimes: a)\, $r = r_{n,L} \geq 2$ may depend on $n$ and $L$   
and $\nu_1 = \infty$; or b)\, $r = r_0\geq 1$ is a fixed  quantity independent of $n$ and $L$.
The first regime corresponds to a piecewise constant (in $x$ and $y$) graphon that generates the DSBM
while the second   regime deals with the situation where $f$ is a piecewise smooth function 
of all three arguments with a finite number of jumps.
In the first case, we set $h=2$, in the second case, we choose $h$ to be a function of $n$ and $L$. 
By minimizing the right-hand side of \fr{graph_upper}, we obtain the following statement.

\begin{corollary} \label{cor:upper_graphon}
Let $\hbLam$ be   obtained as a solution of optimization problem  \fr{graphon_opt} as described above. 
Then, for $\Sig \equiv \Sig(\nu_1, \nu_2, K_1, K_2)$ and $C$ independent of $n$ and $L$, one has 
{\small \be \label{graph_upper_cases}
 \sup_{f \in \Sig} \,  \frac{\EE \|\hbLam - \bLams\|^2}{n^2\,L} \leq 
\lfi \begin{array}{l}
C \min \lfi \frac{1}{L} \lkv \lkr \frac{r}{n}\rkr^2 \log \lkr \frac{n}{r}\rkr \rkv^{\frac{2 \nu_2}{2 \nu_2+1}};   
\lkr \frac{r}{n}\rkr^2 \rfi + \frac{C  \log r}{nL}, \  r = r_{n,L};\\
  \\
C \min \lfi \frac{1}{L}   \lkr \frac{\log L}{n^2}\rkr^{\frac{2 \nu_1 \nu_2}{(\nu_1 +1) (2\nu_2+1)}};   
\lkr \frac{\log L}{n^2}\rkr^{\frac{\nu_1}{\nu_1 +1}} \rfi + \frac{C  \log n}{nL},  \  r = r_0.
\end{array} \right.
\ee}
\end{corollary}

\noindent
In order to assess optimality of the penalized least squares estimator obtained above, we 
derive lower bounds for the minimax risk over the set $\Sig(\nu_1, \nu_2, K_1, K_2)$.
These lower bounds are constructed separately for each of the two regimes.

\begin{theorem} \label{th:lower_graphon}
Let matrix $\bH$ satisfy assumptions \fr{H_assump} and $\nu_2 \geq 1/2$ in \fr{A2}. Then, 
for   $C$ independent of $n$ and $L$, one has 
\be \label{graph_lower_ineq}
\inf_{\hbLam}\ \sup_{f \in \Sig(\nu_1, \nu_2, K_1, K_2)} \  
\PP_{\bLam} \lfi \frac{\|\hbLam - \bLam \|^2}{n^2\,L} \geq   \Del(n,L)  \rfi \geq \frac{1}{4},
\ee
where
\be \label{graph_lower_cases}
\Del(n,L) = \lfi \begin{array}{l}
 C\, \min \lfi \frac{1}{L} \lkv \lkr \frac{r}{n}\rkr^2  \rkv^{\frac{2 \nu_2}{2 \nu_2+1}}; \ 
\lkr \frac{r}{n}\rkr^2 \rfi + \frac{C\, \log r}{nL}, \  r = r_{n,L};\\
  \\
C\,   \min \lfi \frac{1}{L}   \lkr \frac{1}{n^2}\rkr^{\frac{2 \nu_1 \nu_2}{(\nu_1 +1) (2\nu_2+1)}}; \ 
\lkr \frac{1}{n^2}\rkr^{\frac{\nu_1}{\nu_1 +1}} \rfi + \frac{C\, \log n}{nL}, \  r = r_0.
\end{array} \right.
\ee
\end{theorem}

 \noindent
It is easy to see that   the value of $\Del(n,L)$ coincides with the upper bound in 
\fr{graph_upper_cases} up to a  at most a  logarithmic factor of $n/r$ or $L$.
In both cases, the  first quantities in the minimums correspond  to the situation where $f$ is smooth enough as a function of time, so that 
application of transform $\bH$ improves estimation precision by reducing the number of parameters 
that needs to be estimated. The second quantities represent  the case where one needs to keep all elements of  vector $\bd$
and hence application of the transform yields no benefits. The latter can be due to the fact that  $\nu_2$ is too small or $L$ is too low.

The upper and the lower bounds in Theorems \ref{th:upper_graphon} and \ref{th:lower_graphon}
look somewhat similar to the ones appearing in  anisotropic functions estimation (see, e.g., Lepski  (2015)).
Note also that although in the case of a time-independent   graphon ($L=1$), the estimation precision does not improve 
if $\nu_1>1$, this is not true any more in the case of a dynamic graphon. Indeed,  the right-hand sides in 
\fr{graph_lower_cases} become significantly smaller when   $\nu_1, \nu_2$ or $L$ grow.

\begin{remark} \label{rem:DynGraphon}
{\bf (The DSBM and the dynamic graphon).}
{\rm 
Observe that   the definition \fr{graphon} of the dynamic graphon 
assumes that vector $\bzeta$ is independent of $t$.
This is due to the fact that, to the best of our knowledge, the notion of the dynamic graphon
with $\bzeta$ being a function of time has not yet been developed by the probability community. 
For this reason,  we   restrict  our attention to the case where we are certain that, at any time point, the graphon 
describes the limiting behavior of the network  as $n \to \infty$.
Nevertheless, we believe that when the concept of the dynamic graphon is established, 
our techniques will be useful for its  estimation.

In the case of a piecewise constant graphon,  our setting corresponds to the 
situation where the nodes of the network do not switch their group memberships in time, so that $n_0 =0$
in \fr{card_clust}. Therefore,   a piecewise constant graphon ($r = r_{n,L},\ \nu_1 = \infty$) 
is just a particular case of the general DSBM since the latter allows any temporal changes of nodes' memberships. 
However, the dynamic piecewise constant graphon formulation enables us to derive specific minimax convergence 
rates for estimators of $\bLam$ in terms of $n$, $L$ and $r$.
On the other hand, the piecewise smooth graphon ($r = r_0,\ \nu_1 < \infty$) is an entirely different object 
that is not represented by the DSBM.
}
\end{remark}

%%%%%%%%%%%%%%%%%%%%%%%%%%%%%%%%%%%%%%%%%%%%%%%%%%%%%%%%%%%%%%%%%%%%%%%%%%%%%%%%%%%%%%%%%%%%%%%%%

\section{Discussion} 
\label{sec:discussion}

In the present paper we considered estimation of connection probabilities 
in the context of   dynamic network models.  To the best of our knowledge, 
this is the first paper to propose a fully non-parametric  model for the 
time-dependent networks which treats connection probabilities for each group as the
functional data  and allows to exploit the consistency in the group memberships over time.
The paper derives adaptive penalized least squares estimators of the tensor of the connection probabilities 
in a non-asymptotic setting  and shows that the estimators are indeed minimax optimal by
constructing the lower bounds for the risk. This is done via vectorization technique which 
is very useful for the task in the paper and can be very beneficial for solution of other problems such as, 
e.g., inference in   bi-clustering models mentioned in Remark~\ref{rem:biclust}.
In addition, we show that the correct penalty consists of two parts: the portion which accounts 
for the complexity of estimation  and the portion which accounts for the complexity of clustering and is 
proportional to the logarithm of the cardinality 
of the set of clustering matrices. The latter is a novel result and it is obtained by using the 
innovative Packing lemma   (Lemma~4) which can be viewed as a version 
% innovative Packing lemma   (Lemma~\ref{lem:packing}) which can be viewed as a version 
of the Varshamov-Gilbert  lemma for clustering matrices.   
Finally,  the methodologies of the paper  allow a variety of extensions.

\begin{enumerate}

\item
{\bf (Inhomogeneous or non-smooth connection probabilities). } 
Assumption \fr{A2} essentially implies that probabilities of connections are spatially homogeneous 
 and are represented by smooth  functions of time that belong to the same Sobolev class.
The model, however, can be easily generalized. First, by letting  $\bH$ be a wavelet transform and
assuming that for any fixed $x$ and $y$, function $f(x,y,\cdot)$  belongs  to a Besov ball, one can accommodate the case where
$f(x,y,\cdot)$ has   jump discontinuities. Furthermore, by using a weaker version of 
condition \fr{A2}, similarly to how this was done in Klopp and Pensky (2015), we can 
treat the case where   functions $f(x,y,t)$ are spatially inhomogeneous.

\item
 {\bf (Time-dependent  number of nodes). }
One can apply the theory above even when the number of nodes in the network changes from one time instant to another. 
Indeed, in this case we can form a set   that includes all nodes that have ever been in the network and 
denote their number by $n$. Consider a class $\Om_0$ such that all nodes in this class 
have zero probability of interaction with each other or any other node in the network. 
At each time instant, place all nodes that are not in the network into the class $\Om_0$. After that, 
one just needs to modify the optimization procedures  by placing additional restriction  that 
the out-of-the-network nodes indeed belong to class $\Om_0$ and that
$\bG_{0,k,l}=0$ for any $k=0,1,2,\cdots, m$ and $l=1, \cdots, L$.

\item
{\bf (Adaptivity to clustering complexity). }
Although, in the case of the DSBM, our estimator is adaptive
to the unknown number of classes, it requires  knowledge about 
the complexity of the set of clustering matrices.
% how much nodes' memberships can change from one time instant to another.
For example, if at most $n_0$ nodes can change their memberships between two consecutive 
time points  and  $n_0$ is a fixed quantity independent of $n$ and $m$, we 
can replace $n_0$ by $\log n$ that dominates $n_0$ if $n$ is large enough. 
However, if $n_0$ depends on $n$ and $m$,
 development of an adaptive estimator would require an additional investigation.

\end{enumerate}

\vspace{4mm}

%%%%%%%%%%%%%%%%%%%%%%%%%%%%%%%%%%%%%%%%%%%%%%%%%%%%%%%%%%%%%%%%%%%%%%%%%%%%

\noindent
{\bf SUPPLEMENTARY MATERIAL. }  Supplement contains proofs of all statements in the paper

%%%%%%%%%%%%%%%%%%%%%%%%%%%%%%%%%%%%%%%%%%%%%%%%%%%%%%%%%%%%%%%%%%%%%%%%%%%%

%%%%%%%%%%%%%%%%%%%%%%%%%%%%%%%%%%%%%%%%%%%%%%%%%%%%%%%%%%%%%%%%%%%%%%%%%%%%%%%%%%%%%%%%%%%%%%%%%%%%%%%%%%%%

\newpage

\pagestyle{plain}
\setcounter{page}{1}
\pagenumbering{arabic}

%%%%%%%%%%%%%%%%%%%%%%%%%%%%%%%%%%%%%%%%%%%%%%%%%%%%%%%%%%%%%%%%%%%%%%%%%%%%

%  This is the file with proofs
% The root file is DSBM_AOS_March25_2018.tex

\section{Supplemental Material:\ Proofs}
\label{sec:suppl}

\renewcommand{\theequation}{S.\arabic{equation}}
 \setcounter{equation}{0}

%%%%%%%%%%%%%%%%%%%%%%%%%%%%%%%%%%%%%%%%%%%%%%%%%%%%%%%%%%%%%%%%%%%%%%%%%%%%%%%%%%%%%%%%%%%%%%%%%

 %  This is the file with proofs
% The root file is DSBM_AOS_March25_2018.tex

% \section{Proofs}
% \label{sec:proofs}
%\setcounter{equation}{0}

% \subsection{Proofs of the upper bounds for the risk}
%\label{sec:upper}
\subsection{Proof of Theorem \ref{th:oracle}}

%  \noindent
% {\bf Proof of Theorem \ref{th:oracle}. } 
%
Since $(\hbd, \hbC, \hm, \hJ)$ are solutions of optimization problem \fr{opt_problem},
for any $m$, $J$, $\bC$  and $\bd$ one has
\be \label{main_ineq}
\|\ba - \hbC \hbW^T\, \hbd_{(\hJ)}\|^2 + \Pen(|\hJ|,\hm) \leq \|\ba - \bC \bW^T \bd_{(J)} \|^2 + \Pen(|J|,m),
\ee 
% where $\bd_{(J)}$ is the modification of vector 
% $\bd$ where all elements $\bd_j$ with $j \notin J$ are set to zero.
For any $m$, $J$, $\bC$  and $\bd$, it follows from \fr{true_model} that 
\begin{align*} 
& \|\ba -   \bC \bW^T  \bd_{(J)}\|^2  = \|(\ba - \bCs \bWs^T \bds) +  (\bCs \bWs^T \bds -  \bC \bW^T \bd_{(J)})\|^2 \\
& =   \|\bxi\|^2 + \|\bC \bW^T \bd_{(J)} - \bCs \bWs^T \bds\|^2 + 2 \bxi^T (\bCs \bWs^T \bds -  \bC \bW^T \bd_{(J)}).
\end{align*}
Hence, plugging the last identity into the inequality \fr{main_ineq}, derive that
for any $m$, $J$, $\bC$  and $\bd$
\begin{align} \label{main_ineq1}
\|\hbC \hbW^T \hbd_{(\hJ)} - \bCs \bWs^T \bds\|^2 & \leq \|\bC \bW^T \bd_{(J)} - \bCs \bWs^T \bds\|^2 \\
& +    \Del + \Pen(|J|,m) - \Pen(|\hJ|,\hm). \nonumber 
\end{align}
Here $\Del   = 2 |\bxi^T (\hbC \hbW^T \hbd_{(\hJ)} - \bC \bW^T \bd_{(J)})| \leq  \Del_1 + \Del_2 + \Del_3$ with 
\begin{align}
 % \Del & = 2 |\bxi^T (\hbC \hbW^T \hbd_{(\hJ)} - \bC \bW^T \bd_{(J)})| \leq  \Del_1 + \Del_2 + \Del_3, \label{Del_sum}\\
  \Del_1 & = 2 |\bxi^T (\bCs \bWs^T \bds - \bC \bW^T \bd_{(J)})|, \nonumber \\
  \Del_2  & = 2 |\bxi^T (\bI_{NL} - \hPhJ) \bCs \bWs^T \bds|, \quad 
  \Del_3   = 2  \bxi^T \hPhJ\, \bxi,  \label{Del123}
\end{align}
% \be \label{Del_sum}
% \Del   = 2 |\bxi^T (\hbC \hbW^T \hbd_{(\hJ)} - \bC \bW^T \bd_{(J)})| \leq  \Del_1 + \Del_2 + \Del_3,
% \ee
since, due to \fr{projection} and \fr{hbd_solution}, one has  $\hbC \hbW^T \hbd_{(\hJ)} = \hPhJ\, \ba$ 
where $\ba$ is given by \fr{true_model}. Now, we need to find upper bounds for each of the terms in \fr{Del123}.

By Lemma \ref{lem:bern_err} with $\al =1/2$ and any   $t>0$, one has
\be \label{Del1er}
\PP \lfi \Del_1 -  0.5\, \| \bCs \bWs^T \bds - \bC \bW^T \bd_{(J)} \|^2 \leq 4 t  \rfi \geq 1 - 2 e^{-t}.
\ee
Note that  
\begin{align*}
\|\hbC \hbW^T \hbd_{(\hJ)} - \bCs \bWs^T \bds\|^2 & = 
\|\hPhJ\, \bxi\|^2 + \| (\bI_{NL} - \hPhJ)  \bCs \bWs^T \bds\|^2 \\
& \geq  \| (\bI_{NL} - \hPhJ)  \bCs \bWs^T \bds\|^2.
\end{align*}
Therefore,  applying an union bound over $m=1, \cdots, n$, $\bC \in \calC(m,n,L)$ 
and $J$ with $|J| = 1, \cdots, ML,$ we derive that for    any $x >0$ 
\begin{align*} 
&  \PP \lfi \Del_2 -  0.5\, \|\hbC \hbW^T \hbd_{(\hJ)} - \bCs \bWs^T \bds\|^2 < 4x \rfi   \\  
 \geq & \PP \lfi \Del_2 -    0.5\,\| (\bI_{NL} - \hPhJ)  \bCs \bWs^T \bds\|^2 < 4x \rfi    \\
 \geq & 1 - \sum_{m=1}^n\ \sum_{\bC \in \calC(m,n,L)}\ \sum_{j =1}^{ML} \ \sum_{|J|=j} 
\PP \lfi 2 |\bxi^T (\bI_{NL} - \hPhJ) \bC  \bW^T \bds|   \right. \\
- & \left.  0.5\,\| (\bI_{NL} - \hPhJ)  \bC  \bW^T \bds\|^2 < 4x \rfi. 
\end{align*}
Denote
\be \label{RDel}
R(m,J,L) =   \log(|\calC(m,n,L)|) + |J| \log\lkr  m^2 L e\,|J|^{-1} \rkr +  2 \log(m |J|). 
\ee
Then, taking into account that the number of sets $J$ with $|J|=j$ is
$$
{ ML \choose j} \leq \lkr \frac{M L e}{j}\rkr^j \leq \lkr \frac{m^2 L e}{j}\rkr^j
$$
and applying  Lemma \ref{lem:bern_err} with $\al =1/2$ and $x = t + R(m,J,L)$, derive
\bes
\PP \lfi \Del_2 - 0.5\, \|\hbC \hbW^T \hbd_{(\hJ)} - \bCs \bWs^T \bds\|^2 - 4 R(\hm,\hJ,L)
\leq 4t  \rfi \geq 1 -  2 e^{-t} \, \sum_{m=1}^n\  \sum_{j =1}^{ML} m^{-2} j^{-2}.
\ees
Since $\sum_{j=1}^\infty  j^{-2} = \pi^2/6 < \sqrt{3}$, the last inequality yields
\be \label{Del2er}
\PP \lfi \Del_2 - 0.5\, \|\hbC \hbW^T \hbd_{(\hJ)} - \bCs \bWs^T \bds\|^2 - 4 R(\hm,\hJ,L)
\leq 4t  \rfi \geq 1 - 6 e^{-t}.
\ee
Finally, in order to obtain an upper bound for $\Del_3$, apply Lemma \ref{corr:quadr_err}
with $\bA =  \hPhJ$  and again use the union upper bound over $m=1, \cdots, n$, $\bC \in \calC(m,n,L)$ 
and $J$ with $|J| = 1, \cdots, ML$ similarly to the way it was done for $\Del_2$. 
Since for any projection matrix $\PJ$, one has $\|\PJ\|_{op}=1$ and
$\| \PJ\|^2  = |J|$, obtain that for any $t>0$  
\be \label{Del3er}
\PP \lfi \Del_3 - |\hJ| - \frac{3}{2} R(\hm,\hJ,L)
\leq \frac{3t}{2} \rfi \geq 1 -  e^{-t},
\ee
where $R(m,J,L)$ is defined in \fr{RDel}.
Combining \fr{main_ineq1}--\fr{Del3er} and recalling that $\hbte = \hbC \hbW^T \hbd_{(\hJ)}$ and $\btes = \bCs \bWs^T \bds$, 
obtain  that, with probability at least $1 - 9 e^{-t}$, one has
\be \label{in_eq}
% \PP \lfi \|\hbte - \btes\|^2 \leq   \min_{\stackrel{m,J,\bq}{\bC \in \calC(m,n,L)}} 
% \lkv  3\, \| \bC \bW^T \bd_{(J)} - \btes \|^2  +  2\, \Pen(|J|,m)  \rkv + 19\,  t \rfi \geq 1 - 9 e^{-t} 
% \PP \lfi \|\hbte - \btes\|^2 \leq   \min_{\stackrel{m,J,\bq}{\bC \in \calC(m,n,L)}} 
% \lkv  3\, \| \bC \bW^T \bd_{(J)} - \btes \|^2  +   11  R(m,J,L)  + 2|J|  \rkv + 19\,  t \rfi \geq 1 - 9 e^{-t}. 
\|\hbte - \btes\|^2 \leq   \min_{\stackrel{m,J,\bq}{\bC \in \calC(m,n,L)}} 
\lfi  3  \| \bC \bW^T \bd_{(J)} - \btes \|^2  +   11  R(m,J,L)  + 2|J|  \rfi + 19t. 
\ee 
 In order to complete the proof of \fr{oracle_prob}, observe that 
$2 \log(|J|) \leq 2|J|$ and $\log(|\calC(m,n,L)|\geq  2 \log m$ by \fr{clust_assump}.
Therefore, one has
\be \label{pen_cond}
 11  R(m,J,L)  + 2|J| \leq 2 \Pen(|J|,m),
 \ee
and \fr{oracle_prob} follows from \fr{in_eq}, \fr{pen_cond} and the second inequality in \fr{symrel}.
 
Finally, inequality \fr{oracle_expec} can  be proved by noting that for any random variable $\zeta$ one has
$\EE \zeta \leq \int_0^\infty \PP(\zeta >z) dz$ and using it with $\zeta = \|\hbLam - \bLams\|^2$.

 %%%%%%%%%%%%%%%%%%%%%%%%%%%%%%%%%%%%%%%%%%%%%%%%%%%%%%%%%%%%%%%%%%%%%%%%%%%%%%%%%%%%%%%%%%%%%%%%%%%%%%%%%%%%%%%%%%%%
%%%%%%%%%%%%%%%%%%%%%%%%%%%%%%%%%%%%%%%%%%%%%%%%%%%%%%%%%%%%%%%%%%%%%%%%%%%%%%%%%%%%%%%%%%%%%%%%%%%%%%%%%%%%%%%%%%%%

\subsection{ Proofs of Lemma~\ref{lem:balanced},   Theorem~\ref{th:upper_smooth_DSBM} and Corollary~\ref{cor:upper_smooth_DSBM}}
This section contains proofs of Lemma~\ref{lem:balanced},   Theorem~\ref{th:upper_smooth_DSBM} 
and Corollary~\ref{cor:upper_smooth_DSBM}.
\\

%%%%%%%%%%%%%%%%%%%%%%%%%%%%%%%%%%%%%%%%%%%%%%%%%%%%%%%%%%%%%%%%%%%%%%%%%%%%%%%%%%%%%%%%%%%%%%%%%%%%%%%5

\noindent
{\bf Proof of Lemma \ref{lem:balanced}.\ }
% The upper bound in \fr{eq:balanced} is trivial. 
% To prove the lower bound in \fr{eq:balanced},
Since
$|\calZ_{\bal} (m,n,n_0,L,1,1)| \leq |\calZ_{\bal} (m,n,n_0,L,\aleph_1,\aleph_2)|$,
 it is sufficient to prove \fr{card_balanced} for $\aleph_1 = \aleph_2=1$. 
Note that 
\bes 
\log|\calZ_{\bal} (m,n,n_0,L,1,1)| \geq  \log(n!) - m \log [(n/m)!] + (L-1) \log{n \choose n_0}
\ees
since there are ${n \choose n_0}$ ways to select $n_0$ nodes out of $n$ but there are more than 
one way to put them back. Applying Lemma \ref{lem:cardinality} with $\ga =1$ obtain that
$\log(n!) - m \log [(n/m)!] \geq n \log(m)/4$. In addition,
\bes
\log{n \choose n_0} > n_0 \log \lkr \frac{n}{n_0}\rkr = \frac{n_0}{4} \log \lkr \frac{n^4}{n_0^4}\rkr
\geq \frac{n_0}{4} \log \lkr \frac{mne}{n_0}\rkr
\ees
provided $n^4/n_0^4 \geq (mne)/n_0$ which holds under conditions of the lemma due to $m \leq n$.
\\

%%%%%%%%%%%%%%%%%%%%%%%%%%%%%%%%%%%%%%%%%%%%%%%%%%%%%%%%%%%%%%%%%%%%%%%%%%%%%%%%%%%%%%%%%%%%%%%%%%%%%%%5

\noindent
{\bf Proof of Theorem~\ref{th:upper_smooth_DSBM}. }
Let  $\ms$ be the true number of classes,  $\Ms = \ms(\ms +1)/2$.
Let $\bC = \bCs$ be the true clustering matrix and $(\bCs)^T \bCs = (\bSs)^2$ where $(\bSs)^2$ is the diagonal 
matrix with the number of nodes in respective pairs of classes on the diagonal. 
Let $\bQs$ be the true matrix of probabilities of connections for pairs of classes, $\bDs = \bQs \bH$, $\bds = \vect(\bDs)$,
$\btes = \bCs (\bWs)^T \bds$ and $\bWs = \bH \otimes \bI_{\Ms}$. 
We need to find an upper bound for  $\| \bCs  (\bWs)^T \bds_{(J)} - \btes \|^2$ in \fr{oracle_prob}. 
Let $\bQs_{(J)}$ be such that 
\bes
 \vect(\bQs_{(J)}) =  \bqs_{(J)} = (\bWs)^T \bds_{(J)} =  (\bWs)^T (\bWs \bqs)_{(J)}. %, \bU = \bC^T \bS^{-1}.
\ees
Then, by direct calculations, one obtains  
\begin{align*}
\| \bCs  (\bWs)^T \bds_{(J)} - \btes \|^2 & = \lkv \bqs_{(J)}  - (\bSs)^{-2} \bCs \btes\rkv^T  (\bSs)^2  \lkv \bqs_{(J)} - (\bSs)^{-2} \bCs \btes\rkv \\
& +  \|\btes\|^2 - \|(\bSs)^{-1} (\bCs)^T \btes \|^2.
\end{align*}
Since  $\|(\bSs)^{-1} (\bCs)^T \btes\| = \|\btes\|$ and $(\bSs)^{-2} \bCs \btes = \bqs = \vect(\bQs) = \vect(\bDs \bH)$,
% where $\bQs$ and $\bDs$  are the true values of matrices $\bQ$ and $\bD$. 
% Moreover, $\tilbq_J  = \bW^T \bd_{(J)} = \vect(\bQ_J)$.
obtain 
\begin{align*}
\| \bCs  (\bWs)^T \bds_{(J)} - \btes \|^2 & =  \|\bSs (\bQs_{(J)} - \bQs)\|^2 \\
& \leq  \sum_{k=1}^{\Ms} \max_{l}  N_{k}^{(l)}\ \sum_{l=1}^L \lkv (\bQs_{(J)})_{k,l} - \bQs_{k,l} \rkv^2\\
& \leq \aleph_2^2 \lkr \frac{n}{\ms} \rkr^2 \, \sum_{k=1}^{\Ms} \|\bQs_{(J)})_{k,*} - \bQs_{k,*}\|^2 
\end{align*}
Here, $k$ is the index corresponding to a pair of classes $(k_1, k_2)$,   $N_{k_1, k_2}^{(l)}$ is 
defined in formula \fr{eq:Nkl} and the second inequality follows from assumption \fr{balanced}.
In order to complete the proof, note that 
\begin{align*}  
% \sum_{k=1}^{\Ms} \|\bQ_{(J)})_{k,*} - \bQs_{k,*}\|^2 & = 
\sum_{k=1}^{\Ms} \|\bQs_{(J)})_{k,*} - \bQs_{k,*}\|^2 & = \sum_{k=1}^{\Ms}   \|\bDs_{(J_k)})_{k,*} - \bDs_{k,*}\|^2 
% \\ &   
= \sum_{k=1}^{\Ms}   \sum_{l \notin J_k} (\bDs_{k,l})^2.   
\end{align*} 
Therefore, \fr{oracle_prob} implies \fr{oracle_prob_specific}. 
\\

%%%%%%%%%%%%%%%%%%%%%%%%%%%%%%%%%%%%%%%%%%%%%%%%%%%%%%%%%%%%%%%%%%%%

\noindent
{\bf Proof of Corollary~\ref{cor:upper_smooth_DSBM}.}
Observe that it follows from \fr{bias_coef_cond} that one can choose 
$J_k = \{l:\ 1 \leq l \leq L_0\}$ where $L_0 \leq L$. Then, 
\fr{eq:J_union} yields that $|J| = L_0 \Ms$.
Moreover, due to assumption \fr{bias_coef_cond}, obtain 
\bes 
\sum_{k=1}^{\Ms}\, \sum_{l \notin J_k} (\bDs_{k,l})^2 \leq K_0 \Ms L_0^{-2 \nu_0}.
\ees
Note also that if $L_0 = L$, then there is no bias and the sum in \fr{oracle_prob_specific} is identical zero. 
Then, \fr{oracle_prob_specific} becomes
\begin{align*}   
 \|\hbLam - \bLams\|^2  & \leq \min \lfi \Del(n,L,\ms),\ 22 \Ms L   \rfi\\
& + 44 \lkv n \log (\ms)  + 
   n_0 L  \,\log\lkr \frac{\ms ne}{n_0}\rkr \rkv + 38\, t,
\end{align*}
where
\bes % \label{eq:Del_expression}
\Del(n,L,\ms) = \min_{0 \leq L_0 <L} \lfi 6  {K}_0  \aleph_2^2  n^2\, L_0^{-2 \nu_0} 
+ 22 (\ms)^2 L_0 \log\lkr \frac{25\,L}{L_0}\rkr \rfi
\ees 
In order to obtain \fr{oracle_prob_smooth}, minimize the right-hand side of the last expression 
with respect to $L_0<L$ and note that, if $L_0 <L$, then
$\log(L_0) \asymp \log (L) \asymp \log(n/\ms)$.

%%%%%%%%%%%%%%%%%%%%%%%%%%%%%%%%%%%%%%%%%%%%%%%%%%%%%%%%%%%%%%%%%%%%%%%%%%%%%%%%%%%%%%%%%%%%%%%%%%%%%%%%%%%%%%%%%%%%
%%%%%%%%%%%%%%%%%%%%%%%%%%%%%%%%%%%%%%%%%%%%%%%%%%%%%%%%%%%%%%%%%%%%%%%%%%%%%%%%%%%%%%%%%%%%%%%%%%%%%%%%%%%%%%%%%%%%

\subsection{ Proofs of Theorems~\ref{th:DSBM_lower_bounds} and \ref{th:smooth_DSBM_lower_bounds}}
%\label{sec:lower}

This section contains the proofs of the lower bounds for the error.
The lower bounds in both, Theorem \ref{th:DSBM_lower_bounds} and Theorem~\ref{th:smooth_DSBM_lower_bounds},
consist  of two parts, the clustering error and the nonparametric estimation error.
We shall consider those terms separately.
\\  

%%%%%%%%%%%%%%%%%%%%%%%%%%%%%%%%%%%%%%%%%%%%%%%%%%%%%%%%%%%%%%%%%%%%%%%%%%%%%%%%%%%%%%%%%%%%%%%%%%%%%%%%%%%%%%%%%%%%

\noindent
{\bf Proof of Theorem \ref{th:DSBM_lower_bounds}. }
Although the upper bounds for the risk in Corollary~\ref{cor:upper_DSBM}
are derived for the case of general clustering matrices, due to the fact that the balanced model clustering complexity
is the same as complexity of general clustering, we derive the lower bounds for the clustering error for the case 
when $\bC \in  \calZ_{\bal} (m,n,n_0,L,\aleph_1, \aleph_2)$. Moreover, since the case of $\aleph_1 = \aleph_2 =1$
is the most restrictive, we prove the clustering error for this case.
\\

\underline{\bf The clustering error. } 
Without loss of generality, assume that   $\ga m$ and $\ga n$ are integers.  
Assume that connectivity tensor $\bG$ does not change with $l$, so   $\bG_{*,*,l} = \bV$ is an $m\times m$ symmetric matrix.
Let $\bV$ be block diagonal and  such that the diagonal blocks are equal to zero
and the non-diagonal blocks are equal to $\bF$ and $\bF^T$, respectively, so that
$\bV_{k_1,k_2} =0$ if $1 \leq k_1,k_2 \leq (1-\ga)m$ or 
$(1-\ga)m+1 \leq k_1,k_2 \leq m$ and  
$\bV_{k_1,(1-\ga)m +k_2}  = \bF_{k_1,k_2}$ if  $k_1 = 1 ,\cdots (1-\ga)m$, $k_2= (1-\ga)m+1, \cdots,   m$. 
Since components of vectors $\bG_{k_1, k_2, *}$ are constant for any $k_1$, $k_2$, then,
due to condition \fr{H_assump}, each of the vectors $\bH \bF_{k_1, k_2}$ has only one non-zero component, 
so that  the set $J$ has at most $\ga(1-\ga)m^2 < s$ nonzero elements.

Consider a collection of binary vectors $\bom \in \{0,1\}^{(1-\ga)m}$.
By Varshamov-Gilbert Lemma (see Tsybakov (2009)), there exists a subset $\Xi$ of those vectors 
such that for any $\bom, \bom' \in \Xi$ one has 
$\|\bom - \bom'\|_H = \|\bom - \bom'\|^2 \geq (1 - \ga)m/8 \geq m/16$
and $|\Xi| \geq \exp((1-\ga)m/8)$. 
Assume, without loss of generality, that $m$ is large enough, so that 
$\exp((1-\ga)m/8) \geq \ga m$, otherwise, choose a smaller value of $\ga$ 
(inequality $\exp((1-\ga)m/8) \geq \ga m$ is always valid for $\ga \leq 1/9$).
Choose $\ga m$ vectors $\bom$ in $\Xi$, enumerate them as $\bom^{(1)}, \cdots, \bom^{(\ga m)}$
and use them to form columns of matrix $\bF$ as follows:
\be \label{Fmatr}
\bF_{*,j} = 0.5\, \bone + \al \bom^{(j)}, \quad j=1, \cdots, \ga m.
\ee
Then, for any $j, j' = 1, \cdots, \ga m$, obtain
\be \label{dist1}
\|\bF_{*,j} - \bF_{*,j'}\|^2 \geq \al^2   m/16,
\ee
where $\al$ is a positive constant that will be defined later.
Note that for every $l$ and $k$ one has 
\be \label{inequalities}
l\, \log\lkr \frac{k}{l}\rkr \leq \log {k \choose l} \leq l\, \log\lkr \frac{k e}{l}\rkr, \quad 
\log(k!) = k \log k - k + \frac{1}{2} \log(2\pi k) + o(1),
\ee
where the $o(1)$ term is smaller than 1. Therefore, it follows from the first formula in \fr{card_clust} that
\be \label{log_card}
n \log m + n_0(L-1)\, \log \lkr \frac{m n}{n_0} \rkr \leq 
\log|\calZ (m,n,n_0,L)| \leq  
n \log m + n_0(L-1)\, \log \lkr \frac{m n e}{n_0} \rkr  
\ee
The term $n \log m$ in \fr{log_card} is due to the initial clustering while the term 
$n_0(L-1)\, \log \lkr  m n/n_0  \rkr$ is due to temporal changes in the clusters' memberships.

In what follows, we shall utilize clustering functions $z^{(l)}: [n] \to [m]$
corresponding to clustering matrices $\bC^{(l)}$  such that  $z^{(l)} (j) = k$ 
iff at the moment $t_l$ node $j$ belongs to class $\Om_k$, $k=1, \cdots, m$. 
% Denote a set of such clustering matrices at time $t_l$ by $\calZ_l (m,n,n_0)$ and the 
% overall collection   by $\calZ (m,n,n_0,L)$.   Note that $\calZ_1 (m,n,n_0)$ does not depend on   $n_0$,
% so that   $\calZ_1 (m,n,n_0)  \equiv  \calZ_1 (m,n)$.
 \\

%%%%%%%%%%%%%%%%%%%%%%%%%%%%%%%%%%%%%%%%%%%%%%%%%%%%%%%%%%%%%%%%%%%%%%%%%%%%%%%%%%%%%%%%%%%%%%%

{\bf Clustering error due to initial clustering. } 
First consider the case when initial clustering error dominates. 
If $m=2$ or $m$ takes a small value, the proof is almost identical to the proof in Section 3.3 of  Gao \etal (2015).
Hence, we shall skip this part and consider the case when $m$ is large enough,
so that $\ga m \geq 2$.
 
Following  Gao \etal  (2015), we consider   clustering matrices and clustering functions 
independent of $l$, so that $z^{(l)}\equiv z$.
% \bes
% z^{(l)}\equiv z, \quad \calZ (m,n) = \calZ_1 (m,n)  
% \ees 
Consider a sub-collection of clustering matrices $\calF(m,n,\ga) \subset \calM(m,n)$
such that they cluster the first $n(1-\ga)$ nodes into the first $m(1-\ga)$ classes uniformly and sequentially,
$n/m$ nodes in each class, i.e., the first $n/m$ nodes are placed into class  $\Om_1$, 
the second $n/m$ nodes into class  $\Om_2$, and so on.
The remaining $\ga n$ nodes are clustered into the remaining $\ga m$ classes, $n/m$  nodes into each class.
Then, by Lemma \ref{lem:cardinality}, 
$$
\log |\calF(m,n,\ga)| = \log \lkr (\ga n)!\Big/ \lkv (n/m)!\rkv^{\ga m} \rkr \geq  \ga n  \log(\ga m)/4.
$$  
Now, apply Lemma~\ref{lem:packing} with $\ga n$ and $\ga m$, respectively, instead of $n$ and $m$ and 
$r = \ga n/32$. Derive  that there exists a 
subset $\calS(m,n,\ga)$ of the set $\calF(m,n,\ga)$ such that, for any $\bC, \bC' \in \calS(m,n,\ga)$, one has 
$2 \lfi \#j:\  z(j) \neq z'(j)\rfi = \| \bC  - \bC' \|_H   \geq \ga n/32$. Also,  by \fr{AA3}, 
\be    \label{cardSnm}
\log |\calS(m,n,\ga)| \geq \frac{\ga n}{4} \log (\ga m)  -   \frac{\ga n\,  \log (32 m \ga  e)}{32}  \geq
 \frac{\ga n}{16} \log (\ga m).
\ee 
Let $\bLam$ and $\bLam'$ be the tensors of probabilities 
corresponding to,  respectively, clustering matrices 
$\bC, \bC' \in \calS(m,n,\ga)$ with related clustering functions
$z$ and $z'$. Then, by \fr{dist1},
due to the fact that the first $n(1-\ga)$ nodes are clustered uniformly and 
sequentially,   obtain
\begin{align*}
& \| \bLam - \bLam'\|^2 
  =   2 L \, \sum_{i=1}^{n(1-\ga)}\, \sum_{j = n(1-\ga)+1}^n (\bF_{z(i),z(j)} -  \bF_{z'(i),z'(j)})^2  \\
&  =    \frac{2 L n}{m} \, \sum_{k=1}^{m(1-\ga)}\, \sum_{j = n(1-\ga) + 1}^n (\bF_{k,z(j)} -  \bF_{k,z'(j)})^2 \\
& =     \frac{2 L n}{m} \,   \sum_{j = n(1-\ga) + 1}^n \|\bF_{*,z(j)} -  \bF_{*,z'(j)}\|^2  
  \geq  \frac{2 L n}{m} \,   \frac{\al^2\, m}{16} \ \lfi \#j:\  z(j) \neq z'(j)\rfi,   
\end{align*}
so that 
\be \label{d1}
\| \bLam - \bLam'\|^2  \geq 2^{-9}\, L n^2 \al^2 \ga.
\ee
On the other hand, if $\al \leq 1/4$, then, by Lemma~\ref{lem:Gao2015},
obtain that the Kullback divergence is bounded above 
\be \label{kul1}
K (\PP_{\bLam}, \PP_{\bLam'}) \leq 8\, \| \bLam - \bLam'\|^2  \leq
16 \al^2 n^2 L  \ga.
\ee 
Set $\al^2 = C_\al \, \log (\ga m)/nL$ % where $C_\al$ is an absolute constant
 and apply Theorem 2.5 of Tsybakov (2009).
Due to \fr{cardSnm} and \fr{kul1}, if $C_\al$ is a small enough absolute constant,  
$$
16 \al^2 n^2 L \ga = 16 C_\al \, \log (\ga m) n < (1/8)\, (\ga n/16)\, \log (\ga m) \leq (1/8)\,  \log |\calS(m,n,\ga)| 
$$
and conditions of   Theorem 2.5 are satisfied. 
Since  
$L n^2 \al^2 \ga = C_\al  n\,\ga \log (\ga m)$ and $\log (\ga m) \geq C(\ga) \log m$ for some constant $C(\ga)$ dependent on $\ga$ only, 
derive  
\be \label{lower_DSBM1}
\inf_{\hbLam} \sup_{\stackrel{\bG \in \calG_{m,L,s}}{\bC \in \calS(m,n,\ga)}}
\PP_{\bLam} \lfi \frac{\|\hbLam - \bLams\|^2}{n^2\,L} \geq   \frac{C(\ga)\, \log m}{n L}
\rfi \geq 1/4.
\ee

%%%%%%%%%%%%%%%%%%%%%%%%%%%%%%%%%%%%%%%%%%%%%%%%%%%%%%%%%%%%%%%%%%%%%%%%%%%%%%%%%%%%%%%%%%%%%%%%%%%%%%%%%%%%%%%%%%%%

{\bf Clustering error due to  changes in the memberships. }
Now, we consider the case when the clustering error which is due to the temporal changes in memberships dominates
the error of initial clustering. Use the same construction for $\bG$ and $\bF$ as before. 
Consider the following collection of clustering matrices $\calF = \prod_{l=1}^L \calF_l$ where $\calF_l$ is defined as follows.
When $l$ is odd, $\calF_l$ contains only one matrix that clusters nodes uniformly and sequentially, i.e., 
the first $n/m$ nodes go to class $\Om_1$, the second  $n/m$ nodes go to class $\Om_2$
and the last $n/m$ nodes go to class $\Om_m$. If $l$ is even, $\calF_l =  \calP(m,n,n_0,\ga)$ where $\calP(m,n,n_0,\ga)$ 
is the set of clustering matrices  that corresponds to a perturbation of the uniform sequential 
clustering with at most $n_0$ nodes moved to different classes 
in the manner described below. Let $k_0$ be an integer such that
\be \label{k0}
k_0 \leq n_0/(\ga m) < k_0 +1.
\ee
If $k_0 =0$, then $n_0 < \ga m$ and we choose $n_0$ clusters out of the last $\ga m$ clusters,
remove one element from each of those clusters and then put those $n_0$ elements back  
in such a manner that every element goes to a different cluster and 
no elements goes back to its own cluster. If $k_0 \geq 1$, we remove $k_0$ elements from each of the 
last $\ga m$ clusters and then put each of those $k_0$-tuples back, one tuple per cluster, 
so that none of the tuple goes back to its own cluster.
Then, $\log |\calF| = [L/2] \log|\calP(m,n,n_0, \ga)|$ where $[L/2] \geq (L-1)/2$ is the largest integer not exceeding $L/2$ and 
\bes % \label{F_card}
\log|\calP(m,n,n_0, \ga)| = \lfi
\begin{array}{ll}
\log {\ga m \choose n_0} + n_0 \log(n/m) + \log[(n_0-1)!], & \mbox{if}\quad k_0 =0;\\
\ga m \log {n/m \choose k_0} + \log[(\ga m -1)!], & \mbox{if}\quad k_0 \geq 1.
\end{array} \right.
\ees
If $n_0 < \ga m$, so that $k_0 =0$,  then, by \fr{inequalities}, obtain that
$\log|\calP(m,n,n_0, \ga)| \geq n_0 \log(\ga m/n_0) + n_0 \log(n/m) = n_0 \log(\ga n/n_0)$.
If $n_0 \geq \ga m$ and $k_0 \geq 1$, then,  
by \fr{inequalities}, obtain
$\log|\calP(m,n,n_0, \ga)| \geq \ga   m k_0\,  \log(n/(m k_0))$.
Since $k_0 +1 \leq 2 k_0$, obtain that $k_0 \geq n_0/(2 m \ga)$. Hence, for any $k_0 \geq 0$ 
\be \label{calF_card}
\log |\calP(m,n,n_0, \ga)| \geq \frac{n_0}{2} \, \log \lkr \ga n/n_0  \rkr.
\ee
For every even value of $l$, apply Lemma \ref{lem:packing} with $\ga n$ and $\ga m$, 
respectively, instead of $n$ and $m$ and $r = n_0/40$ obtaining that there exists a 
subset $\calS_l (m,n,n_0,\ga)$ of the set $\calP(m,n,n_0, \ga)$ such that, for any 
$\bC^{(l)}, \bC'^{(l)} \in  \calS_l (m,n,n_0,\ga)$, one has 
\be \label{d_l}
  \| \bC^{(l)}  - \bC'^{(l)} \|_H   \geq  n_0/40, \quad  l=2k,\ k=1, \cdots, [L/2].
\ee 
By \fr{calF_card} and Lemma~\ref{lem:packing}, for every even $l$, one has 
\begin{align*}
\log |\calS_l (m,n,n_0,\ga)| \geq    \frac{n_0}{2}\, \log\lkr \frac{\ga n}{n_0}\rkr -   
\frac{n_0}{40}\, \log \lkr \frac{80 n e m \ga^2}{n_0}\rkr \geq \frac{n_0}{40}\,  \log \lkr \frac{n e m}{n_0}\rkr
\end{align*}
since, due to $n_0 \leq   4/3\, \ga n\,m^{-1/9}$ and $(80 e^2)^{1/18} \leq 0.75$, one has 
$(\ga n/n_0)^{20} \ge (80 e^2)\, [m\, (\ga n /n_0)]^2$, so that 
$$
20\, \log(\ga n/n_0) \ge  \log(80 n^2 m^2 e^2 \ga^2/n_0^2) = \log(n e m/n_0) +
 \log (80 n e m \ga^2/n_0). 
$$   
% if $\ga n/n_0  \geq m^{1/9} (80 e^2)^{1/18} \leq 0.75\, m^{1/10}$. 
For odd values of $l$, let $\calS_l (m,n,n_0,\ga)$ contain just one clustering matrix corresponding to the 
uniform sequential clustering. 
Now, consider the set $\calS  (m,n,n_0,\ga,L) = \prod_{l=1}^{L} \calS_l (m,n,n_0,\ga)$ with   
\be  \label{cardSnmL}
% \calS  (m,n,n_0,\ga,L) = \prod_{l=1}^{L} \calS_l (m,n,n_0,\ga) \ \mbox{with} \ 
\log|\calS  (m,n,n_0,\ga,L)| \geq  \frac{(L-1)n_0}{80}\, \log \lkr \frac{n e m}{n_0}\rkr.
\ee 
Let $\bC = (\bC^{(1)}, \bC^{(2)},\cdots,\bC^{(L)})$ and $\bC' = (\bC'^{(1)}, \bC'^{(2)},\cdots,\bC'^{(L)})$  
be two sets of clustering matrices with $\bC^{(l)}, \bC'^{(l)} \in  \calS_l (m,n,n_0,\ga)$ and let $z = (z_1, \cdots, z_L)$ and 
 $z' = (z'_1, \cdots, z'_L)$  be the corresponding clustering functions.
Let $\bLam$ and $\bLam'$ be the tensors of probabilities 
corresponding to sets of clustering matrices $\bC, \bC' \in \calS  (m,n,n_0,\ga,L)$. 
% and related clustering functions $z$ and $z'$. 
%
Then, similarly to the previous case, using \fr{dist1}, derive
\begin{align*}
& \| \bLam - \bLam'\|^2 
  =   2 \sum_{l=1}^{[L/2]} \, \sum_{i=1}^{n(1-\ga)}\, \sum_{j = n(1-\ga)+1}^n (\bF_{z_{2l}(i),z_{2l}(j)} -  \bF_{z'_{2l}(i),z'_{2l}(j)})^2  \\
& =    \frac{2n}{m} \, \sum_{l=1}^{[L/2]} \, \sum_{k=1}^{m(1-\ga)}\, \sum_{j = n(1-\ga)+1}^n (\bF_{k,z_{2l}(j)} -  \bF_{k,z'_{2l}(j)})^2 \\
& =       \frac{2n}{m} \, \sum_{l=1}^{[L/2]} \,  \sum_{j = n(1-\ga)+1}^n \|\bF_{*,z_{2l}(j)} -  \bF_{*,z'_{2l} (j)}\|^2  
  \geq    \frac{n \al^2}{8} \, \sum_{l=1}^{[L/2]}   \|\bC^{(2l)} - \bC'^{(2l)}\|_H, 
\end{align*}
so that  by \fr{d_l}, 
\be \label{d11}
\| \bLam - \bLam'\|^2  \geq  (L-1) n n_0 \al^2 /1280.
\ee
Again, similarly to the previous case, if $\al \leq 1/4$, then by Lemma~\ref{lem:Gao2015},
obtain that the Kullback divergence is bounded above 
\be \label{kul2}
K (\PP_{\bLam}, \PP_{\bLam'}) \leq 8\, \| \bLam - \bLam'\|^2  \leq
8  L n n_0 \al^2.
\ee 
Set 
%\bes
$\al^2 = C_\al n^{-1}\, \log (n e m/n_0)$
%\ees
 where $C_\al$ is an absolute constant  and apply Theorem 2.5 of Tsybakov (2009).
Observe that if $C_\al$ is small enough, then, due   to \fr{cardSnmL} and \fr{kul2},  conditions of this theorem are satisfied, hence,  
\be \label{lower_DSBM2}
\inf_{\hbLam} \sup_{\stackrel{\bG \in \calG_{m,L,s}}{\bC \in \calS  (m,n,n_0,\ga,L)}}
\PP_{\bLam} \lfi \frac{\|\hbLam - \bLams\|^2}{n^2\,L} \geq  C(\ga) \,   \frac{n_0}{n^2}\,\log \lkr \frac{n e m}{n_0} \rkr
 \rfi \geq 1/4.
\ee

%%%%%%%%%%%%%%%%%%%%%%%%%%%%%%%%%%%%%%%%%%%%%%%%%%%%%%%%%%%%%%%%%%%%%%%%%%%%%%%%%%%%%%%%%%%%%%%%%%%%%%%%%%%%%%%%%%%%

\underline{\bf The nonparametric estimation error. } 
Consider uniform sequential clustering with $n/m$ nodes in each group and 
group memberships remaining the same for all $l=1, \cdots, L$. 
% Let tensor $\bG$ be such that $\bG_{k,k,l} =0$ for any $k = 1, \cdots, m$, $l=1, \cdots, L$. 
Let $\bQ  \in \RR^{M\times L}$ be the matrix with columns   $\bq^{(l)}$,
$l=1, \cdots, L$, defined in Section \ref{sec:vectorization}.
Denote $\bV = \bQ \bH^T \in \RR^{M\times L}$ and recall that for $\bG \in \calG_{m,L,s}$,   by   \fr{vecQHT}
and \fr{transformed}, matrix $\bV$  should have at most $s$ nonzero entries.

Let $k_0 = \min(s/2, M)$.  Choose $k_0$ rows among $M$ rows of matrix $\bV$
and denote this set by $\calX$. If $k_0 = M$,  set  $\calX  = \{1, \cdots, M\}$.
For $k \in \calX$, set $\bV_{k,1} \neq 0$. 
We have already distributed $k_0$ non-zero entries and have $s -k_0$ entries left.
We distribute those entries into  the $k_0$ rows $\bV_{k,*}$ where $k \in \calX$.
Let   
\be \label{k0s0} 
s_0 = [(s -k_0)/k_0] = [s/k_0] -1 \quad \mbox{with} \quad s/2 \leq k_0 s_0 <s,
\ee 
where $[s/k_0]$ is the largest integer no larger than $s/k_0$.
Consider a set of binary vectors $\bom \in \{0,1\}^{L}$ with exactly $s_0$ ones in each 
vector. By Lemma 4.10 of Massart (2007), there exists a subset $\calT$ of those vectors  such that 
for any $\bom, \bom' \in \calT$, one has
\bes
\|\bom - \bom'\|_H \geq s_0/2  \quad \mbox{and} \quad 
\log|\calT| \geq 0.233\, s_0\, \log(L/s_0).
\ees 
Denote $\tilde{\calT} = \bigotimes_{k \in \calX} \calT_k$, where $\calT_k$ is a copy of the set $\calT$
corresponding to row $k$ of matrix $\bV$.
For $\bom^{(k)} \in \calT_k$, set 
\be \label{bV}
\bV_{k,*} = (\sqrt{L}/2, \cdots, 0) + \al m/n \,  \bom^{(k)}, \ \mbox{if}\ k \in \calX,
\quad  \bV_{k,*} = 0 \ \mbox{if}\  k \notin \calX.
\ee
It is easy to see that matrix $\bV$ has at most $s$ nonzero entries as required.

Let $\bV$ and $\bV'$ be matrices corresponding to sequences $\bom^{(k)}$ and $\bom'^{(k)}$  
in $\calT_k$, $k \in \calX$. Let $\bLam$ and $\bLam'$ be the tensors corresponding to $\bV$ and $\bV'$.
% It is pretty straightforward to check that $s/2 \leq k_0 s_0 \leq s$.   
Then, due to \fr{k0s0} and the uniform sequential clustering,   \fr{bV} implies that 
\begin{align*} 
\|\bLam - \bLam'\|^2 
& \geq   \lkr \frac{n}{m}\rkr^2 \|\bV - \bV'\|^2 
\geq   \lkr \frac{n}{m}\rkr^2 \sum_{k \in \calX} \|\bV_{k,*} - \bV'_{k,*}\|^2  
  \geq   \frac{\al^2  k_0 s_0}{2} \geq  \frac{\al^2  s}{4};\\
\|\bLam - \bLam'\|^2 & \leq    4\, \lkr \frac{m}{n}\rkr^2 \al^2 \lkr \frac{n}{m}\rkr^2 k_0 s_0 \leq 4 \al^2 s.
\end{align*}
Set $\al^2 = C_\al \log (L/s_0)$.    It is easy to check that, due  to assumptions \fr{H_assump} and \fr{Js_cond},
one has $\bQ_{ij} \in [1/4,3/4]$ for any $i$ and $j$. Hence,   by Lemma~\ref{lem:Gao2015}, obtain
$K (\PP_{\bLam}, \PP_{\bLam'}) \leq 8\, \| \bLam - \bLam'\|^2  \leq 32 \al^2 s$. 
 If $C_\al \leq 2^{-8} \cdot 0.233$, then $32 \al^2 s < (1/8) \log(\tilde{\calT})$ and 
conditions of Theorem 2.5 of Tsybakov (2009) hold.

Finally, in order to obtain the last term in \fr{lower_DSBM}, examine $L/s_0$.
If $s < 2M$, then $k_0 = s/2$, $s_0=1$ and  $L/s_0 = L = L m^2/m^2 \geq \ga L m^2/s$.
If $s \geq 2M$, then $k_0 = M$, $s_0\leq s/M$ and $L/s_0 \geq LM/s \geq Lm^2/(2s) \geq \ga L m^2/s$. 
Since for some constant $C(\ga)>0$ independent of $L$ and $m$, one has $\log (\ga L m^2/s) \geq C(\ga) \log(L m^2/s)$, obtain
\be \label{lower_DSBM3}
\inf_{\hbLam} \sup_{\stackrel{\bG \in \calG_{m,L,s}}{\bC \in \calS  (m,n,n_0,\ga,L)}}
\PP_{\bLam} \lfi \frac{\|\hbLam - \bLams\|^2}{n^2\,L} \geq   C_\ga \,   \frac{s}{n^2 L}\,\log \lkr \frac{L  m^2}{s} \rkr
 \rfi \geq 1/4.
\ee

Finally, in order to obtain the lower bound in \fr{lower_DSBM} observe that, for any $a,b,c \geq 0$, 
one has $\max(a,b,c)\leq a+b+c \leq 3 \max(a,b,c)$     and then combine
\fr{lower_DSBM1}, \fr{lower_DSBM2} and \fr{lower_DSBM3}.
\\

%%%%%%%%%%%%%%%%%%%%%%%%%%%%%%%%%%%%%%%%%%%%%%%%%%%%%%%%%%%%%%%%%%%%%%%%%%%%%%%%%%%%%%%%%%%%%%%%%%%%%%%%%%%%%%%%%%%%
%%%%%%%%%%%%%%%%%%%%%%%%%%%%%%%%%%%%%%%%%%%%%%%%%%%%%%%%%%%%%%%%%%%%%%%%%%%%%%%%%%%%%%%%%%%%%%%%%%%%%%%%%%%%%%%%%%%%

\noindent
{\bf Proof of Theorem \ref{th:smooth_DSBM_lower_bounds}. }
Note that for Theorem~\ref{th:DSBM_lower_bounds} we proved the lower bounds for the clustering error in the most restrictive case when 
$\bC \in  \calZ_{\bal} (m,n,n_0,L,1,1)$. Moreover, in this proof, the connection probabilities are set to be constant over time,  hence
due to condition  \fr{H_assump}, $\bDs_{k,l}=0$ for $l \geq 2$ in assumption \fr{bias_coef_cond}, so \fr{bias_coef_cond} holds
for any $K_0$ and $\nu_0$. Therefore, the lower bounds for the risk due to clustering errors hold in this case and coincide with the lower bounds in 
Theorem~\ref{th:DSBM_lower_bounds}. For this reason, we only need to prove the lower bounds that are due to the nonparametric estimation error.
\\

%%%%%%%%%%%%%%%%%%%%%%%%%%%%%%%%%%%%%%%%%%%%%%%%%%%%%%%%%%%%%%%%%%%%%%%%%%%%%%%%%%%%%%%%%%%%%%%%%%%%%%%%%%%%%%%%%%%%

\underline{\bf The nonparametric estimation error. } 
Consider a set up where the nodes are grouped into $m$ classes,
$n/m$ nodes in each class, so the model is fully balanced.
Let $\bG_{k,k,l} =0$ for any $k=1, \cdots, m$ and $l=1, \cdots, L$.

Consider an even number $L_0$ such that $1 \leq L_0 \leq L/2$ and a set of vectors $\bom \in \{ 0,1\}^{L_0}$ 
with exactly $L_1 = L_0/2$ nonzero entries.
By Lemma 4.10 of Massart (2007), there exists a subset $\calT$ of those vectors such that 
for any $\bom, \bom' \in \calT$ one has 
\be \label{new_bombom}
\|\bom - \bom'\|_H \geq L_0/4, \quad
\log |\calT| \geq 0.233 (L_0/2) \log 2 \geq 0.08 L_0.
\ee 
Denote $\calK = \{(k_1, k_2):\ 1 \leq k_1 < k_2 \leq m\}$ and 
let  $\calT_{k_1, k_2}$ be the copies of $\calT$ for $(k_1, k_2) \in \calK$. 
Denote $\tilde{\calT} = \bigotimes_{(k_1,k_2) \in \calK}\, \calT_{k_1, k_2}$ and
observe that 
\bes
\log|\tilde{\calT}| = 0.5\, m(m-1)\, \log |\calT| \geq 0.02 m^2\, L_0.
\ees
Then, $\tilde{\bom} \in \tilde{\calT}$ are   binary tensors with elements $\bom_{l}^{(k_1, k_2)}$, $l=1, \cdots, L_0$,
and  $(k_1,k_2) \in \calK$.
Consider a set of matrices  $\bD^{(\tilde{\bom})}$ indexed by  $\tilde{\bom}$
such that for the index $k = 1, \cdots, M,$ corresponding to $(k_1,k_2) \in \calK$, one has
\bes
\bD^{(\tilde{\bom})}_{k,1} = \sqrt{L}/2; \quad 
\bD^{(\tilde{\bom})}_{k, l} =  
\al \, \bom_{l-L_0}^{(k_1, k_2)},\ l=L_0+1, \cdots, 2L_0  , \  k =(k_1,k_2) \in \calK.
% \quad \mbox{if}
% \beta_{k_1-1}  < x \leq \beta_{k_1}, \ \beta_{k_2-1} < y \leq \beta_{k_2}.
\ees  
% where $\tilde{\bom}$ is a binary matrix with elements $\bom_{l}^{(k_1, k_2)}$,$l=1, \cdots, L_0$
% and  $(k_1,k_2) \in \calK$.
%
In order condition \fr{bias_coef_cond} is satisfied, we set
\be \label{new_rhoL0}
\al^2 \leq C_1 L_0^{-(2 \nu_0 +1)}  \quad \mbox{with} \quad C_1 \leq  \min(K_2 2^{1 - 2\nu_0},1/8).
\ee 
Denote by $\bLam$ and $\bLam'$ the probability tensors corresponding, respectively,  to $\tilde{\bom}$ and 
$\tilde{\bom'}$ in $\tilde{\calT}$. Then, due to \fr{new_bombom} and the symmetry, 
\beqns
\|\bLam - \bLam'\|^2 &\geq &    \al^2 \sum_{k_1 =1}^m \sum_{k_2 = k_1 +1}^m 
\|\bom^{(k_1, k_2)} - \bom'^{(k_1, k_2)}\|_H \lkr \frac{n}{m}\rkr^2  \geq \frac{\al^2 n^2 L_0}{8};\\
\|\bLam - \bLam'\|^2 &\leq &    \al^2 m(m-1) (n/m)^2 L_0 \leq \al^2 n^2 L_0.
\eeqns
Note that  one has $1/4 \leq \bQ_{k,l} \leq 3/4$ provided
$\|\bH^T \bom ^{(k_1, k_2)}\|_\infty \leq 1/4$ for $(k_1, k_2) \in \calK$.
By  Assumption \fr{H_assump}, the latter is guaranteed by 
$\al^2 L_0^2/L \leq 1/4$, so that, due to $L  \geq 2 L_0$ and $\nu_0 \geq 1/2$,
it is ensured by \fr{new_rhoL0}. Then, by  Lemma~\ref{lem:Gao2015}, one has 
\bes
K (\PP_{\bLam}, \PP_{\bLam'}) \leq 8\, \| \bLam - \bLam'\|^2  \leq
8  \al^2 n^2 L_0 \leq \log |\tilde{\calT}|/8   
\ees  
provided 
\bes % \label{rho_con1}
\al^2 \leq C_2 (m/n)^2,
\ees
 where $C_2$ is an absolute constant.
Therefore, application of Theorem 2.5 of Tsybakov (2009) yields that 
\bes
\inf_{\hbLam} \sup_{\stackrel{\bG \in \calG_{m,L,s}}{\bC \in \calZ_{\bal}}}
\PP_{\bLam}   \lfi \frac{\|\hbLam - \bLam \|^2}{n^2\,L}   \geq \Del(n,L) \rfi \geq \frac{1}{4} 
\ees
with $\Del(n,m,L) = C\,\al^2 L_0/L$ where $C$ is an absolute constant.

Now, we denote $C_3^2 = (C_2/C_1) 2^{-(2 \nu_0 +1)}$ and consider two cases.   
If $n \leq C_3 m L^{\nu_0 + 1/2}$, choose $L_0 = [(C_1/C_2) (n/m)^2]^{1/(2 \nu_0 +1)}$ 
which leads to $\al^2 = C_1 L_0^{-(2 \nu_0 +1)} = C_2 (m/n)^2$.  
It is easy to check that $L_0 \geq 2$ and that $L_0 \leq L/2$, so that  
\bes 
\Del(n,m,L) =   \frac{C}{L} \lkv \lkr \frac{m}{n}\rkr^2   \rkv^{\frac{2 \nu_0}{2 \nu_0+1}}.
\ees
If  $n > C_3 m L^{\nu_0 + 1/2}$, choose $\al^2 = C_2 (m/n)^2$ and set $L_0 = L/2$. Then, 
\fr{new_rhoL0} holds and $\Del(n,m,L) = C (m/n)^2$ which completes the proof of the lower bound
in this case. 
\\
 
%%%%%%%%%%%%%%%%%%%%%%%%%%%%%%%%%%%%%%%%%%%%%%%%%%%%%%%%%%%%%%%%%%%%%%%%%%%%%%%%%%%%%%%%%%%%%%%%%5

%%%%%%%%%%%%%%%%%%%%%%%%%%%%%%%%%%%%%%%%%%%%%%%%%%%%%%%%%%%%%%%%%%%%%%%%%%%%%%%%%%%%%%%%%%%%%%%%%

\subsection{Proof of Theorem \ref{th:sparse} } 

% Let $\bC$ be a clustering matrix and $\PCJ$ be the projection matrix on the column space of matrix $(\bC \bW^T)_J$.
Denote $\bW^T = \bV$ and recall that $\bS^2 = \bC^T \bC$ where $\bS^2$ is the diagonal matrix with entries $N_k^{(l)}$, 
the number of nodes in the pair of classes $k = (k_1, k_2)$ at time $t_l$. 
Note that  $(\bC \bW^T)_J = \bC \bV_J$, so that $\PCJ = \bC \bV_J [\bV_J^T \bC^T \bC \bV_j]^{-1} \bV_J^T \bC^T$
can be written as 
\bes 
\PCJ = \bC \bS^{-1} \PSJ (\bC \bS^{-1})^T \quad \mbox{with} \quad \PSJ = \bS \bV_J (\bV_J^T \bS^2 \bV_J)^{-1} \bV_J^T \bS.
\ees
Here, $\PSJ$ is the projection matrix on the column space of matrix $\bS \bV_J$. 

Since $(\hbd, \hbC, \hm, \hJ)$ are solutions of optimization problem \fr{opt_problem_balanced},
for any   $J$  and $\bds = \vect(\bDs)$,  one has
\begin{align}  \label{new_main_ineq}
\|\ba - \hbC \hbW^T\, \hbd_{(\hJ)}\|^2 & + 2  \rho_n\,  \Pen(|\hJ|,\hm) \\
& \leq \|\ba - \bCs \bWs^T \bds_{(J)} \|^2 + 2  \rho_n\, \Pen(|J|,\ms), \nonumber
\end{align} 
It follows from \fr{true_model} that for any $m,J,\bd$ and $\bC$
\begin{align*} 
& \|\ba -   \bC \bW^T  \bd_{(J)}\|^2  = \|(\ba - \bCs \bWs^T \bds) +  (\bCs \bWs^T \bds -  \bC \bW^T \bd_{(J)})\|^2 \\
& =   \|\bxi\|^2 + \|\bC \bW^T \bd_{(J)} - \bCs \bWs^T \bds\|^2 + 2 \bxi^T (\bCs \bWs^T \bds -  \bC \bW^T \bd_{(J)}).
\end{align*}
Hence, plugging the last identity into the inequality \fr{new_main_ineq}, derive that
% for any $m$, $J$, $\bC$  and $\bd$
\begin{align} \label{new_main_ineq1}
\|\hbC \hbW^T \hbd_{(\hJ)} - \bCs \bWs^T \bds\|^2 & \leq \|\bCs \bWs^T \bds_{(J)} - \bCs \bWs^T \bds\|^2 \\
& +    \Del + \rho_n\, \Pen(|J|,\ms) - \rho_n\, \Pen(|\hJ|,\hm), \nonumber 
\end{align}
where, $\Del   = 2 |\bxi^T (\hbC \hbW^T \hbd_{(\hJ)} - \bCs \bWs^T \bds_{(J)})|$. 
Note that, due to \fr{projection} and \fr{hbd_solution}, one has  $\hbC \hbW^T \hbd_{(\hJ)} = \hPChJ\, \ba$ 
with $\ba$ is given by \fr{true_model}, $\bCs \bWs^T \bds_{(J)} = \PCJs \btes$ and $\bCs \bWs^T \bds = \btes$.  
Therefore,
\begin{align}
  & \Del    \leq  \Del_1 + \Del_2 + \Del_3, \label{new_Del123}\\
  \Del_1   = 2 |\bxi^T \PCJso \btes|, \quad %\nonumber \\
  & \Del_2    = 2 |\bxi^T \hPChJo \btes|, \quad  
  \Del_3   = 2  \bxi^T \hPChJ\, \bxi.  \nonumber
\end{align}
% Now, we need to find upper bounds for each of the terms in \fr{new_Del123}.

%%%%%%%%%%%%%%%%%%%%%%%%%%%%%%%%%%%%%%%%%%%%%%
%
In order to obtain  an upper bound for $\Del_1$ and $\Del_2$, note that by Bernstein inequality, for any $\bC$, $m$, $\bte$, $J$ and
for any $x>0$ with probability at least $1 - 2 e^{-x}$, one has 
\bes
2 |\bxi^T \PCJo \bte| \leq 2 \sqrt{2 x \rho_n\, \|\PCJo \bte \|^2 } + 4/3\,  \|\PCJo \bte \|_{\infty}\, x.
\ees
Due to $2 ab \leq a^2 + b^2$, obtain that with probability at least $1 - 2 e^{-x}$
\be \label{eq:Bernstein}
2 |\bxi^T \PCJo \bte| \leq 0.5\, \|\PCJo \bte \|^2 + 4\, \lkr \rho_n\,   +    \|\PCJo \bte \|_{\infty}/3 \rkr x
\ee 
%
%%%%%%%%%%%%%%%%%%%%%%%%%%%%%%%%%%%%%%%%%%%%%%
%
Applying \fr{eq:Bernstein}  to $\Del_1$ with  $\bte = \btes$ and  $x=t>0$,   using \fr{uni_sparse_cond} and keeping in mind that 
$\|\bLams\|_\infty = \|\btes\|_\infty \leq \rhons \leq \rho_n$, obtain
\be \label{new_Del1er}
\PP \lfi \Del_1 -  0.5\, \| \bCs \bWs^T \bds - \bCs \bWs^T \bds_{(J)} \|^2  - 4\,  \rho_n (1 + B_0/3) t \leq  0  \rfi \geq 1 - 2 e^{-t}.
\ee
%
%%%%%%%%%%%%%%%%%%%%%%%%%%%%%%%%%
%
In order to obtain an upper bound for $\Del_2$, note that, similarly to the proof of Theorem~\ref{th:oracle}, one has  
\bes
\|\hbC \hbW^T \hbd_{(\hJ)} - \bCs \bWs^T \bds\|^2   
% + \|\hPhJ\, \bxi\|^2 + \| (\bI_{NL} - \hPhJ)  \bCs \bWs^T \bds\|^2 
\geq  \| (\bI_{NL} - \hPChJ)  \bCs \bWs^T \bds\|^2 = \|\hPChJo \bCs \bWs^T \bds\|^2.
\ees
Hence,  due to condition \fr{uni_sparse_cond}, since $\|\btes\|_\infty \leq \rhons\leq \rho_n$, obtain
\begin{align*} 
& \PP \lfi \Del_2 -  0.5\,  \|\hbC \hbW^T \hbd_{(\hJ)} - \bCs \bWs^T \bds\|^2 \leq
4 \lkr   \rho_n    + 1/3\,  \|\hPChJo \btes \|_{\infty} \rkr   x   \rfi \geq \\ 
& \PP \lfi \Del_2 -   0.5\,  \| \hPChJo  \bCs \bWs^T \bds\|^2 - 4\,  \rho_n (1 + B_0/3) x \leq 0 \rfi   
\end{align*}
Set  $x  = t + R(\hm,\hJ,L)$ where $R(m,J,L)$ is defined in \fr{RDel}.
Applying inequality \fr{eq:Bernstein} with $\bte = \btes$ % \bCs \bWs^T \bds$ 
together with the union bound  over $m=1, \cdots, n$, $\bC \in \calZ_{\bal} (m,n,n_0,L,\aleph_1,\aleph_2)$ 
and $J$ with $|J| = 1, \cdots, ML,$ we derive % that with probability at least $1 - 6 e^{-t}$, one has     
\begin{align}  \label{new_Del2er}
\PP \lfi
\Del_2 \right. & -  0.5\,  \|\hbC \hbW^T \hbd_{(\hJ)} - \bCs \bWs^T \bds\|^2  \\
& \left. - 4\,  \rho_n (1    + B_0/3)
(t + R(\hm,\hJ,L)) \leq 0  \rfi \geq 1 - 6 e^{-t}. \nonumber
\end{align} 
%
%%%%%%%%%%%%%%%%%%%%%%%%%%%%%%%%%
%  
For an upper bound for $\Del_3$,  write   
$$
\Del_3   = 2  \bxi^T \hbC \hbS^{-1} \hPShJ (\hbC \hbS^{-1})^T \bxi    = 2  \boeta^T   \hPShJ    \boeta.
$$
Here, for any fixed $\hbC = \bC$, due to Corollary~\ref{cor:subgaus}, vector $\boeta = (\bC \bS^{-1})^T \bxi$
has independent sub-Gaussian components such that,  for any index $i$ corresponding to a pair of nodes 
$k = (k_1, k_2)$ at time $t_l$,  one has 
$$
\| \boeta_{i} \|_{\psi_2}^2 \leq   1.5\, e^2  \, \max( \|\btes\|_{\infty}, [N^{(l)}_{k}]^{-1}).
$$
Observe that, due to \fr{eq:Nkl} and \fr{balanced},    
$$
N^{(l)}_{{k_1},{k_2}} \geq  0.5\, n^{(l)}_{k_1} n^{(l)}_{k_2} \geq 0.5\, \aleph_1^2 (n/m)^2.
$$
Therefore, 
$$
\| \boeta_{i} \|_{\psi_2}^2 \leq  \rho_{m,n} \quad \mbox{where} \quad  \rho_{m,n}= 
\frac{3e^2}{2} \, \max \lkr \rhons,  \frac{2}{\aleph_1^2}\, \frac{m^2}{n^2}\rkr.
$$   
Therefore, by Lemma 5.5 of Vershynin (2012), obtain that there exists an absolute constant $C_0$ such that 
for any vector  $\bt$  one has 
$$
\EE [\exp(\boeta^T \bt)] \leq \exp(C_0 \rho_{m,n} \|\bt\|^2). 
$$
Applying Lemma~\ref{lem:subgaus_quadr_err} with $\sig^2 = 2 C_0 \rho_{m,n}$, derive that for any fixed 
$\bC$ and $J$, one has
\bes  
\PP \lfi 2 \boeta  \PSJ \boeta   -  4 C_0 \rho_{m,n}\, \lkr 2 |J| +   3 t \rkr \leq 0 \rfi \geq  1- e^{-t}.
\ees
Again, taking a union bound over $m=1, \cdots, n$, $\bC \in \calZ_{\bal} (m,n,n_0,L,\aleph_1,\aleph_2)$ 
and $J$ with $|J| = 1, \cdots, ML,$ we derive that, due to $\rhons \leq \rho_n$, for some absolute constant $\tCo$  one has
\begin{align}  \label{new_Del3er}
\PP \lfi \Del_3 - \tCo \max \lkr \rho_n, \frac{\hm^2}{n^2} \rkr\,  [R(\hm,\hJ,L) + t] \leq 0 \rfi 
\geq  1- e^{-t}. 
\end{align}
where $R(m,J,L)$ is defined in \fr{RDel}.
The rest of the proof  of \fr{oracle_prob_sparse1} is very similar to the proof of Theorem~\ref{th:oracle}.
In order to establish \fr{oracle_prob_sparse2}, follow the arguments of the proof of Corollary~\ref{cor:upper_smooth_DSBM}.

%%%%%%%%%%%%%%%%%%%%%%%%%%%%%%%%%%%%%%%%%%%%%%%%%%%%%%%%%%%%%%%%%%%%%%%%%%%%%%%%%%%%%%%%%%%%%%%%%
%%%%%%%%%%%%%%%%%%%%%%%%%%%%%%%%%%%%%%%%%%%%%%%%%%%%%%%%%%%%%%%%%%%%%%%%%%%%%%%%%%%%%%%%%%%%%%%%%

\subsection{ Proof  of Theorem~\ref{th:upper_graphon}, Corollary~\ref{cor:upper_graphon} and 
Theorem~\ref{th:lower_graphon} }

This section contains proofs for the upper and the lower bounds for the risks in the case  of  graphon estimation.
\\

%%%%%%%%%%%%%%%%%%%%%%%%%%%%%%%%%%%%%%%%%%%%%%%%%%%%%%%%%%%%%%%%%%%%%%%%%
  
\noindent
{\bf Proof of Theorem \ref{th:upper_graphon}. }
To prove \fr{graph_upper}, we approximate the graphon by the DSBM and use inequality \fr{oracle_expec} in  Theorem \ref{th:oracle}.
In order to find an upper bound for the bias term, we need to cluster the nodes and create an approximate 
connectivity tensor $\bQ$. For this purpose, let $h \geq 1$ be a positive integer and denote
$\kappa_j = 1 + [(\beta_j - \beta_{j-1})h]$ where $[x]$ is the largest integer
no larger than $x$.  For $k = 1, \cdots, \kappa_j$ and $j = 1, \cdots, r$, consider a set of intervals 
\bes
U_{j,k}= (U_{j,k}^{(1)}, U_{j,k}^{(2)}] \quad \mbox{with} \quad  U_{j,k}^{(1)} =  \beta_{j-1} + (k-1)/h, \quad 
U_{j,k}^{(2)} = \min\{\beta_{j-1} + k/h, \beta_j\}.
% \quad k = 1, \cdots, \kappa_j; \ j = 1, \cdots, r.
\ees
Intervals $U_{j,k}$ subdivide every interval $(\beta_{j-1}, \beta_j]$ 
into $\kappa_j$ sub-intervals of  length at most $1/h$ and the total number of intervals 
is equal to 
\bes 
m = \sum_{j=1}^r \kappa_j.
\ees
Re-number the intervals consecutively as $U_1, \cdots, U_m$ and observe that since 
$(\beta_j - \beta_{j-1})h \leq \kappa_j \leq 1 +  (\beta_j - \beta_{j-1})h$,
one has 
\be \label{sr}
h \leq m \leq h+r.
\ee
The value $m$ in \fr{sr}  acts   as a number of classes.  Indeed, if $\zeta_i \in U_k$,
we place node $i$ into class $\Om_k$ and set $\bZ_{i, j} = \II(j=k)$.

Let $\bTe^*$ be the true tensor of connection probabilities. Set $\bPhis = \bTes \bH^T$, 
$\bQ  = (\bZ^T \bZ)^{-1} \bZ^T \bTe^*$ and  $\bV  = \bQ  \bH^T$. 
Since $\bH$ is an orthogonal matrix,  the bias term  $\| \bC \bW^T \bd_{(J)} - \btes \|^2$ 
in the oracle inequality \fr{oracle_expec} is equal to 
\begin{align}  \nonumber
& \| \bC \bW^T \bd_{(J)} - \btes \|^2   =   \|\bZ \bV^{(\rho)} \bH - \bTes \|^2 = \|\bZ \bV^{(\rho)} - \bTe \bH^T\|^2\\
\label{bias_all} 
& =   \|\bZ \bV^{(\rho)} -  \bPhis \|^2 =     \|\bZ \bV^{(\rho)} -  \bPhis^{(\rho)}\|^2
+ \| \bPhis - (\bPhis)^{(\rho)}\|^2.
\end{align} 
The first term in the right-hand side of \fr{bias_all} describes how well the first $L^\rho$ columns of matrix $\bV^{(\rho)}$
represent the first $L^\rho$ columns $(\bPhis)^{(\rho)}$ of matrix $\bPhis$. The second term 
$\| \bPhis - (\bPhis)^{(\rho)}\|^2 = \| \bTes - (\bPhis)^{(\rho)} \bH \|^2$ evaluates 
how well $\bTes$ is represented by $L^\rho$ columns of its coefficients in the transform $\bH$.

The upper bound  for $\|\bZ \bV^{(\rho)} -  \bPhis^{(\rho)}\|^2$   can be found by repeating the calculations in Lemma 2.1 of Gao \etal (2015) 
with the only difference that $f$ is replaced by $v_l$ and there is an additional sum over $l = 1, \cdots L^\rho$.
Then, we obtain 
\be \label{bias1}
 \|\bZ \bV^{(\rho)} -  \bPhis^{(\rho)}\|^2  \leq K_1  2^{2 \nu_1} n^2 L^\rho\,  h^{-2 \nu_1}.
\ee
On the other hand, if $\rho <1$, then by Assumption {\bf A},
\be \label{bias2}
\| \bPhis - (\bPhis)^{(\rho)}\|^2 \leq \sum_{i_1=1}^n\, \sum_{i_2 =1}^n\ \sum_{l = L^\rho +1}^L \bv_l^2 (\zeta_{i_1}, \zeta_{i_2}) 
\leq K_2 n^2 L^{-2 \nu_2 \rho}
\ee
and $\| \bPhis - (\bPhis)^{(\rho)}\|^2 =0$ if $\rho=1$. 
Now, note that  for  given $m$ and $\rho$,  one has $|J| \leq m^2 L^\rho$, so that 
\be \label{gra_pen} 
\Pen(|J|,m) \leq C \lkv n \log m +   m^2 L^\rho \log(25\, L^{1-\rho})\rkv.
\ee  
Therefore,   \fr{sr}--\fr{gra_pen}  yield  \fr{graph_upper}.
\\

%%%%%%%%%%%%%%%%%%%%%%%%%%%%%%%%%%%%%%%%%%%%%%%%%%%%%%%%%%%%%%%%%%%%%%%%%%%%%%%%%%%%%%%%%%%%%%%%%%%%%%%%%%%%%%%%%%%%

%\medskip

\noindent
{\bf Proof of Corollary \ref{cor:upper_graphon}.  }
Note that if    $\nu_1 =\infty$, one can set $h=2$, so that the first term in \fr{graph_upper} is equal to zero, 
$h+r \leq 3 r$ and 
\bes  
\frac{\EE \|\hbLam - \bLams\|^2}{n^2\,L} \leq C \lfi \frac{I(\rho <1)}{L^{2 \rho \nu_2+1}} + 
\lkr\frac{r}{n}\rkr^2 \frac{1 + (1 - \rho)\log L}{L^{1-\rho}}  + \frac{\log r}{n L}\rfi.
\ees
Minimizing this expression with respect to $\rho \in [0,1]$ obtain the result in  \fr{graph_upper_cases} for $r = r_{n,L}$.
If $r=r_0$, then 
\be \label{risk2}
\frac{\EE \|\hbLam - \bLams\|^2}{n^2\,L} \leq C \lfi  
\frac{L^{\rho-1}}{h^{2 \nu_1}} + \frac{I(\rho <1)}{L^{2 \rho \nu_2+1}} + \frac{h^2\log L}{n^2 L^{1-\rho}}  
+ \frac{\log h}{n  L}  \rfi.
\ee 
Minimizing \fr{risk2} with respect to $h$ and $\rho$, obtain that the   values $h^*$ and $L^* = L^{\rho^*}$  delivering the minimum in 
\fr{risk2} are such that    $h^* \asymp  (n^2/\log L)^{1/(2(\nu_1+1)}$.
Hence, for some absolute constants $C$, one has  $\log (h^*) \leq C \log n$  and 
\bes
L^* =   L^{\rho^*} = \min \lfi L, C\, \lkr  n^2/\log L \rkr^{\frac{\nu_1}{(\nu_1 +1)(2 \nu_2+1)}}\rfi. 
\ees
Therefore,  \fr{graph_upper_cases} holds for $r = r_0$.
\\

%%%%%%%%%%%%%%%%%%%%%%%%%%%%%%%%%%%%%%%%%%%%%%%%%%%%%%%%%%%%%%%%%%%%%%%%%%%%%%%%%%%%%%%%%%%%%%%%%%%%%%%%%%%%%%%%%%%%
%%%%%%%%%%%%%%%%%%%%%%%%%%%%%%%%%%%%%%%%%%%%%%%%%%%%%%%%%%%%%%%%%%%%%%%%%%%%%%%%%%%%%%%%%%%%%%%%%%%%%%%%%%%%%%%%%%%%

% \subsection{Proof of Theorem \ref{th:lower_graphon}}
%\label{sec:lower}

 \noindent
 {\bf Proof of Theorem \ref{th:lower_graphon}. } 
We consider the cases when $r = r_{n,L}$ and $r=r_0$,
corresponding to piecewise constant and piecewise smooth graphon, separately.
\\

\underline{\bf Piecewise constant graphon. }
Assume, without loss of generality, that $j = n/r$ is an integer.  
Consider a set up where the nodes are grouped into $r$ classes and 
values of $\zeta_j$'s are fixed:   
\bes
\zeta_{k j + i} = \beta_k + (i-1/2)  (\beta_{k+1} - \beta_k)/j,
\quad k = 0, \cdots, r-1, \ i = 1, \cdots, j.
\ees
Then, there are $j$ nodes in each class.   
Let $\bG_{k,k,l} =0$ for any $k=1, \cdots, r$ and $l=1, \cdots, L$.

Consider an even number $L_0$ such that $1 \leq L_0 \leq L/2$ and a set of vectors $\bom \in \{ 0,1\}^{L_0}$ 
with exactly $L_1 = L_0/2$ nonzero entries.
By Lemma 4.10 of Massart (2007), there exists a subset $\calT$ of those vectors such that 
for any $\bom, \bom' \in \calT$ one has 
\be \label{bombom}
\|\bom - \bom'\|_H \geq L_0/4, \quad
\log |\calT| \geq 0.233 (L_0/2) \log 2 \geq 0.08 L_0.
\ee 
Denote $\calK = \{(k_1, k_2):\ 1 \leq k_1 < k_2 \leq r\}$ and 
let  $\calT_{k_1, k_2}$ be the copies of $\calT$ for $(k_1, k_2) \in \calK$. 
Denote $\tilde{\calT} = \bigotimes_{(k_1,k_2) \in \calK}\, \calT_{k_1, k_2}$ and
observe that 
\bes
\log|\tilde{\calT}| = r(r-1)/2 \log |\calT| \geq 0.02 r^2\, L_0.
\ees
Then, $\tilde{\bom} \in \tilde{\calT}$ are   binary tensors with elements $\bom_{l}^{(k_1, k_2)}$,$l=1, \cdots, L_0$
and  $(k_1,k_2) \in \calK$.
Consider a set of functions $f^{(\tilde{\bom})}$ indexed by  $\tilde{\bom}$
such that, for $\beta_{k_1-1}  < x \leq \beta_{k_1}$ and $\beta_{k_2-1} < y \leq \beta_{k_2}$,
their coefficients in the transform $\bH$ are given by 
\bes
\bv_1^{(\tilde{\bom})} (x,y) = \sqrt{L}/2; \quad \bv_l^{(\tilde{\bom})}  (x,y) =  
\al  \bom_{l-L_0}^{(k_1, k_2)},\ l=L_0+1, \cdots, 2L_0  , \  (k_1,k_2) \in \calK.
% \quad \mbox{if}
% \beta_{k_1-1}  < x \leq \beta_{k_1}, \ \beta_{k_2-1} < y \leq \beta_{k_2}.
\ees  
% where $\tilde{\bom}$ is a binary matrix with elements $\bom_{l}^{(k_1, k_2)}$,$l=1, \cdots, L_0$
% and  $(k_1,k_2) \in \calK$.
%
Then, Assumption \fr{A1} holds. In order condition \fr{A2} is satisfied, we set
\be \label{rhoL0}
\al^2 \leq C_1 L_0^{-(2 \nu_2 +1)}  \quad \mbox{with} \quad C_1 \leq  \min(K_2 2^{1 - 2\nu_2},1/8).
\ee 
Denote by $\bLam$ and $\bLam'$ the probability tensors corresponding, respectively,  to $\tilde{\bom}$ and 
$\tilde{\bom'}$ in $\tilde{\calT}$. Then, due to \fr{bombom} and the symmetry, 
\beqns
\|\bLam - \bLam'\|^2 &\geq &  2 \al^2 \sum_{k_1 =1}^r \sum_{k_2 = k_1 +1}^r 
\|\bom^{(k_1, k_2)} - \bom'^{(k_1, k_2)}\|_H \lkr \frac{n}{r}\rkr^2  \geq \frac{\al^2 n^2 L_0}{8};\\
\|\bLam - \bLam'\|^2 &\leq &    \al^2 r(r-1) (n/r)^2 L_0 \leq \al^2 n^2 L_0.
\eeqns
Note that  one has $1/4 \leq f(x,y,t) \leq 3/4$ provided
$\|\bH^T \bom ^{(k_1, k_2)}\|_\infty \leq 1/4$ for for $(k_1, k_2) \in \calK$.
By  Assumption \fr{H_assump}, the latter is guaranteed by 
$\al^2 L_0^2/L \leq 1/4$, so that, due to $L  \geq 2 L_0$ and $\nu_2 \geq 1/2$,
it is ensured by \fr{rhoL0}. Then, by Lemma~\ref{lem:Gao2015}, one has 
\bes
K (\PP_{\bLam}, \PP_{\bLam'}) \leq 8\, \| \bLam - \bLam'\|^2  \leq
8  \al^2 n^2 L_0 \leq \log |\tilde{\calT}|/8   
\ees  
provided 
\be \label{rho_con1}
\al^2 \leq C_2 (r/n)^2
\ee
 where $C_2$ is an absolute constant.
Therefore, application of Theorem 2.5 of Tsybakov (2009) yields \fr{graph_lower_ineq}
with 
$\Del(n,L) = C\,\al^2 L_0/L$ 
where $C$ is an absolute constant.

Now, we denote $C_3^2 = (C_2/C_1) 2^{-(2 \nu_2 +1)}$ and consider two cases.   
If $n \leq C_3 r L^{\nu_2 + 1/2}$, choose $L_0 = [(C_1/C_2) (n/r)^2]^{1/(2 \nu_2 +1)}$ 
which leads to $\al^2 = C_1 L_0^{-(2 \nu_2 +1)} = C_2 (r/n)^2$.  
It is easy to check that $L_0 \geq 2$ and that $L_0 \leq L/2$, so that  
\bes 
\Del(n,L) =   \frac{C}{L} \lkv \lkr \frac{r}{n}\rkr^2   \rkv^{\frac{2 \nu_2}{2 \nu_2+1}}.
\ees
If  $n > C_3 r L^{\nu_2 + 1/2}$, choose $\al^2 = C_2 (r/n)^2$ and set $L_0 = L/2$. Then, 
\fr{rhoL0} holds and $\Del(n,L) = C (r/n)^2$ which completes the proof of the lower bound 
when $r = r_{n,L}$.
\\
 
%%%%%%%%%%%%%%%%%%%%%%%%%%%%%%%%%%%%%%%%%%%%%%%%%%%%%%%%%%%%%%%%%%%%%%%%%%%%%%%%%%%%%%%%%%%%%%%%%5

\underline{\bf Piecewise smooth graphon. }
Since $r=r_0$ is a fixed quantity, without loss of generality, we set $r=1$.
Let  $\zeta_j = j/n, j=1, \cdots, n$, be fixed. Let $h$   be a positive integer, $1 \leq h \leq n$,
and denote $\del = 1/h$. Consider a kernel function $F(x)$ such that 
$F(x)$ is $n_1 > \nu_1$ times continuously differentiable  and for any $x,x' \in \RR$ and some $C_F >0$
\be \label{kern_properties}
\supp (F) = (-1/2; 1/2), \quad |F(x) - F(x')| \leq C_F |x  - x'|^{\nu_1}.
\ee
It follows from \fr{kern_properties} that $|F(x)| \leq C_F$ for any $x$.
Denote % $u_k = (k-1)\del + \del/2$, $k=1, \cdots, h$ and 
\be \label{Psi}
\Psi_{k_1, k_2} (x,y) = h^{-\nu_1} F(h(x - u_{k_1})) F(h(y - u_{k_2}))
%\quad \mbox{where} \quad u_k = (k-1/2)\del, \ k=1, \cdots, h.
\ee
where $u_k = (k-1/2)\del, \ k=1, \cdots, h$.
It is easy to see that $\Psi_{k_1, k_2} (x,y)= \Psi_{k_2, k_1} (y,x)$ 
and, for different pairs  $(k_1, k_2)$,  functions  $\Psi_{k_1, k_2} (x,y)$
have disjoint supports. Similar to the case of the piecewise constant graphon, 
consider an even  $L_0 \leq L/2$ and a set of vectors $\bom \in \{ 0,1\}^{L_0}$ with exactly $L_0/2$ nonzero entries.
By Lemma 4.10 of Massart (2007), there exists a subset $\calT$ of those vectors such that \fr{bombom} holds
for any $\bom, \bom' \in \calT$.
Let again $\calT_{k_1, k_2}$ be the copies of $\calT$ for $(k_1, k_2) \in \calK$ 
and denote $\tilde{\calT} = \prod_{(k_1,k_2) \in \calK} \calT_{k_1, k_2}$
where $\calK = \{(k_1, k_2):\ 1 \leq k_1 < k_2 \leq h\}$.
Then, 
\be \label{card11}
\log|\tilde{\calT}|= h(h+1)/2\, \log|{\calT}| \geq 0.04\, L_0 h^2. 
\ee
For any $\tilde{\bom} \in \tilde{\calT}$ and $l=L_0+1, \cdots, 2L_0$, define 
\begin{align} 
& \bv_1^{(\tilde{\bom})} (x,y) = \sqrt{L}/2; \label{bv_val} \\
& \bv_l^{(\tilde{\bom})} (x,y) =  
\al   \sum_{k_1 =1}^h \sum_{k_2 = k_1 + 1}^h \bom_{l-L_0}^{(k_1, k_2)}\, [\Psi_{k_1, k_2} (x,y) + \Psi_{k_2, k_1} (x,y)].
\nonumber 
%\quad l=L_0+1, \cdots, 2L_0.
 \end{align}   
It is easy to see that $\bv^{(\tilde{\bom})} (x,y) = \bv^{(\tilde{\bom})} (y,x)$ for any $x,y \in [0,1]$. 
Now we need to check that conditions \fr{A1} and \fr{A2} hold.

Note that for any $x,y,x',y'$, due to \fr{kern_properties}, obtain
\begin{align}
& |\Psi_{k_1, k_2} (x,y)   - \Psi_{k_1, k_2} (x',y')|   \leq 
h^{-\nu_1} \lkv |F(h(x-u_{k_1})- F(h(x'-u_{k_1})| |F(h(y-u_{k_2})| \right. \nonumber \\
& \left. +  |F(h(y-u_{k_2})- F(h(y'-u_{k_2})| |F(h(x'-u_{k_1})| \rkv  
  \leq C_{\psi}  \lkv |x-x'| + |y-y'| \rkv^{\nu_1}, \label{ineqPsi}
\end{align}
where constant $C_{\psi}$ depends  only on $C_F$ and $\nu_1$.
Since functions $\Psi_{k_1, k_2} (x,y)$ have disjoint supports for different pairs of indices  $(k_1, k_2)$, 
the sums below have at most four nonzero terms. Then, \fr{ineqPsi} implies that 
\begin{align*}
& \left| \bv_l^{(\tilde{\bom})} (x,y) - \bv_l^{(\tilde{\bom})} (x',y')\right| 
\leq \al \sum_{k_1, k_2=1}^h   |\Psi_{k_1, k_2} (x,y)   - \Psi_{k_1, k_2} (x',y')|\\
& \leq 4 \al\,  C_{\psi}\  \lkv |x-x'| + |y-y'| \rkv^{\nu_1},
\end{align*}
so that \fr{A1} holds if $\al \leq K_1/(4   C_{\psi})$.  
Also, it is easy to check that, by \fr{kern_properties}, one has
$[\bv_l^{(\tilde{\bom})} (x,y)]^2 \leq C_v^2  \al^2 h^{-2\nu_1}$ where $C_v$ depends  only on $C_F$ and $\nu_1$.
Therefore, 
\bes
\sum_{l = L_0+1}^{2L_0} (l-1)^{2 \nu_2} [\bv_l^{(\tilde{\bom})} (x,y)]^2  \leq C_v^2 \al^2 h^{-2\nu_1} L_0^{2 \nu_2+1}.
\ees  
 Hence, both assumptions, \fr{A1} and \fr{A2} are valid provided
\be \label{rhoineq}
\al \leq \min \lkr K_1/(4   C_{\psi}),\ \sqrt{K_2}/C_v\, h^{\nu_1} L_0^{-(\nu_2+1/2)}\rkr.
\ee
Denote by $\bLam$ and $\bLam'$ the probability tensors corresponding to $\tilde{\bom}$ and 
$\tilde{\bom'}$ in $\tilde{\calT}$, respectively. 
Let $\tilde{\bv} (x,y)$ be a vector with  $\tilde{\bv}_1 (x,y) = 0$ and $\tilde{\bv}_l (x,y) = \bv_l (x,y)$
for $l\geq 2$. 
By \fr{H_assump} and \fr{bv_val}, since $\Psi_{k_1, k_2} (x,y)$ have disjoint supports,
we derive that for any $\tilde{\bom} \in \tilde{\calT}$ one has
\bes
\left| \Lam_{i,j,l}  - 1/2\right| \leq    \| \bH \tilde{\bv}  (x,y) \|_\infty
\leq  \|\tilde{\bv}  (x,y) \|_1/\sqrt{L} \leq C_v \al L_0 h^{-\nu_1}/\sqrt{L}.
\ees
Hence, $\bLam_{i,j,l}, \bLam'_{i,j,l} \in [1/4;3/4]$ provided
\be \label{rho_as1}
\al \leq \sqrt{L}\, h^{\nu_1}/(4 C_v  L_0).
\ee  
Then, by   Lemma~\ref{lem:Gao2015}, since each $\bom \in   \calT$ has exactly $L_0/2$ nonzero entries  
\begin{align*}
K (\PP_{\bLam}, \PP_{\bLam'}) & \leq    8 \|\bLam - \bLam'\|^2  
\leq 8 \al^2  L_0/2\ \sum_{i,j}^n \sum_{k_1, k_2}^h 4\, \lkv \Psi_{k_1, k_2} \lkr \frac{i}{n}, \frac{j}{n} \rkr \rkv^2\\
& =   16 \al^2  h^{-2\nu_1} L_0 \ \sum_{i,j}^n \sum_{k_1, k_2}^h F^2 (h(i/n-u_{k_1})) F^2 (h(j/n-u_{k_2})) \\
& =   16 \al^2  h^{-2\nu_1} L_0 \lkv  \sum_{i=1}^n \sum_{k=1}^h F^2(h(i/n-u_k))\rkv^2.
\end{align*}  
Here, 
\begin{align} \label{F_integral}
&  \sum_{i=1}^n \sum_{k=1}^h F^2(h(i/n-u_k))   = 
\sum_{k=1}^h \quad  \sum_{u_k-\del/2 <i/n \leq u_k+\del/2} F^2(h(i/n-u_k))   \\
& \approx   n \sum_{k=1}^h \int_{u_k-\del/2}^{u_k+\del/2} F^2(h(i/n-u_k)) dx 
   =   n \int_{-1/2}^{1/2}  F^2(z) dz,  \nonumber
\end{align} 
so that 
\be \label{KL11}
K (\PP_{\bLam}, \PP_{\bLam'}) \leq 16\, \|F\|_2^4\, \al^2  h^{-2\nu_1} n^2 L_0
\ee 
where $\|F\|_2$ is the $L^2$-norm of $F$.
On the other hand,  due to the first inequality in \fr{bombom} and \fr{F_integral}, obtain  
\begin{align*}
&\frac{\|\bLam - \bLam'\|^2}{n^2\, L} \geq    \frac{2 \al^2}{n^2\, L} \, \sum_{l=L_0+1}^{2L_0} \, \sum_{k_1 =1}^r \sum_{k_2 = k_1 +1}^r  
\sum_{i=1}^n \sum_{j=1}^n [\bom_l^{(k_1, k_2)} - {\bom'_l}^{(k_1, k_2)}]^2 [\Psi_{k_1, k_2} (i/n, j/n)]^2\\
&\geq  \frac{\al^2 h^{-2\nu_1} L_0}{4\, n^2 \, L} \lkv \sum_{i=1}^n \sum_{k=1}^h F^2(h(i/n-u_k)) \rkv^2 
\geq  \frac{\al^2 h^{-2\nu_1} L_0}{8\, L} \|F\|_2^4.
\end{align*} 
 Application of Theorem 2.5 of Tsybakov (2009) yields that
\fr{graph_lower_ineq} holds with 
$\Del(n,L) = C\, \al^2 h^{-2\nu_1} L_0/L$ provided $K (\PP_{\bLam}, \PP_{\bLam'}) \leq 1/8 \log|\tilde{\calT}|$,
which, due to \fr{card11} and \fr{KL11}, is guaranteed by 
\bes
4 \|F\|_2^4 \al^2  h^{-2\nu_1} n^2 L_0 \leq (1/8)\,  (0.04 L_0 h^2)/2
\ees
and leads to  the following restriction on $\al$:
\be \label{rho_as2}
\al \leq  h^{\nu_1 +1} n^{-1} /(40 \sqrt{2} \|F\|_2^2).
\ee
Set 
\be \label{hL0}
h = n^{\frac{1}{\nu_1+1}},\   L_0 = \min\lkr n^{\frac{2 \nu_1}{(\nu_1+1)(2 \nu_2 +1)}},\frac{L}{2} \rkr,\   
n_L = \lkr \frac{L}{2} \rkr^\frac{(\nu_1+1)(2 \nu_2 +1)}{2 \nu_1} 
\ee 
and consider two cases.

If $n \leq n_L$, then $L_0$ is given by the first expression 
in \fr{hL0} and    inequalities \fr{rhoineq}, \fr{rho_as1} and \fr{rho_as2}  hold if $\al^2$ is a small enough 
absolute constant that depends on $\nu_1, \nu_2, K_1$ and $K_2$ only. In this case, 
$\Del(n,L) = C\, L^{-1}    n ^{-\frac{4 \nu_1 \nu_2}{(\nu_1 +1) (2\nu_2+1)}}$
and \fr{graph_lower_cases} is valid. 

If $n > n_L$, then $L_0=L/2$ and again all inequalities \fr{rhoineq}, \fr{rho_as1} and \fr{rho_as2}  hold if $\al^2$ is a small enough 
absolute constant that depends on $\nu_1, \nu_2, K_1$ and $K_2$ only.
In this case, $\Del(n,L) = C\,  n^{-\frac{2\nu_1}{\nu_1 +1}}$ which completes the proof.

%%%%%%%%%%%%%%%%%%%%%%%%%%%%%%%%%%%%%%%%%%%%%%%%%%%%%%%%%%%%%%%%%%%%%%%%%%%%%%%%%%%%%%%%%%%%%
%%%%%%%%%%%%%%%%%%%%%%%%%%%%%%%%%%%%%%%%%%%%%%%%%%%%%%%%%%%%%%%%%%%%%%%%%%%%%%%%%%%%%%%%%%%%% 

\subsection{Proofs of supplementary statements}
\label{sec:supplement}

In this section, we formulate and prove supplementary lemmas used in the proofs of other statements.
\\

%%%%%%%%%%%%%%%%%%%%%%%%%%%%%%%%%%%%%%%%%%%%%%%%%%%%%%%%%%%%%%%%%%%%%%%%%%%%%%%%%%%%%%%%%%%%%%%%%%%%%%%5

% lemma 2
\begin{lemma} \label{lem:bern_err}
Let $a_i$ be independent Bernoulli$(\te_i)$ variables and 
consider vectors $\ba$ and $\bte$ with components $a_i$ and $\te_i$,
respectively. Then, for any vector $\bz$ and any positive $t$ and   $\al$ one has
\be \label{expineq_bernoulli}
\PP \lkr 2 |\bz^T(\ba - \bte)| >  \af \| \bz \|^2  +  2t/\al \rkr \leq 2 e^{-t}, \quad
\EE \lkv \exp (\bz^T(\ba - \bte) \rkv \leq \exp(\|\bz\|^2/8).
\ee
\end{lemma}

\noindent
{\bf Proof. } Validity of the  Lemma follows from Hoeffding inequality (see, e.g., Massart (2007)).
\\

%%%%%%%%%%%%%%%%%%%%%%%%%%%%%%%%%%%%%%%%%%%%%%%%%%%%%%%%%%%%%%%%%%%%%%%%%%%%%%%%%%%%%%%%%%%%%%%%%%%%%%%%%

% lemma 3
\begin{lemma} \label{lem:subgaus_quadr_err}
Let $\bxi$ be a vector with independent components, $\EE \bxi = 0$ and such that for any vector 
$\bz$ and some $\sig>0$ one has
\be \label{subgaus_cond}
\EE \lkv \exp (\bz^T \bxi) \rkv \leq \exp(\|\bz\|^2 \sig^2/2).
\ee
Then, for any matrix  $\bA$ and any positive $t$   one has
\be \label{subgaus_expineq_quadr}
\PP \lfi \|\bA\bxi \|^2  > \sig^2 \lkr 2 \| \bA \|^2 +   3 \,\| \bA \|_{op}^2\, t \rkr \rfi \leq  e^{-t}.
\ee
\end{lemma}

\noindent
{\bf Proof. } Denote   $\bSig = \bA^T \bA$ and use Theorem 2.1 of Hsu \etal (2012).
Obtain, for any $t>0$,  
\bes
\PP \lkv \|\bA \bxi \|^2  > \sig^2 \lkr \Tr(\bSig) + 2   \sqrt{t\, \Tr(\bSig^2)} 
+2 t\, \|\bSig\|_{op}  \rkr \rkv \leq e^{-t}.  
\ees
Note that $\Tr(\bSig^2) \leq \|\bSig\|_{op} \Tr(\bSig)$ and $2 \sqrt{xy} \leq x+y$ for any positive $x$ and $y$.
In order to complete the proof, recall that $\|\bSig\|_{op} = \| \bA \|_{op}^2$ and $\Tr(\bSig) = \| \bA \|^2$.
\\

%%%%%%%%%%%%%%%%%%%%%%%%%%%%%%%%%%%%%%%%%%%%%%%%%%%%%%%%%%%%%%%%%%%%%%%%%%%%%%%%%%%%%%%%%%%%%%%%%%%%%%%%%

\begin{corollary} \label{corr:quadr_err}
Let $a_i$ be independent Bernoulli$(\te_i)$ variables and 
consider vectors $\ba$ and $\bte$ with components $a_i$ and $\te_i$,
respectively. Then, for any matrix  $\bA$ and any positive $t$   one has
\be \label{expineq_quadr}
\PP \lkr \|\bA (\ba - \bte)\|^2  > \frac{  \| \bA \|^2}{2} +  \frac{3t\,\| \bA \|_{op}^2}{4}\rkr \leq  e^{-t}.
\ee
\end{corollary}

\noindent
{\bf Proof. } Denote $\bxi = \ba-\bte$ and note that, due to \fr{expineq_bernoulli}, condition \fr{subgaus_cond}
holds with $\sig = 1/2$.
\\

%%%%%%%%%%%%%%%%%%%%%%%%%%%%%%%%%%%%%%%%%%%%%%%%%%%%%%%%%%%%%%%%%%%%%%%%%%%%%%%%%%%%%%%%%%%%%%%%%%%%%%%%%

% lemma 4
\begin{lemma} \label{lem:packing} {\bf (The packing lemma). }
Let $\calZ (m,n) \subseteq \calM (m,n)$ be a collection of clustering matrices.
Then, there exists a subset $S_{n,m} (r) \subset \calZ (m,n)$ such  that 
for $\bC_1, \bC_2 \in \calZ (m,n)$ one has 
$\| \bC_1 - \bC_2 \|_H =  \| \bC_1 - \bC_2 \|^2 \geq r$
and 
$\log |S_{n,m} (r)| \geq  \log|\calZ (m,n)| - r \log (n e m/r)$.
\end{lemma}

\noindent
{\bf Proof. } For any clustering matrix $\bC$ define the $r$-neighborhood of $\bC$ as 
\bes
\calB (\bC, r) = \lfi \tilbC \in \calZ (m,n):\ \|\tilbC - \bC\|_H \leq r \rfi.
\ees
Let $\tilS_{n,m}(r)$ be the covering set of $\calZ (m,n)$ and $|\tilS_{n,m}(r)|$
be the covering number, the smallest number of closed balls of radius $r$ whose union
covers $\calZ (m,n)$.  Let $|S_{n,m}(r)|$ be the packing number of $\calZ (m,n)$,
the largest number of points with the distance at least $r$ between them. 
It is known (see Pollard (1990), page 10) that
\be \label{pollard_ineq}
|\tilS_{n,m}(r)| \leq |S_{n,m}(r)| \leq |\tilS_{n,m}(r/2)|
\ee 
Note that $|\calZ (m,n)| \leq \sum_{\bC \in \tilS_{n,m}(r)} |\calB (\bC, r)| \leq 
|\tilS_{n,m}(r)| \max_{\bC \in \tilS_{n,m}(r)} |\calB (\bC, r)|$, so that
\be \label{cardF}
|\tilS_{n,m}(r)| \geq  |\calZ (m,n)| \Bigg/  \max_{\bC \in \tilS_{n,m}(r)} |\calB (\bC, r)| 
\ee 
and, also, 
\be \label{cardB}
|\calB (\bC, r)| \leq {n \choose r} m^r \leq \lkr \frac{ne}{r}\rkr^r m^r = \lkr \frac{n e m}{r}\rkr^r.
\ee
Now, combining \fr{pollard_ineq} -- \fr{cardB}, obtain
$\log |S_{n,m}(r)| \geq \log |\tilS_{n,m}(r)|  \geq \log |\calZ (m,n)| - r \log(n e m/r)$
which completes the proof.
\\

%%%%%%%%%%%%%%%%%%%%%%%%%%%%%%%%%%%%%%%%%%%%%%%%%%%%%%%%%%%%%%%%%%%%%%%%%%%%%%%%%%%%%%%%%%%%%%%%%%%%%%%%%

% lemma 5
\begin{lemma} \label{lem:cardinality}
Let $\ga m$ and $n/m$ be positive integers. Then, for $n_0 \geq 1$ and $\ga m \geq 2$, one has
\begin{align}
& \log \lfi (\ga n)! \rfi - \ga m \, \log \lfi \lkr n/m \rkr !\rfi \geq \frac{\ga n}{4} \log(\ga m); \label{AA1}\\
& \log \lfi {m \ga \choose n_0} \lkr \frac{n}{m} \rkr^{n_0} (n_0-1)!  \rfi \geq n_0 \log\lkr \frac{\ga n}{n_0} \rkr 
+ \frac{n_0}{2} \log \lkr \frac{n_0}{2}\rkr; \label{AA2}\\
& 6 \ga n \log(\ga m)  -  \ga n  \log(32\, m   \ga \, e) \geq 0. \label{AA3}
\end{align}
\end{lemma}

\noindent
{\bf Proof. } Note that due to \fr{inequalities}, one has 
\bes
A_1  = \log \lfi (\ga n)! \rfi - \ga m \, \log \lfi \lkr n/m \rkr !\rfi  
\approx \ga n \log (\ga m) - \frac{\ga m}{2} \log \lkr\frac{2 \pi n}{m}\rkr + \frac{1}{2} \log(2\pi \ga n).
\ees 
Consider a function
\bes
F(x) = \frac{3 \ga n}{4} \log x  - \frac{x}{2} \log \lkr \frac{2 \pi \ga n}{x} \rkr + \frac{1}{2} \log(2 \pi \ga n).
\ees 
It is easy to check that $F(1)=0$ and that 
$F'(x) =  0.75\, \ga n/x  - 0.5 \log(2 \pi \ga n/(x e))$. Replacing $2 \pi \ga n/(x e)$ in $F'(x)$ by $z$ and noting that the inequality
$3e/(4 \pi) z > \log z$ is true for any $z>0$, we confirm that $F'(x) >0$. so that $F(x) >0$ for any $x \geq 1$. 
Finally, in order to prove \fr{AA1}, observe that $A_1 = F(\ga m) +  \ga n \log(\ga m)/4$.

For the sake of proving \fr{AA2},  note that for every $n_0 \geq 1$ one has 
$\log[(n_0-1)!] \geq 0.5\, n_0 \log(n_0/2)$ and apply the first inequality in \fr{inequalities}.

The validity of  inequality \fr{AA3} follows from $\ga m \geq 2$ and the fact that $\log(64 e) < 6$.  
\\

%%%%%%%%%%%%%%%%%%%%%%%%%%%%%%%%%%%%%%%%%%%%%%%%%%%%%%%%%%%%%%%%%%%%%%%%%%%%%%%%%%%%%%%%%%%%%%%%%%%%%%%%%

% lemma 6
\begin{lemma} \label{lem:subgaussian}
Let $a_i$ be independent  Bernoulli$(\te_i)$ variables with $0\leq \te_i \leq \rho < 1/2$, $i=1, \cdots, n$. 
Denote $\xi_i = a_i - \te_i$  and 
\be \label{sum_subgaus}
\zeta = \frac{1}{\sqrt{n}} \ \sum_{i=1}^n \xi_i.
\ee
Then, $\zeta$ is a sub-Gaussian random variable with sub-Gaussian norm 
\be \label{subgaus}
\|\zeta\|_{\psi_2} = \sup_{p \geq 1} \lkr p^{-1/2}\, \EE |\zeta|^p\rkr^{1/p} \leq \frac{3 e^2}{2} \, \max \lkr \rho, \frac{1}{n} \rkr. %, \quad p \geq 1.
\ee
\end{lemma}

\noindent
{\bf Proof. } First, observe that we need to check inequality of the type \fr{subgaus} only for an even integer $p=2m$.
Recall that $\|\zeta\|_{\psi_2} \leq K_0$ provided 
$\EE |\zeta|^p \leq (K_0 \sqrt{p})^p = (K_0^2 p)^{p/2}$ for every $p \geq 1$. Suppose that  
 $\EE |\zeta|^{2m}  \leq (K_1 \sqrt{2m})^{2m} = (2 K_1^2 m)^m$ for every integer $m \geq 1$ and some  some $K_1 >0$.    
For any $p \geq 1$ choose $m$ such that $2m-2 < p \leq 2m$ which implies  $2m < p+2$. Then,
$$
\EE |\zeta|^p \leq \lkr \EE  \zeta^{2m} \rkr^{\frac{p}{2m}} \leq \lkr (2 K_1^2 m)^m \rkr^{\frac{p}{2m}} \leq \lkr   K_1^2 (p+2) \rkr^{\frac{p}{2}}
$$
Observing that $p+2 \leq 3p$ for $p\geq 1$, obtain that  $\EE |\zeta|^p \leq  [(\sqrt{3} K_1)^2 p]^{p/2}$. 
Hence, \fr{subgaus} is valid provided that for any $m \geq 1$
\be \label{2m_ineq}
\EE (\zeta)^{2m} = \EE \lkr \frac{1}{\sqrt{n}} \ \sum_{i=1}^n \xi_i \rkr^{2m} 
\leq  (e^2 m \rho_n)^m, \quad \rho_n= \max(\rho, 1/n). % ,\quad m \geq 1.
\ee
Note that 
\bes  
\EE (\zeta)^{2m} =  \sum \frac{\EE(\xi_{i_1}^{m_1} \ldots \xi_{i_l}^{m_l})}{n^m}
\ees
where the sum is taken over all positive integers $l, m_1, \ldots, m_l, i_1, \ldots, i_l$ such that $2m = m_1 + m_2 + m_l$
and $i_1, \ldots, i_l$ are distinct. Note that $\EE \xi_i =0$ for every $i$, hence, in the sum the terms with $m_l=1$ are 
equal to zero. Therefore, 
\be \label{2m_ineq1}
\EE (\zeta)^{2m} \leq  \sum_{l}  S(m,l) P(n,l) \frac{\EE (\xi_{i_1}^{m_1} \ldots \xi_{i_l}^{m_l})}{n^m}
\ee
where $S(m,l)$ is the number of partitions of $2m$ unlabeled objects into $l$ distinct subparts of size at least two,
and $P(m,l)$ is the number of ways of choosing $l$ variables out of $n$ (order matters since powers can be different):
$$
S(m,l) = {2m-l-1 \choose l-1}, \quad P(n,l) = \frac{n!}{(n-l)!} \leq n^l.
$$
Here, $2m - l-1 \geq l-1$, so that $l \leq m$. Observe that due to
$$
\EE \xi_i^k =  \te_i(1-\te_i) [(-\te_i)^{k-1} + (1-\te_i)^{k-1}] \leq \rho,
$$ 
\fr{2m_ineq1} yields 
\bes 
\EE (\zeta)^{2m} \leq  \sum_{l=1}^{\min(m,n)}  {2m-l-1 \choose l-1} \frac{\rho^l}{n^{m-l}}.
\ees
By considering the cases $\rho \geq 1/n$ and $\rho < 1/n$ separately, it is easy to show that
$$
\frac{\rho^l}{n^{m-l}} \leq [\max(\rho, n^{-1})]^m  = \rho_{n}^m,
$$
so that 
\bes 
\EE (\zeta)^{2m} \leq \rho_{n}^m\ \sum_{l=0}^{\min(m,n)}  {2m-l  \choose l}  \leq 
\rho_{n}^m\ \sum_{l=0}^{\infty}  \frac{(2m)^l}{l!} \leq \rho_{n}^m\, e^{2m} = (e^2 \rho_{n})^m
\ees
which implies \fr{2m_ineq} and, hence, \fr{subgaus}.
\\

%%%%%%%%%%%%%%%%%%%%%%%%%%%%%%%%%%%%%%%%%%%%%%%%%%%%%%%%%%%%%%%%%%%%%%%%%%%%%%%%%%%%%%%%%%%%%%%%%%%%%%%5

\begin{corollary} \label{cor:subgaus}
Let $\ba \in [0,1]^n$ be a random vector with 
 independent Bernoulli$(\te_i)$ components. Let $\bxi = \ba - \bte$ and $\bZ \in \calM (M,N)$ be a clustering matrix
with $\bS^2 = \bZ^T \bZ = \diag(N_1, \cdots, N_M)$. Define $\boeta = \bS^{-1} \bZ^T \bxi$. 
Then, vector $\boeta$ has independent sub-Gaussian 
components with  
\bes
\| \boeta_k \|_{\psi_2}^2 \leq   \frac{3e^2}{2} \, \max( \|\te\|_{\infty}, N_k^{-1}).
\ees
\end{corollary}

\noindent
{\bf Proof. } Validity of the Lemma follows from the fact that each $\boeta_k$ is of the form \fr{sum_subgaus}
and is evaluated using a distinct set of components of vector $\bxi$.
\\

\vspace{4mm}

%%%%%%%%%%%%%%%%%%%%%%%%%%%%%%%%%%%%%%%%%%%%%%%%%%%%%%%%%%%%%%%%%%%%%%%%%%%%%%%%%%%%%%%%%%%%%%%%%%%%%%%%%

\begin{lemma} \label{lem:Gao2015} {\bf (Proposition 4.2 of Gao, Lu and Zhou (2015)) }
Let $\bte = \{ \te_{i,j}\} \in [0,1]^{n \times n}$. Let $\PP_{\te_{i,j}}$ denote the probability of
Bernoulli $(\te_{i,j})$ and the probability $\PP_{\bte}$ stands for the product measure 
$\PP_{\bte} = \otimes_{i,j} \PP_{\te_{i,j}}$.  Then, for any $\bte, \bte' \in [1/4, 3/4]^{n \times n}$ 
one has
$$
K(\PP_{\bte},\PP_{\bte'}) \leq 8 \, \sum_{i,j} (\te_{i,j} - \te'_{i,j})^2.
$$
\end{lemma}

%%%%%%%%%%%%%%%%%%%%%%%%%%%%%%%%%%%%%%%%%%%%%%%%%%%%%%%%%%%%%%%%%%%%%%%%%%%%%%%%%%%%%%%%%%%%%%%%%%%%%%%%%%%%

\end{document}